\documentstyle[bezier,12pt]{article}
\newcommand{\be}[1]{\begin{equation}\label{#1}}
\newcommand{\ee}{\end{equation}}
\newcommand{\rav}{\langle\,\Sigma\mid\rr\,\rangle}
\newcommand{\pres}[2]{\langle\,#1\mid #2\,\rangle}

\newcommand{\iv}{^{-1}}

\newcommand{\bs}{w}

\newcommand{\topp}{\mathop{\mbox{\bf top}}}
\newcommand{\bott}{\mathop{\mbox{\bf bot}}}

\newcommand{\rr}{{\cal R} }
\newcommand{\pp}{{\cal P} }
\newcommand{\qq}{{\cal Q} }
\newcommand{\sss}{{\cal S} }
\newcommand{\zz}{{\bf Z}}
\newcommand{\ve}{\varepsilon}
\newcommand{\eps}{\varepsilon}
\newcommand{\prf}{{\bf Proof. } }
\newcommand{\dg}{{\cal D}(\pp,\bs)}
\renewcommand{\wr}{\mathrel{\rm wr}}
\newcommand{\No}{{\char'362}}
\newcommand{\supp}{\mathop{\rm supp}}
\newcommand{\comp}{\mathop{\bf comp\,}}
\newcommand{\disto}{\mathop{\rm disto\,}}
\newcommand{\la}{\langle\,}
\newcommand{\ra}{\,\rangle}
\newcommand{\grp}[1]{\mathop{\rm gp\,}\la{#1}\ra}
\newcommand{\by}{\bar y}

\begin{document}

\title
{On subgroups of R.\,Thompson's group $F$\\
and other diagram groups}
\author{V.\,S.~Guba,\ \,M.\,V.\,Sapir}
\date{}

\maketitle

\begin{abstract}
In this paper, we continue our study of the class of diagram groups.
Simply speaking, a diagram is a labelled plane graph bounded by a pair of
paths (the top path and the bottom path). To multiply two diagrams, one
simply identifies the top path of one diagram with the bottom path of the
other diagram, and removes pairs of ``reducible" cells. Each diagram group
is determined by an alphabet $X$, containing all possible labels of edges,
a set of relations ${\cal R}=\{\,u_i=v_i\mid i=1,2,...\,\}$, containing all
possible labels of cells, and a word $w$ over $X$ -- the label of the top
and bottom paths of diagrams. Diagrams can be considered as 2-dimensional
words, and  diagram groups can be considered as 2-dimensional analogue of
free groups. In our previous paper, we showed that the class of diagram
groups contains many interesting groups including the famous R.\,Thompson
group $F$ (it corresponds to the simplest set of relations $\{\,x=x^2\,\}$),
closed under direct and free products and some other constructions. In this
paper we study mainly subgroups of diagram groups. We show that not every
subgroup of a diagram group is itself a diagram group (this answers a
question from the previous paper). We prove that every nilpotent subgroup of
a diagram group is abelian, every abelian subgroup is free, but even the
Thompson group contains solvable subgroups of any degree. We also study
distortion of subgroups in diagram groups, including the Thompson group.
It turnes out that centralizers of elements and abelian subgroups are always
undistorted, but the Thompson group contains distorted soluble subgroups.
\end{abstract}

\setlength{\unitlength}{0.001in}
\newtheorem{thm}{\quad Theorem}
\newtheorem{lm}[thm]{\quad Lemma}
\newtheorem{cy}[thm]{\quad Corollary}
\newtheorem{df}[thm]{\quad Definition}
\newtheorem{ex}[thm]{\quad Example}
\newtheorem{rk}[thm]{\quad Remark}
\newtheorem{prob}{\quad Problem}

\section*{Introduction}
\label{Intro}

$$
\ \
$$

This paper is devoted to further study of the so called {\em diagram groups}.
The definition of diagram groups was first given by Meakin and Sapir in 1995.
Their student Vesna Kilibarda obtained first results about diagram groups
in her thesis~\cite{KilDiss} (see also her paper~\cite{Kil}). Further results
about diagram groups have been obtained in our paper~\cite{GuSa97}. Here we survey
the main results of that paper (see~\cite{GuSa97} for details).
\vspace{1ex}

Diagram groups reflect certain important properties of semigroup
presentations. For instance we showed that three definitions of asphericity
given by Pride~\cite{Pr95b} are in fact equivalent and are equivalent to
the
triviality of all diagram groups over the presentation. One can say that
diagram groups measure the non-asphericity of semigroup presentations.

On the other hand, it turned out that the class of diagram groups is
interesting even if we forget about its connection with semigroup
presentations. These groups have nice algorithmic properties: the word
problem in every diagram group is solvable in time $O(n^{2+\ve})$ for every
$\ve>0$. This does not depend on whether the word problem for the
corresponding semigroup presentation is solvable or not. If the word problem
of this semigroup presentation is solvable then the conjugacy problem is
solvable in the corresponding diagram group.

If a group is representable by diagrams (i.e. it is a subgroup of a diagram
group) then one can use geometry of planar graphs to deduce certain
properties of the group. Diagrams can be viewed as ``2-dimensional words"
and in~\cite{GuSa97},  we developed a calculus called ``combinatorics on
diagrams", which is parallel to the well known combinatorics on words (see
Lothaire~\cite{Loth83}).

Geometry of diagrams allows one to consider many homomorphism from diagram
groups into the group of piecewise linear homeomorphisms of the real line.
Thus we have a connection between groups representable by diagrams and groups
representable by piecewise linear functions. This connection can be used in
both directions.

We showed that the class of diagram groups is wide. It contains the free
groups, free abelian groups, the R.\,Thompson group $F$ and its
generalizations found by Brown~\cite{Bro87}. This class is closed under
finite direct products, arbitrary free products and some other constructions.
Note that the Thompson group is the diagram group over the following simple
presentation $\pres{x}{x^2=x}$. In~\cite{GuSa97} we obtained several
previously unknown results about Thompson's group, essentially using its
representation as a diagram group.

\begin{itemize}
\item The conjugacy problem in $F$ is solvable.
\item The centralizer of each element of $F$ is a
finite direct product of groups each of which is either a copy of $F$
or an infinite cyclic group $\zz$.
\end{itemize}

Let us give a short summary of the content of this paper.
\vspace{1ex}

Section~\ref{Prelim} contains the list of the main concepts used in this
paper.
\vspace{0.5ex}

In Section~\ref{Constr}, we introduce the concept of diagram product of
groups. It is defined as the fundamental group of a certain 2-complex of
groups. Theorem~\ref{DiagProd}, the main result of this section, states that
the class of diagram groups is closed under diagram products. It turns out
that all ``products" considered before (the free product, the direct product,
etc.) are particular cases of the diagram product. Examples \ref{ex-dir} --
\ref{ex-DirPow}, \ref{ex-wr}, \ref{kh} show applications of
Theorem~\ref{DiagProd}. In partricular we prove that the class of diagram
groups is closed under countable direct powers (Theorem~\ref{dir-pow}),
wreath products with $\zz$ (Theorem~\ref{wr}), and a certain special
construction ${\cal O}(G,H)$ (Theorem~\ref{skew}), whose role will be clear
later.
\vspace{0.5ex}

In Section~\ref{NilAb}, we show that nilpotent subgroups of diagram groups
are abelian (Corollary~\ref{nilp}), and abelian subgroups are free abelian
(Theorem~\ref{FrAb}).
We finish this section with a description of sets of pairwise commuting
diagrams (Theorem~\ref{prws-cmmt}).
\vspace{0.5ex}

In Section~\ref{Sol}, we prove that the Thompson group $F$ contains subgroups
isomorphic to the restricted wreath product of two infinite cyclic groups and
soluble subgroups of arbitrary degree (this was first proved by Brin). It
turns out that for any subgroup of piecewise linear functions (including $F$)
there exists a dichotomy: either it contains $\zz\wr\zz$ or it is abelian
(Theorem~\ref{PLIZwrZ}). This implies, in particular, that a non-abelian
subgroup of the group of piecewise linear functions cannot be a one-relator
group (Corollary~\ref{OneRelPL}). This result strengthens the well known fact
that the group of piecewise linear functions does not contain free
non-abelian subgroups.

In Section~\ref{Sol}, we also give necessary and sufficient conditions for a
diagram group to contain a copy of $\zz\wr\zz$ as a subgroup
(Theorem~\ref{3equiv}). We study the question when a diagram group over some
semigroup presentation $\pp$ contains a copy of the Thompson group $F$. We
prove that if the semigroup given by $\pp$ contains an idempotent then a
diagram group over this presentation contains a copy of $F$.
(Theorem~\ref{FSub}). The interesting question of whether the converse
statement holds is open.
\vspace{0.5ex}

In Section~\ref{SubConj}, we present a counterexample to the Subgroup
Conjecture. This conjecture stated that every subgroup of a diagram
group is a diagram group itself. It was motivated by the similarity between
diagram groups and free groups. At first we thought that the conjecture is
easy to disprove and the derived subgroup $F'$ of $F$ is a counterexample.
But it turned out that $F'$ is a diagram group (Theorem~\ref{CommF}). This
solves several problems from~\cite{GuSa97}. We asked whether a diagram group
can coincide with its derived subgroup and whether every diagram group has an
LOG-presentation. Corollary~\ref{simp} gives a positive answer to the first
question and a negative answer to the second question.

But the main result of this Section is Theorem~\ref{dispr} which shows
that the one-relator group $\pres{x,y}{xy^2x=yx^2y}$ is not a diagram group
but is isomorphic to a subgroup of a diagram group. In the proof, we use the
construction ${\cal O}(G,H)$ from Section~\ref{Constr}. This gives a
counterexample to the Subgroup Conjecture. Nevertheless, we think that in
many partricular cases this conjecture is true and we pose several open
question in this regard.

At the end of this section we study the following series of groups
$$
G_n=\pres{x_1,\ldots,x_n}
{[x_1,x_2]=[x_2,x_3]=\cdots=[x_{n-1},x_n]=[x_n,x_1]=1}.
$$
For $n\le4$ these groups are diagram groups; we prove (Theorem~\ref{PartCom}),
that for odd $n\ge5$ this is not so. It is not known whether these groups are
representable by diagrams. If so, this will give us new counterexamples to
the Subgroup Conjecture.
\vspace{0.5ex}

The last Section~\ref{Disto} is devoted to the distortrion of subgroups in
diagram groups. In a recent paper~\cite{Bur} Burillo prroved that for every
natural $n$ the Thompson group $F$ contains subgroups isomorphic to
$F\times\zz^n$ and quasi-isometrically (without distortrion) embedded into
$F$. A similar fact is true for an embedding of $F\times F$. We prove
(Theorem~\ref{q-i-C}) that every centralizer of an element in $F$ is embedded
into $F$ without distortion. Centralizers of elements of $F$ can be arbitrary
finite direct products of copies of $F$ and copies of $\zz$. Burillo also
proved that every cyclic subgroup of $F$ is embedded without distortion
(this fact is also an immediate corollary of Lemma~15.29 in~\cite{GuSa97}).
We prove (Theorem~\ref{Znqi}) that not only cyclic but arbitrary finitely
generated abelian subgroups of {\bf any} diagram groups are undistorted.
Finally we found solvable subgroups of $F$ which are distorted.
Theorem~\ref{Disto-nd} shows that for every natural $d\ge2$ there exists a
finitely generated solvable subgroup $K_d$ in $F$ such that its distortion
function is at least $n^d$.
\vspace{2ex}

{\bf Acknowledgements.} The authors thank M.~Brin and S.~Pride for helpful
discussions of the results of this paper.

\section{Preliminaries}
\label{Prelim}

For an alphabet $\Sigma$ let $\Sigma^+$ denote the free semigroup over
$\Sigma$, and let $\Sigma^*$ denote the free monoid.  Elements of the free
monoid are called {\em words}. The identity element, i.e. the empty
word, is denoted by $1$.

Let $\pp=\rav$ be a presentation of a semigroup where $\Sigma$ is an
alphabet, $\rr$ is a set of pairs of non-empty words over $\Sigma$.
The semigroup $S$ given by $\pp$ is the factor-semigroup $\Sigma^+/\!\sim$
where $\sim$ is the smallest congruence on $\Sigma^+$ containing $\rr$.
Elements of $\Sigma$ are called {\em generators}, pairs $(u,v)\in\rr$
written also as $u=v$ are called {\em defining relations}. Left and right
parts of defining relations are called {\em defining words}. We shall assume
that all presentations are anti-symmetric that is if $u=v\in\rr$ then
$v=u\not\in \rr$. In particular $\rr$ does not contain relations $u=u$.

With any semigroup presentation $\pp$, we associate the following graph
$\Gamma(\pp)$.  The vertices are all words in $\Sigma^+$. Edges are the
elements of $\Sigma^*\times\rr^{\pm1}\times\Sigma^*$. We shall denote edges
by $(x,u\to v,y)$, where $x,y\in\Sigma^*$, and either $(u,v)\in\rr$ or
$(v,u)\in\rr$. If $e=(x,u\to v,y)$ then the {\em inverse edge} is defined by
$e\iv=(x,v\to u,y)$. The {\em initial vertex\/} of $e$ is $\iota(e)=xuy$ and
the {\em terminal  vertex\/} is $\tau(e)=xvy$. Thus vertices of this graph
are words and edges are elementary transformations of words (i.e.
substitutions of defining words by their pairs). The graph $\Gamma(\pp)$
describes all derivations over $\pp$: two non-empty words $w_1$, $w_2$ are
{\em equal modulo\/} $\pp$ if and only if there exists a path in the graph
connecting $w_1$ and $w_2$. This path is called a {\em derivation} of $w_2$
from $w_1$.

With every derivation over $\pp$, one can associate a geometric object,
a {\em semigroup diagram\/} over $\pp$.  Semigroup diagrams were first
introduced by E.\,V.\,Kashintsev~\cite{Kash70} and then rediscovered by
Remmers~\cite{Rem} and others(see~\cite{Higg92,Sta87}). We do not give an
exact definition here (see also~\cite{GuSa97}), the definition will be clear
from the following example.

\begin{ex}
\label{ex1}
{\rm
Let $\pp=\langle\,a,b,c\mid abc=ba,bca=cb,cab=ac\,\rangle$. Consider the
following derivation over $\pp$:
$$
(1,ac\to cab,cb^2ca)(c,abc\to ba,b^2ca)(cbab,bca\to cb,1)(cb,abc\to ba,b).
$$
The corresponding diagram $\Delta$ over $\pp$ is the following:

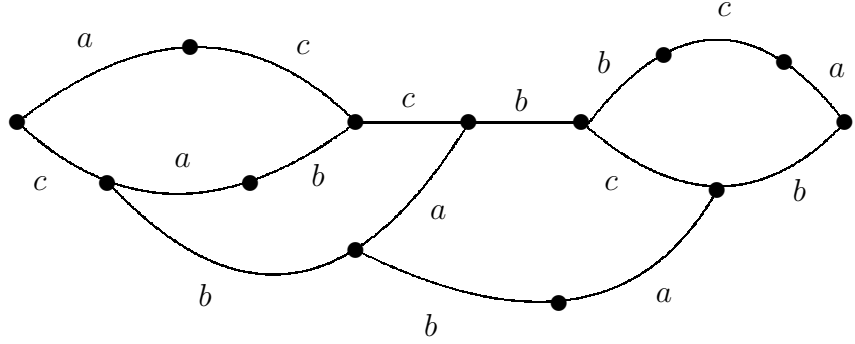
\begin{center}
\unitlength=1mm
\linethickness{0.4pt}
\begin{picture}(112.00,45.00)
\put(1.00,29.00){\circle*{2.00}}
\put(46.00,29.00){\circle*{2.00}}
\put(61.00,29.00){\circle*{2.00}}
\put(76.00,29.00){\circle*{2.00}}
\put(111.00,29.00){\circle*{2.00}}
\put(46.00,29.00){\line(1,0){15.00}}
\put(61.00,29.00){\line(1,0){15.00}}
\bezier{240}(1.00,29.00)(26.00,49.00)(46.00,29.00)
\bezier{236}(1.00,29.00)(21.00,10.00)(46.00,29.00)
\bezier{224}(77.00,29.00)(94.00,51.00)(111.00,29.00)
\bezier{196}(76.00,29.00)(95.00,12.00)(111.00,29.00)
\put(24.00,39.00){\circle*{2.00}}
\put(87.00,38.00){\circle*{2.00}}
\put(103.00,37.00){\circle*{2.00}}
\put(94.00,20.00){\circle*{2.00}}
\put(13.00,21.00){\circle*{2.00}}
\put(32.00,21.00){\circle*{2.00}}
\bezier{320}(13.00,21.00)(39.00,-7.00)(61.00,29.00)
\put(46.00,12.00){\circle*{2.00}}
\bezier{268}(46.00,12.00)(80.00,-5.00)(94.00,20.00)
\put(73.00,5.00){\circle*{2.00}}
\put(10.00,40.00){\makebox(0,0)[cc]{$a$}}
\put(39.00,39.00){\makebox(0,0)[cc]{$c$}}
\put(53.00,32.00){\makebox(0,0)[cc]{$c$}}
\put(68.00,32.00){\makebox(0,0)[cc]{$b$}}
\put(79.00,37.00){\makebox(0,0)[cc]{$b$}}
\put(95.00,44.00){\makebox(0,0)[cc]{$c$}}
\put(110.00,36.00){\makebox(0,0)[cc]{$a$}}
\put(4.00,21.00){\makebox(0,0)[cc]{$c$}}
\put(23.00,24.00){\makebox(0,0)[cc]{$a$}}
\put(41.00,22.00){\makebox(0,0)[cc]{$b$}}
\put(26.00,6.00){\makebox(0,0)[cc]{$b$}}
\put(57.00,17.00){\makebox(0,0)[cc]{$a$}}
\put(80.00,21.00){\makebox(0,0)[cc]{$c$}}
\put(105.00,20.00){\makebox(0,0)[cc]{$b$}}
\put(56.00,2.00){\makebox(0,0)[cc]{$b$}}
\put(87.00,6.00){\makebox(0,0)[cc]{$a$}}
\end{picture}
\end{center}
}
\end{ex}

Let us introduce the terminology associated with diagrams. Every diagram
over $\pp$ is a planar graph. It has {\em vertices}, {\em edges} and
{\em cells}. In Example~\ref{ex1}, the diagram $\Delta$ has 13 vertices,
16 edges and 4 cells. The number of cells is equal to the length of the
corresponding derivation. Each positive edge has a {\em label} from
$\Sigma$, positive edges are oriented from left to right. The label of an
edge $e$ is denoted by $\varphi(e)$. We shall consider only positive paths
in $\Delta$, that is paths consisting of positive edges. For every path $p$
in a diagram $\Delta$, its {\em label} $\varphi(p)$ is the word read on the
path. Any diagram $\Delta$ has the {\em initial verrtex\/} $\iota(\Delta)$
and the {\em terminal vertex\/} $\tau(\Delta)$, the {\em top path\/}
$\topp(\Delta)$ and the {\em bottom path\/} $\bott(\Delta)$ connecting the
initial and terminal vertices. The diagram $\Delta$ lies between its top and
bottom paths. This notation is illustrated by the following example.

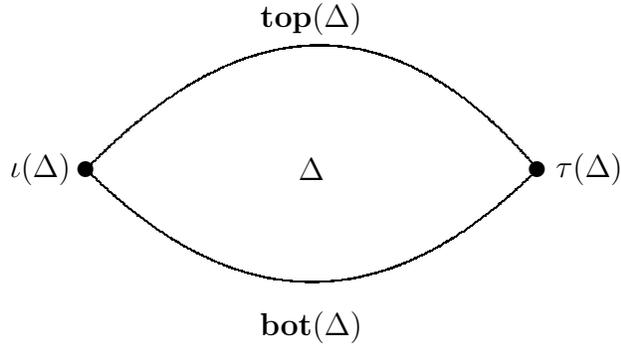
\begin{center}
\unitlength=1mm
\special{em:linewidth 0.4pt}
\linethickness{0.4pt}
\begin{picture}(72.00,50.00)
\put(7.00,22.00){\circle*{2.00}}
\put(67.00,22.00){\circle*{2.00}}
\bezier{356}(7.00,22.00)(39.00,55.00)(67.00,22.00)
\bezier{340}(7.00,22.00)(37.00,-8.00)(67.00,22.00)
\put(37.00,22.00){\makebox(0,0)[cc]{$\Delta$}}
\put(1.00,22.00){\makebox(0,0)[cc]{$\iota(\Delta)$}}
\put(74.00,22.00){\makebox(0,0)[cc]{$\tau(\Delta)$}}
\put(37.00,42.00){\makebox(0,0)[cc]{$\topp(\Delta)$}}
\put(37.00,1.00){\makebox(0,0)[cc]{$\bott(\Delta)$}}
\end{picture}
\end{center}

Notice that paths $\topp(\Delta)$ and $\bott(\Delta)$ can have common edges. Every cell $\pi$ of a diagram is a diagram itself, so we can define the notation
$\iota(\pi)$, $\tau(\pi)$, $\topp(\pi)$, $\bott(\pi)$ and the corresponding
concepts. If words $u$ and $v$ are labels of the top and the bottom paths of a
cell $\pi$ then either $u=v$ or $v=u$ is a defining relation (that is it
belongs to $\rr$). In this case we call $\pi$ a $(u,v)$-cell.

For every non-empty word $w$, there exists a {\em trivial diagram} $\ve(w)$
without cells whose top and bottom paths coincide and have label $w$.

We do not distinguish isotopic diagrams. The notation $\Delta_1\equiv\Delta_2$
means that $\Delta_1$ and $\Delta_2$ are isotopic.

If the label of $\topp(\Delta)$ is $w_1$ and the label of $\bott(\Delta)$ is
$w_2$ then $\Delta$ is called a $(w_1,w_2)$-diagram. Let $w_1$, $w_2$, $w_3$
 be any three vertices of the graph $\Gamma(\pp)$ and let $p_i$ ($i=1,2$) be
paths in the graph $\Gamma(\pp)$ from $w_i$ to $w_{i+1}$. By $\Delta_i$, we
denote the diagram corresponding to the path $p_i$ ($i=1,2$). It is easy to
see that the product of paths $p_1$ and $p_2$ in the graph corresponds to the
diagram $\Delta$ obtained from $\Delta_1$ and $\Delta_2$ by identifying the
bottom path of $\Delta_1$ and the top path of $\Delta_2$. The resulting
diagram $\Delta$ will be called the {\em composition\/} of diagrams
$\Delta_1$ and $\Delta_2$, we denote it by $\Delta_1\circ\Delta_2$. Thus
$\circ$ is a partial operation on the set of all diagrams over $\pp$.
The composition of a $(w_1,w_2)$-diagram and a $(w_2,w_3)$-diagram is a
$(w_1,w_3)$-diagram. For every word $w\in\Sigma^+$, the set of all
$(w,w)$-diagrams over $\pp$ is a semigroup with respect to the operation
$\circ$. Diagrams of this form will be called {\em spherical\/} diagrams
with {\em base} $w$. This semigroup has the identity element $\ve(w)$.
We define also another associative operation on the set of all diagrams over
$\pp$. Namely the {\em sum\/} $\Delta_1+\Delta_2$ of diagrams $\Delta_1$ and
$\Delta_2$ is the diagram obtained by identifying $\tau(\Delta_1)$ and
$\iota(\Delta_2)$. These two operations are illustrated by the following
figure:

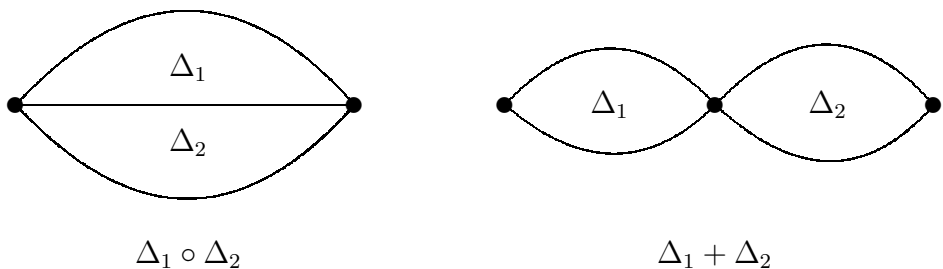
\begin{center} 
\unitlength=1mm
\special{em:linewidth 0.4pt}
\linethickness{0.4pt}
\begin{picture}(124.41,55.00)
\put(1.00,30.00){\circle*{2.00}}
\put(46.00,30.00){\circle*{2.00}}
\put(1.00,30.00){\line(1,0){45.00}}
\bezier{320}(1.00,30.00)(24.00,55.00)(46.00,30.00)
\bezier{332}(1.00,30.00)(24.00,5.00)(46.00,30.00)
\put(24.00,35.00){\makebox(0,0)[cc]{$\Delta_1$}}
\put(24.00,25.00){\makebox(0,0)[cc]{$\Delta_2$}}
\put(24.00,10.00){\makebox(0,0)[cc]{$\Delta_1\circ\Delta_2$}}
\put(66.00,30.00){\circle*{2.00}}
\put(94.00,30.00){\circle*{2.00}}
\put(123.00,30.00){\circle*{2.00}}
\bezier{164}(66.00,30.00)(80.00,45.00)(94.00,30.00)
\bezier{152}(66.00,30.00)(81.00,17.00)(94.00,30.00)
\bezier{172}(94.00,30.00)(109.00,46.00)(123.00,30.00)
\bezier{168}(94.00,30.00)(110.00,15.00)(123.00,30.00)
\put(80.00,30.00){\makebox(0,0)[cc]{$\Delta_1$}}
\put(109.00,30.00){\makebox(0,0)[cc]{$\Delta_2$}}
\put(94.00,10.00){\makebox(0,0)[cc]{$\Delta_1+\Delta_2$}}
\end{picture}
\end{center}

Suppose that a diagram $\Delta$ contains a $(u,v)$-cell and a $(v,u)$-cell
such that the top path of the first cell is the bottom path of the second
cell. Then we say that these two cells form a {\em dipole}. In this case we
can remove these two cells by first removing their common path, and then
identifying the bottom path of the first cell with the top path of the
second cell. A diagram is called {\em reduced} if it does not contain
dipoles. One can get a reduced diagram from any diagram by removing dipoles.
Kilibarda~\cite{KilDiss} proved that every diagram has a unique reduced form.
We call two diagrams $\Delta_1$ and $\Delta_2$ {\em equivalent}, written as
$\Delta_1\cong\Delta_2$, if their reduced forms are the same. It is easy to
see that if $\Delta_1\cong \Delta_2$, $\Delta_3\cong\Delta_4$ then
$\Delta_1\circ\Delta_3\cong \Delta_2\circ\Delta_4$ and
$\Delta_1+\Delta_3\cong\Delta_2+\Delta_4$.

Therefore on the set $\dg$ of all equivalence classes of $(\bs,\bs)$-diagrams
one can define a product, by setting
$[\Delta_1]\cdot[\Delta_2]=[\Delta_1\circ\Delta_2]$, where square brackets
denote equivalence classes.

The product of a $(w_1,w_2)$-diagram $\Delta$ and the $(w_2,w_1)$-diagram
$\Delta'$, which is a mirror image of $\Delta$ is obviously equivalent to
the trivial diagram $\ve(w_1)$. The diagram $\Delta'$ will be denoted by
$\Delta\iv$.

For simplicity we shall call equivalent diagrams {\em equal}, use ``=" instead
of ``$\cong$", and drop square brackets and the multiplication sign. So for
every $\Delta\in\dg$ $\Delta\Delta\iv=\ve(\bs)$. As a result $\dg$ turns out
to be a group which is called the {\em diagram group} over the semigroup
presentation $\pp$ with base $\bs$. Since every equivalence class contains
a unique reduced diagram, one can assume that $\dg$  consists of reduced
diagrams with the natural multiplication ($\Delta_1\Delta_2$ is the reduced form
of $\Delta_1\circ\Delta_2$).

In what follows, the term {\em diagram group} means a diagram group over
some presentation with some base.
\vspace{2ex}

We shall use the standard notation for conjugation in groups: $a^b=b\iv ab$, and
for the commutator: $[a,b]=a\iv a^b=a\iv b\iv ab$. If $A$ and $B$ are
subgroups of a group $G$ then $[A,B]$ denotes the subgroup generated by all
commutators $[a,b]$ where $a\in A$, $b\in B$.

Now let us shortly describe some results about diagram groups obtained
earlier.

The diagram group corresponding to the presentation
$\pp=\langle\,x\mid xx=x\,\rangle$ with base $x$ is the famous R.\,Thompson's
group $F$, which has the following presentation:
$$
\langle\,x_0,x_1,\ldots\mid x_j^{x_i}=x_{j+1}\ (j>i)\,\rangle.
$$
(see~\cite[Example 6.4]{GuSa97}.

This group has several interesting propertries and is studied by
mathematicians working in different areas of mathematics ($\lambda$-calculus,
functional analysis, homological algebra, homotopy theory, group theory).
It was discovered by R.\,Thompson in 1965, and was rediscovered later by
other authors. \cite{CFP} presents a survey of results about $F$.
Since $F$ is one of the most important diagram groups, and since we are
going to present some new results about it in this paper, let us recall some
known properties of this group. These propertries can be found in~\cite{CFP},
\cite{GuSa97} and~\cite{Gu98}.

\begin{enumerate}
\item
The group $F$ is isomorphic to the group of all increasing continuous
piecewise linear maps of the interval $[0,1]$ onto itself such that the
singularities occur at finitely many dyadic points (points of the form
$m/2^n$) and all slopes are powers of 2. The group operation is the
composition of functions (we shall write function symbols to the right of
the argument).
\item In the previous paragraph, one can replace the interval $[0,1]$
by $[0,+\infty]$, adding the assumption that the slop on $+\infty$ is $1$.
The resulting group is also isomorphic to $F$.
\item
$F$ does not satisfy any non-trivial identity.
\item
$F$ does not contain any free non-abelian subgroups.
Every subgroup of $F$ either is abelian or contains an infinite direct
power of $\zz$.
\item
$F$ is finitely presented, it has a presentation with two generators and two
defining relations. The word problem and the conjugacy problem are solvable
in $F$. It has a polynomial isoperimetric function~\cite{Gu98}.
\end{enumerate}
\vspace{1ex}

There exists a clear connection between representation of elements of $F$ by
diagrams and normal form of elements in $F$. Recall~\cite{CFP} that every
element in $F$ is uniquely representable in the following form:
\be{NormForm}
x_{i_1}^{s_1}\ldots x_{i_m}^{s_m}x_{j_n}^{-t_n}\ldots x_{j_1}^{-t_1},
\ee
where $i_1\le\cdots\le i_m\ne j_n\ge\cdots\ge j_1$;
$s_1,\ldots,s_m,t_1,\ldots t_n\ge0$, and if $x_i$ and $x_i\iv$ occur in
(\ref{NormForm}) for some $i\ge0$ then either $x_{i+1}$ or $x_{i+1}\iv$
also occurs in~(\ref{NormForm}).
This form is called the {\em normal form} of elements in $F$. (Note that
in ~\cite{GuSa97a}, we constructed another normal form for elements of $F$,
our normal forms are locally testable.)

Let us show how given an $(x,x)$-diagram over $\pp=\pres{x}{xx=x}$ one can
get the normal form of the element represented by this diagram. We simply
describe the procedure providing no proofs. Details can be deduced
from~\cite[Example 6.4]{GuSa97}.

\begin{ex}
\label{TGNF}
{\rm
Every diagram $\Delta$ over $\pp$ can be divided by its longest positive
path from its intitial vertex to its terminal vertex into two parts,
{\em positive} and {\em negative}, denoted by $\Delta^+$ and $\Delta^-$,
respectively. So $\Delta=\Delta^+\circ\Delta^-$. It is easy to prove by
induction on the number of cells that all cells in $\Delta^+$ are
$(x,x^2)$-cells, all cells in $\Delta^-$ are $(x^2,x)$-cells. This implies
that the numbers of cells in $\Delta^+$ and in $\Delta^-$ are the same.
Denote this number by $k$. Let us number the cells of $\Delta^+$ by
numbers from $1$ to $k$ by taking every time the ``rightmost" cell,
that is, the cell which is to the right of any other cell attached to the
bottom path of the diagram formed by the previous cells. The first cell is
attached to the top path of $\Delta^+$ ($=\topp(\Delta)$). The $i$th cell in
this sequence of cells corresponds to an edge of the graph $\Gamma(\pp)$,
which has the form $(x^{\ell_i},x\to x^2,x^{r_i})$, where $\ell_i$ ($r_i$)
is the length of the path from the initial (resp. terminal) vertex of the
diagram (resp. the cell) to the initial (resp. terminal) vertex of the cell
(resp. the diagram), and contained in the bottom path of the diagram formed by the first $i-1$ cells. If $\ell_i=0$ then we label this cell by 1. If
$\ell_i\ne0$ then we label this cell by the element $x_{r_i}$ of $F$.
Multiplying the labells of all cells, we get the ``positive" part of the
normal form. For example, the diagram on the next picture

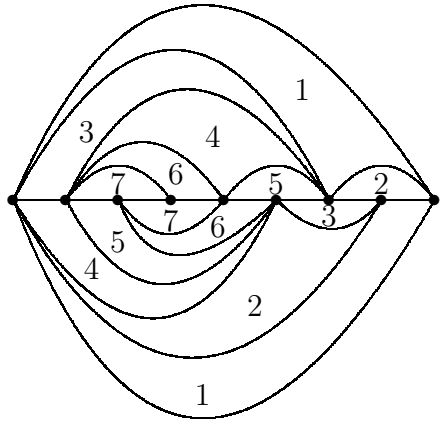
\begin{center} 
\unitlength=0.7mm
\special{em:linewidth 0.4pt}
\linethickness{0.4pt}
\begin{picture}(83.00,90.00)
\put(2.00,43.00){\circle*{2.00}}
\put(12.00,43.00){\circle*{2.00}}
\put(22.00,43.00){\circle*{2.00}}
\put(32.00,43.00){\circle*{2.00}}
\put(42.00,43.00){\circle*{2.00}}
\put(52.00,43.00){\circle*{2.00}}
\put(62.00,43.00){\circle*{2.00}}
\put(72.00,43.00){\circle*{2.00}}
\put(82.00,43.00){\circle*{2.00}}
\put(2.00,43.00){\line(1,0){10.00}}
\put(12.00,43.00){\line(1,0){10.00}}
\put(22.00,43.00){\line(1,0){10.00}}
\put(32.00,43.00){\line(1,0){10.00}}
\put(42.00,43.00){\line(1,0){10.00}}
\put(52.00,43.00){\line(1,0){10.00}}
\put(62.00,43.00){\line(1,0){10.00}}
\put(72.00,43.00){\line(1,0){10.00}}
\bezier{132}(62.00,43.00)(72.00,56.00)(82.00,43.00)
\bezier{132}(62.00,43.00)(52.00,56.00)(42.00,43.00)
\bezier{132}(12.00,43.00)(22.00,56.00)(32.00,43.00)
\bezier{212}(12.00,43.00)(26.00,65.00)(42.00,43.00)
\bezier{392}(12.00,43.00)(33.00,85.00)(62.00,43.00)
\bezier{516}(2.00,43.00)(33.00,100.00)(62.00,43.00)
\bezier{676}(2.00,43.00)(34.00,117.00)(82.00,43.00)
\bezier{132}(22.00,43.00)(31.00,30.00)(42.00,43.00)
\bezier{212}(22.00,43.00)(30.00,22.00)(52.00,43.00)
\bezier{120}(52.00,43.00)(63.00,32.00)(72.00,43.00)
\bezier{304}(12.00,43.00)(29.00,11.00)(52.00,43.00)
\bezier{412}(2.00,43.00)(30.00,-2.00)(52.00,43.00)
\bezier{556}(2.00,43.00)(35.00,-17.00)(72.00,43.00)
\bezier{740}(2.00,43.00)(34.00,-40.00)(82.00,43.00)
\put(57.00,64.00){\makebox(0,0)[cc]{1}}
\put(72.00,46.00){\makebox(0,0)[cc]{2}}
\put(16.00,56.00){\makebox(0,0)[cc]{3}}
\put(40.00,55.00){\makebox(0,0)[cc]{4}}
\put(52.00,46.00){\makebox(0,0)[cc]{5}}
\put(33.00,48.00){\makebox(0,0)[cc]{6}}
\put(22.00,46.00){\makebox(0,0)[cc]{7}}
\put(38.00,6.00){\makebox(0,0)[cc]{1}}
\put(48.00,23.00){\makebox(0,0)[cc]{2}}
\put(62.00,40.00){\makebox(0,0)[cc]{3}}
\put(17.00,30.00){\makebox(0,0)[cc]{4}}
\put(22.00,35.00){\makebox(0,0)[cc]{5}}
\put(41.00,38.00){\makebox(0,0)[cc]{6}}
\put(32.00,39.00){\makebox(0,0)[cc]{7}}
\end{picture}
\end{center}

\noindent
the positive part is equal to $x_0x_2^2x_4x_5$ (cells 1 and 3 were labelled
by the identity element).

In order to find the ``negative" part of the normal form, consider
$(\Delta^-)\iv$, number its cells as above and label them as above. In our
example, we get the word $x_1x_3^2x_4$ (cells 1, 2, 4 are labelled by 1).
Thus the ``negative" part of the normal form is $(x_1x_3^2x_4)\iv$, and it
remains to multiply the positive and negative parts. In our example,
the normal form is
$x_0x_2^2x_4x_5x_4\iv x_3^{-2}x_1\iv$.
}
\end{ex}

Diagrams presented below are generators $x_0$, $x_1$ of the group $F$. They
generate the whole $F$.

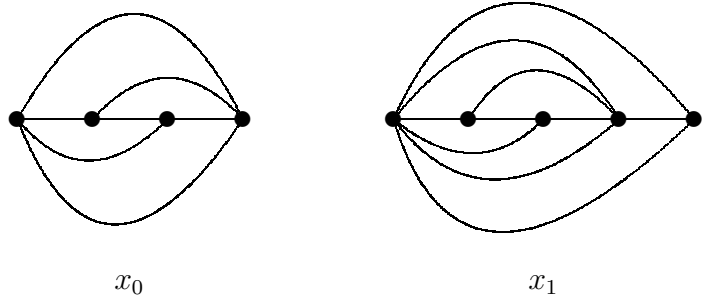
\begin{center} 
\unitlength=1mm
\special{em:linewidth 0.4pt}
\linethickness{0.4pt}
\begin{picture}(94.00,50.00)
\put(3.00,24.00){\circle*{2.00}}
\put(13.00,24.00){\circle*{2.00}}
\put(23.00,24.00){\circle*{2.00}}
\put(33.00,24.00){\circle*{2.00}}
\put(53.00,24.00){\circle*{2.00}}
\put(63.00,24.00){\circle*{2.00}}
\put(73.00,24.00){\circle*{2.00}}
\put(83.00,24.00){\circle*{2.00}}
\put(93.00,24.00){\circle*{2.00}}
\put(3.00,24.00){\line(1,0){10.00}}
\put(13.00,24.00){\line(1,0){10.00}}
\put(23.00,24.00){\line(1,0){10.00}}
\put(53.00,24.00){\line(1,0){10.00}}
\put(63.00,24.00){\line(1,0){10.00}}
\put(73.00,24.00){\line(1,0){10.00}}
\put(83.00,24.00){\line(1,0){10.00}}
\bezier{120}(13.00,24.00)(23.00,35.00)(33.00,24.00)
\bezier{120}(3.00,24.00)(12.00,13.00)(23.00,24.00)
\bezier{256}(3.00,24.00)(19.00,52.00)(33.00,24.00)
\bezier{256}(3.00,24.00)(14.00,-4.00)(33.00,24.00)
\bezier{132}(63.00,24.00)(71.00,37.00)(83.00,24.00)
\bezier{108}(53.00,24.00)(64.00,15.00)(73.00,24.00)
\bezier{208}(53.00,24.00)(70.00,45.00)(83.00,24.00)
\bezier{176}(53.00,24.00)(65.00,8.00)(83.00,24.00)
\bezier{296}(53.00,24.00)(67.00,55.00)(93.00,24.00)
\bezier{296}(53.00,24.00)(62.00,-6.00)(93.00,24.00)
\put(18.00,2.00){\makebox(0,0)[cc]{$x_0$}}
\put(73.00,2.00){\makebox(0,0)[cc]{$x_1$}}
\end{picture}
\end{center}

One can ask several natural general questions about diagram groups.
Which groups are diagram groups? Which groups are representable by diagrams
(are subgroups of diagram groups)? Do these two classes coinside? How to
compute a diagram group over a given presentation $\pp$ and with a given
base?

There exists a well developed technology for computing diagram groups.
The starting point for computing diagram groups is Kilibarda's theorem
about fundamental groups of Squier complexes. In order to formulate this
important result, let us define the structure of a 2-complex on the graph
$\Gamma(\pp)$ for any presentation $\pp$.

First notice that although for every path in $\Gamma(\pp)$ there exists
a unique diagram associated with this path, the same diagram can be associated
with many paths. Consider the following typical case:

Let $\pp=\rav$ where $\rr$ contains two defining relations $\ell_i=r_i$
($i=1,2$), and let $u$, $v$, $z$ be arbitrary words in $\Sigma^*$. Consider
the following paths in $\Gamma(\pp)$:
\be{path1}
(u,\ell_1\to r_1,z\ell_2v)(ur_1z,\ell_2\to r_2,v),
\ee
\be{path2}
(u\ell_1z,\ell_2\to r_2,v)(u,\ell_1\to r_1,zr_2v).
\ee
It is easy to see that $(u\ell_1z\ell_2v,ur_1zr_2v)$-diagrams
corresponding to these paths are equal. This diagram is shown on the
following picture:

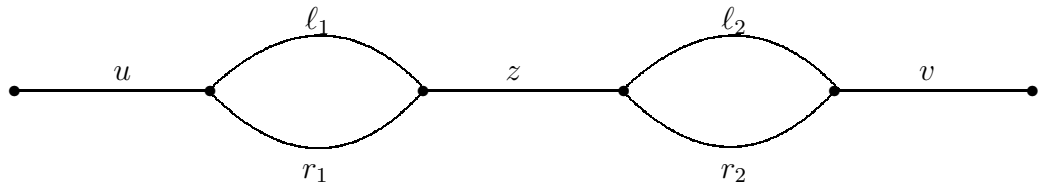
\begin{center} 
\unitlength=1.00mm
\special{em:linewidth 0.4pt}
\linethickness{0.4pt}
\begin{picture}(137.66,29.33)
\put(1.67,14.67){\line(1,0){25.67}}
\bezier{164}(27.33,14.67)(42.67,29.33)(56.33,14.67)
\put(56.33,14.67){\line(1,0){25.67}}
\bezier{164}(82.33,14.67)(97.67,29.33)(111.33,14.67)
\put(111.00,14.67){\line(1,0){25.67}}
\bezier{168}(27.33,14.67)(42.33,-0.67)(56.00,14.67)
\bezier{164}(82.33,14.67)(97.67,29.33)(111.33,14.67)
\bezier{164}(82.33,14.67)(97.00,-0.33)(110.67,14.67)
\put(1.67,14.67){\circle*{1.33}}
\put(27.67,14.67){\circle*{1.33}}
\put(56.00,14.67){\circle*{1.33}}
\put(82.67,14.67){\circle*{1.33}}
\put(110.67,14.67){\circle*{1.33}}
\put(137.00,14.67){\circle*{1.33}}
\put(16.00,17.00){\makebox(0,0)[cc]{$u$}}
\put(68.00,17.00){\makebox(0,0)[cc]{$z$}}
\put(123.00,17.00){\makebox(0,0)[cc]{$v$}}
\put(42.00,24.00){\makebox(0,0)[cc]{$\ell_1$}}
\put(97.33,24.00){\makebox(0,0)[cc]{$\ell_2$}}
\put(97.33,3.67){\makebox(0,0)[cc]{$r_2$}}
\put(41.67,3.67){\makebox(0,0)[cc]{$r_1$}}
\end{picture}
\end{center}

This situation hints to a homotopy relation on the set of paths in the graph
$\Gamma(\pp)$: paths~(\ref{path1}) and~(\ref{path2}) should be called
homotopic. In order to define the homotopy relation we need the structure
of a 2-complex on $\Gamma(\pp)$. For every 5-tuple
$(u,\ell_1=r_1,z,\ell_2=r_2,v)$, where $u,v,z\in\Sigma^*$, $(\ell_1=r_1),
(\ell_2=r_2)\in\rr$ we have a $2$-cell whose defining path is $p_1p_2^{-1}$
where $p_1$, $p_2$ are the paths~(\ref{path1}) and~(\ref{path2}),
respectively. The resulting $2$-complex is called the {\em Squier complex\/}
of the semigroup presentation $\pp$. It is denoted by ${\cal K}(\pp)$. It
was implicitely defined by Squier in~\cite{Sq94}. The same complex was
independently constructed by Kilibarda~\cite{KilDiss,Kil} and
Pride~\cite{Pr95a}. The important role of this complex is justified by the
fact that equal diagrams over $\pp$ correspond to homotopic paths in
${\cal K}(\pp)$. The following Kilibarda's theorem~\cite{KilDiss, Kil} plays
an important role in this paper: {\em The diagram group\/} $\dg$ {\em is
isomorphic to the fundamental group\/} $\pi_1({\cal K},\bs)$ {\em of the
Squier complex} ${\cal K}={\cal K}(\pp)$.
\vspace{2ex}

\section{Diagram Product of Groups}
\label{Constr}

In~\cite{GuSa97}, we considered several group-theoretical operations such
that the class of diagram groups is closed under them. These operations
were: finite direct products (result due to Kilibarda~\cite{KilDiss}),
any free products, and also some special operation $\bullet$ which we used
for constructing an example of a diagram group that was finitely generated
but not finitely presented. In this Section we introduce a quite general
operation on groups, the {\em diagram product}. We show that the class of
diagram groups is closed under this operation. All the above listed
constructions are partial cases of this new operation. We also consider some
concrete applications of this construction. They will be essentially used in
the later Sections.
\vspace{1ex}

Let us recall the definition of a graph. A {\em graph\/} (in the sense of
Serre~\cite{Se80}) is an ordered tuple $\Gamma=\la V,E,\iv,\iota,\tau\ra$
where $V$, $E$ are disjoint sets, $\iv$ is an involution on $E$, $\iota$,
$\tau$ are mappings from $E$ to $V$. The following axioms hold:

\begin{itemize}
\item $e\iv\ne e$ for any $e\in E$;
\item $\iota(e\iv)=\tau(e)$, $\tau(e\iv)=\iota(e)$.
\end{itemize}

Elements of the sets $V$ and $E$ are called {\em vertices\/} and {\em edges}
of the graph respectively. If $e\in E$, then $\iota(e)$ is called the
{\em initial vertex\/} of the edge $e$, and $\tau(e)$ is called the
{\em terminal vertex\/} of the edge $e$.

A {\em path} on the graph $\Gamma$ is either a vertex, or a nonempty
sequence of edges $e_1$, $e_2$, \dots, $e_n$ such that
$\tau(e_i)=\iota(e_{i+1})$ for each $i=1,\ldots,n-1$. Usually a path is
written in the form $p=e_1e_2\ldots e_n$. If a path $p$ consists of a
vertex $v$, then it is called an {\em empty path} and we denote it by $1_v$.
If $p=e_1e_2\ldots e_n$ is a path, then the {\em inverse path\/} $p^{-1}$ is
the path $e_n^{-1}e_{n-1}^{-1}\ldots e_1^{-1}$. An empty path coincides with
its inverse. A path $p$ is called {\em closed} whenever $\iota(p)=\tau(p)$.

An {\em orientation\/} on the graph $\Gamma$ is a subset $E^+$ of the
set $E$ of all edges such that, for any edge $e\in E$, there is exactly
one of the edges $e$, $e\iv$ that belongs to $E^+$. The edges in $E^+$ are
called {\em positive\/} and the edges in $E^-=E\setminus E^+$ are called
{\em negative}. A path on an oriented graph is called {\em positive\/}
whenever it involves positive edges only. (An empty path is always
positive.) For any path $p$, there are defined its {\em initial vertex\/}
$\iota(p)$ and its {\em terminal vertex\/} $\tau(p)$: if $p=1_v$, then
$\iota(p)=\tau(p)=v$; if $p=e_1\ldots e_n$, then $\iota(p)=\iota(e_1)$,
$\tau(p)=\tau(e_n)$. For any two paths $p$, $q$ such that $\tau(p)=\iota(q)$,
one can naturally define a {\em product\/} $p\cdot q$ of the paths $p$ and
$q$: for $p=e_1\ldots e_n$, $q=f_1\ldots f_m$ we put
$p\cdot q= e_1\ldots e_nf_1\ldots f_m$. If $p$ ($q$) is empty, then
$p\cdot q=q$ ($p\cdot q=p$).

An {\em orineted graph\/} is by definition a graph $\Gamma$ with a fixed
orientation $E^+$. It is clear that any graph admits an orientation.
\vspace{1ex}

The concept of a graph of groups will play an important role. Let us have
an oriented graph $\Gamma$, where $E^+$ is the set of positive edges. We
say that a {\em graph of groups\/} structure on the graph $\Gamma$ is given
whenever to each edge $e\in E^+$ we assign a group $G_e$, to each vertex
$v\in V$ we assign a group $G_v$, and we fix embeddings
$\iota_e\colon G_e\to G_{\iota(e)}$, $\tau_e\colon G_e\to G_{\tau(e)}$
for any $e\in E^+$.

In the construction described below, we will have a 2-complex structure on
$\Gamma$ together with the graphs of groups structure. This means that we
have a set $F$ which is disjoint from $V$ and $E$. This set is called the
{\em set of\/} $2$-{\em cells}. We also have a mapping that assigns a closed
path in $\Gamma$ to each element in $F$. This path is called the {\em defining
path\/} of the $2$-cell. Given a 2-complex, we define the homotopy relation
on the set of paths in a standard way. Also one can define the concept of the
fundamental group of ${\cal K}$ with basepoint $w$. We denote this group by
$\pi_1({\cal K},w)$.

We will consider $2$-complexes that have a graph of groups structure on their
$1$-skeletons. We shall call such structures $2$-{\em complexes of groups}.
The concept of a $2$-complex of groups already exists and it is used widely
in many papers (see~\cite{Haef}). Every 2-complex of groups in our sense is a
2-complex of groups in the sense of~\cite{Haef}, but not vice versa.
(In general, a $2$-complex of groups is a structure that has not only vertex
groups $G_v$ ($v\in V$) and edge groups $G_e$ ($e\in E$) but also cell groups
of the form $G_f$ ($f\in F$) that are assigned to $2$-cells. In our situation
all the cell groups $G_f$ are trivial.)

So, let $\cal G$ be a $2$-complex of groups. Now we define the fundamental
group of $\cal G$. One can define it in different ways. We will define
it as a fundamental group of an ordinary $2$-complex ${\cal_K}({\cal G})$
with a basepoint. Here is the description of the complex
${\cal_K}({\cal G})$.

We add new edges and new $2$-cells to the $2$-complex ${\cal K}$. For any
vertex $v\in V$ and for any element $g\in G_v$ we add an edge denoted by
$g_v$ that has $v$ as both initial and terminal vertex. The new $2$-cells
are of two types. $2$-cells of the first type have defining paths
$g_vh_v(gh)_v\iv$ for any vertex $v\in V$ and for any elements $g,h\in G_v$.
$2$-cells of the second type have defining paths
$e\iv g_{\iota(e)}eh_{\tau(e)}\iv$, where $g=\iota_e(x)$, $h=\tau_e(x)$,
$x\in G_e$, $e\in E^+$. (The $2$-cells of the second type correspond to all
pairs of the form $(e,x)$, where $e\in E^+$, $x\in G_e$.) Recall that
$\iota_e$, $\tau_e$ are embeddings of the group $G_e$ into the groups
$G_\iota(e)$, $G_{\tau(e)}$, respectively. The $2$-complex obtained from
${\cal K}$ by adding new $2$-cells will be denoted by ${\cal_K}({\cal G})$.
For any vertex $v\in V$, the fundamental group $\pi_1({\cal_K}({\cal G}),v)$
will be called the {\em fundamental group\/} of $2$-{\em complex of groups}
${\cal G}$ with {\em basepoint\/} $v$. It will be denoted by
$\pi_1({\cal G},v)$.

A standard way to compute the fundamental group of a $2$-complex
(see~\cite{Still80}) can be also applied to compute the fundamental group
of a $2$-complex of groups. Let we have a structure of a graph of groups
${\cal G}$ on the $1$-skeleton of a $2$-complex ${\cal K}$. Consider the
connected component ${\cal K}_w$ of ${\cal K}$ that contains vertex $w$ and
let us choose some maximal subtree ${\cal T}_w$ in this component. It will
also be a maximal subtree of the connected component of the new
$2$-complex ${\cal_K}({\cal G})$ that contains vertex $w$. Let us take
the set of positive edges $E^+_w$ of this connected component together
with elements of the form $g_v$, $v\in V_w$, $g\in G_v$, where $V_w$ is the
set of vertices of this connected component. The union of these sets is
the set of generators of the fundamental group. For defining relations, we
take all relations of the form $e=1$, where $e$ belongs to the tree
${\cal T}$, and also all relations of the form $r=1$, where $r$ is the
defining path of any $2$-cell of the $2$-complex ${\cal_K}({\cal G})$. The
group given by the described presentation is isomorphic to the fundamental
group $\pi_1({\cal G},w)$ of our $2$-complex of groups.

Let us give one more equivalent description. It is clear that the
fundamental group $\pi_1({\cal K},w)$ with basepoint $w$ of the original
$2$-complex ${\cal K}$ can be computed in the same way, choosing ${\cal T}_w$
as a maximal subtree of the connected component of $w$. Now each edge
$e$ of this component uniquely defines an element in $\pi_1({\cal K},w)$,
namely, the equivalence class of the path $p_{\iota(e)}ep_{\tau(e)}\iv$,
where $p_v$ denotes the geodesic path from $w$ to $v$ in the subtree
${\cal T}_w$. Then
\be{fgSq}
\left.\pi_1({\cal G},w)\cong\mathop{\ast}\limits_{v}G_v
\ast\pi_1({\cal K},w)\right/{\cal N},
\ee
where the free product of groups $G_v$ is taken over all vertices $v$ of
the connected component of ${\cal K}$ that contains $w$, and ${\cal N}$ is
the normal closure of the following relations:
\be{descr}
g_{\iota(e)}^e=h_{\tau(e)}\mbox{  for every }e\in E^+_w,\ x\in G_e,
\mbox{ where }g=\iota_e(x),\ h=\tau_e(x).
\ee
This description will be often used below.
\vspace{1ex}

Let us give one of the main definitions.

\begin{df}
\label{diagprod}
{\rm Let $X$ be an alphabet, let $H_x$ ($x\in X$) be an arbitrary family of
groups and let $\qq=\pres{X}{\sss}$ be a semigroup presentation, $w\in X^+$.
Consider the Squier complex ${\cal K}={\cal K}(\qq)$ and introduce the
structure of graph of groups on its $1$-skeleton in the following way.
Let $E^+$ be the set of all positive edges of ${\cal K}$, that is the set
of triples of the form $e=(u,s\to t,v)$, where $u,v\in X^*$, $(s=t)\in\sss$.
For any word $z=x_1x_2\ldots x_n$, where $x_1,x_2,\ldots,x_n\in X$, let
$$
H_u=H_{x_1}\times H_{x_2}\times\cdots H_{x_n};
$$
if $u$ is empty, then $H_u=1$. For any vertex $u\in{\cal K}$ let $G_u=H_u$;
for any edge $e=(u,s\to t,v)\in E^+$, where $u,v\in X^*$, $(s=t)\in\sss$,
let $G_e=H_u\times H_v$. We have the natural embeddings
$\iota_e\colon G_e\to G_{\iota(e)}$ as the embedding
$H_u\times H_v\to H_{usv}=H_u\times H_s\times H_t$ and
$\tau_e\colon G_e\to G_{\tau(e)}$ as an embedding
$H_u\times H_v\to H_{utv}=H_u\times H_t\times H_t$. This gives us a
$2$-complex of groups, which will be denoted by ${\cal K}_H$. The fundamental
group $\pi_1({\cal K}_H,w)$ of this $2$-complex of groups with basepoint $w$
is called the {\em diagram product\/} of the family
$H_X=\{\,H_x\ (x\in X)\,\}$ of groups over the presentation
$\qq=\pres{X}{\sss}$ with base $w$. It will be denoted by
${\cal D}(H_X;\sss,w)$.
}
\end{df}

It is easy to see that the diagram product ${\cal D}(H_X;\sss,w)$
coincides with the diagram group ${\cal D}(\qq,w)$ in the case when
the groups $H_x$ are trivial for all $x\in X$.
\vspace{0.5ex}

The main result about this construction is that the diagram product of
any family of diagram groups over a semigroup presentation is again a diagram
group. Let us formulate this result in its general form giving the
description of the presentation, for which the corresponding diagram
group is the diagram product.

\begin{thm}
\label{DiagProd}
Let $\qq=\pres{X}{\sss}$ be a semigroup presentation, $w\in X^+$.
To each $x\in X$ we assign a diagram group $G_x={\cal D}(\pp_x,w_x)$,
where $\pp_x=\pres{\Sigma_x}{\rr_x}$ are semigroup presentations,
$w_x\in\Sigma_x^+$ $(x\in X)$. Let $A=\{\,a_x\mid x\in X\,\}$ be some
alphabet. Assume that the alphabets $X$, $A$, $\Sigma_x$ $(x\in X)$ are
disjoint. Let
$$
\Sigma=\bigcup\limits_{x\in X}\Sigma_x,\ \
\rr=\bigcup\limits_{x\in X}\rr_x,\ \
{\cal W}=\bigcup\limits_{x\in X}\{\,x=a_xw_xa_x\,\}.
$$
Consider the presentation
$$
\pp=\pres{X\cup\Sigma\cup A}{\sss\cup\rr\cup{\cal W}}.
$$

We claim that the diagram group ${\cal D}(\pp,w)$ is isomorphic to the
diagram product ${\cal D}(G_X;\sss,w)$ of the family
$G_X=\{\,G_x\ (x\in X)\,\}$ of groups over the presentation
$\qq=\pres{X}\sss$ with base $w$.

In particular, the diagram product of any family of groups over a
semigroup presentation is a diagram group.
\end{thm}

Let us consider a geometric description of this construction. In the
above notation, each group $G_x$ is isomorphic to the diagram group
over the presentation $\hat\pp_x$ that is obtained from $\pp_x$ by
adding letters $x$, $a_x$ to the alphabet $\Sigma_x$ and adding new relation
$x=a_xw_xa_x$ to the set $\rr_x$. We have $G_x\cong{\cal D}(\hat\pp_x,x)$.
Let us take any $(w,w)$-diagram over $\qq$ and consider some of its edges.
Let $x\in X$ be its label. We can cut the diagram along this edge and insert
any $(x,x)$-diagram over $\hat\pp_x$ in the resulting hole. We can do this
with all edges of the diagram. (If we insert trivial $(x,x)$-diagram, then
nothing changes.) After these transformations, we obtain some
$(w,w)$-diagram over $\pp$. One can show that diagrams obtained in this way
form the whole group ${\cal D}(\pp,w)$. From this point of view, diagrams
that represent elements in ${\cal D}(\pp,w)$ are obtained from
$(w,w)$-diagrams over $\qq$ by insertions of $(x,x)$-diagrams, which
represent elements in $G_x$.
\vspace{1ex}

\prf Let us construct the Squier complex for the presentation $\pp$. By
Kilibarda's Theorem, the fundamental group of this complex with basepoint
$w$ is isomorphic to the diagram group ${\cal D}(\pp,w)$. Our goal is to
transform the Squier complex into some new $2$-complex with the same
fundamental group. We will need to check that it will be isomorphic to
the fundamental group of a certain $2$-complex of groups, that is, to
the diagram products of our groups.
\vspace{1ex}

Let ${\cal K}_w(\pp)$ be the connected component of the Squier complex
over $\pp$, which contains the vertex $w$. A vertex $v$ in the same
component is an arbitrary word that equals $w$ modulo $\pp$. This word can
be uniquely decomposed into the product $v=v_1\ldots v_\mu$ ($\mu=\mu(v)$)
of subwords $v_1$, \dots, $v_\mu$ in such a way that each of them will be
either a letter in $X$, or a word of the form $a_xua_x$, where $u$ is
a word over $\Sigma_x$ that equals $w_x$ modulo $\pp_x$. The words $v_1$,
\dots, $v_\mu$ will be called the {\em factors\/} of the word $v$. To each
factor $v_i$ ($1\le i\le\mu$), we assign a letter in the alphabet $X$. This
letter will be denoted by $\pi(v_i)$. If $v_i\in X$, then we put
$\pi(v_i)=v_i$ and if $v_i$ has a form $a_xua_x$, where $u$ equals $w_x$
modulo $\pp_x$, then we put $\pi(v_i)=x$. Let us extend the function $\pi$,
putting by definition $\pi(v)=\pi(v_1)\ldots\pi(v_\mu)$ for any word $v$
that equals $w$ modulo $\pp$. The function $\pi$ will be called the
{\em projection}.

There is a natural two-sided action of the free monoid
$M=\left(X\cup\Sigma\cup A\right)^*$ on the Squier complex. It can be defined
by the following rule: for any $m_1,m_2\in M$ and for any vertex $v$ of the
Squier complex, let $m_1\ast v\ast m_2$ be the vertex $m_1vm_2$ and let
$m_1\ast e\ast m_2$ be the edge $(m_1u,p\to q,vm_2)$, for any edge
$e=(u,p\to q,v)$. The images of a given subgraph of ${\cal K}(\pp)$ under
this action will be called the {\em shifts\/} of this subgraph.

Let $\hat\rr_x=\rr_x\cup\{x=a_xw_xa_x\}$. We introduce a presentation
$\hat\pp_x=\pres{\Sigma_x,x,a_x}{\hat\rr_x}$ for all $x\in X$. It is clear
that ${\cal D}(\hat\pp_x,x)\cong G_x$. Let $v$ be a vertex of the complex
${\cal K}_w(\pp)$ and let $v=v_1\ldots v_\mu$ be the decomposition of
the word $v$ into factors. It is obvious that for any $1\le i\le\mu$, the
letter $\pi(v_i)=x$ is equal to the word $v_i$ modulo $\hat\pp_x$. Also
it is clear that $v$ equals $\pi(v)$ modulo $\pp$. Note that Squier
complexes for presentations $\hat\pp_x$ and their shifts can be regarded
as subcomplexes of the Squier complex of $\pp$.

For each $x\in X$ we choose a maximal subtree ${\cal T}_x$ in the connected
component ${\cal K}(\hat\pp_x,x)$ of the Squier complex of the presentation
$\hat\pp_x$ that contains the vertex $x$. Let us also choose a maximal
subtree ${\cal T}_Q$ in the connected component ${\cal K}(\qq,w)$ of the
Squier complex of the presentation $\qq$ that contains vertex $w$. Let $v$ be
an arbitrary vertex of the complex ${\cal K}_w(\pp)$ and let
$v=v_1\ldots v_\mu$ be the decomposition of $v$ into factors. Let
$x_i=\pi(v_i)$ ($1\le i\le\mu$). Consider the following subgraphs of
${\cal K}_w(\pp)$:
\be{sbgrph}
1\ast{\cal T}_{x_1}\ast x_2\ldots x_\mu,\
v_1\ast{\cal T}_{x_2}\ast x_3\ldots x_\mu, \ldots,\
v_1\ldots v_{\mu-1}\ast{\cal T}_{x_\mu}\ast1.
\ee
Now consider the subgraph ${\cal T}$ of ${\cal K}_w(\pp)$ that is a union
of subgraphs~(\ref{sbgrph}) for all vertices $v$ that are equal to $w$
modulo $\pp$, together with the subgraph ${\cal T}_Q$. Let us prove that
${\cal T}$ is a maximal subtree of ${\cal K}_w(\pp)$.

First of all we will establish that ${\cal T}$ is a connected subgraph
that contains all vertices of ${\cal K}_w(\pp)$, that is, for any vertex
$v$ of our component, we will find a path in ${\cal T}$ from $w$ to $v$.
Let $v=v_1\ldots v_\mu$ be the decomposition of $v$ into a product of
factors. By $p$ we denote the geodesic path in ${\cal T}_Q$ from $w$ to
$\pi(v)=x_1\ldots x_\mu$, where $x_i=\pi(v_i)$ for all $i$ from $1$ to $\mu$.
For each $i$, let $p_i$ be the geodesic path from $x_i$ to $v_i$ in the graph
${\cal T}_{x_i}$. For each $i$ from $1$ to $\mu$ we consider the path
$\tilde p_i=v_1\ldots v_{i-1}\ast p_i\ast x_{i+1}\ldots x_\mu$. Obviously,
it connects vertices $v_1\ldots v_{i-1}x_i\ldots x_\mu$ and
$v_1\ldots v_ix_{i+1}\ldots x_\mu$ in the graph ${\cal T}$. The product
$p\tilde p_1\ldots\tilde p_\mu$ is a path in ${\cal T}$ from $w$ to
$v_1\ldots v_\mu=v$.

Now let us prove that ${\cal T}$ has no nontrivial cycles. We argue by
contradiction. Suppose that a nontrivial cycle exists. If it does not
consist of edges that belong to subgraphs of the form~(\ref{sbgrph}),
then it has an edge from ${\cal T}_Q$. Since ${\cal T}_Q$ has no nontrivial
cycles, our cycle must contain edges from subgraphs of the
form~(\ref{sbgrph}). Hence our cycle has a nontrivial cyclic subpath $\rho$
that is a loop at some vertex in ${\cal K}_w(\qq)$ and all its edges are
from subgraphs of the form~(\ref{sbgrph}). Let us make a simple but important
observation: the endpoints of each edge that belong to any shift of the
subcomplex ${\cal K}(\hat\pp_x)$ ($x\in X$), have equal projection. Indeed,
applying relations of the form $x=a_xw_xa_x$ does not change the projection,
and applying relations from $\rr_x$ occurs within a factor of the form
$a_xua_x$, where $u$ is a word over $\Sigma_x$. This also does not change the
projection, (in the last case one can see the role of the auxiliary alphabet
$A$). Thus projections of all vertices of the cyclic path $\rho$ coincide,
that is, $\rho$ is a nontrivial cycle in ${\cal T}$ that consists of edges
from subgraphs of the form~(\ref{sbgrph}). So in any case there is a
nontrivial cycle $\rho$ with the above property. Without loss of generality,
one can assume that $\rho$ does not contain occurrences of adjacent edges
that are mutually inverse.

Let $v=v_1\ldots v_\mu$ be the decomposition of $v$ into factors, where
$\rho$ is the loop at $v$. Each edge of the path $\rho$ touches exactly
one of the factors, as we could see above. Let $j$ be the greatest number
such that an edge of $\rho$ touches $j$th factor. Let $x_i=\pi(v_i)$
($1\le i\le\mu$). It follows from the structure of subgraphs~(\ref{sbgrph})
that $v_i=x_i\in X$ for all $j<i\le\mu$. Since the $j$th factor occurs
in the process of application of relations from $\hat\rr_j$ to it, the
path $\rho$ or one of its cyclic shifts has a maximal subpath $\rho'$ that
consists of edges that touch the $j$th factor only. Let $v'$ and $v''$ be
the initial and the terminal points of $\rho'$ respectively. We claim
that $v'=v''$. Suppose this is not true. It is clear that $v'$ and $v''$
differ by the $j$th factor only and so one can say that the $j$th
factor is not equal to $x_j$ either in $v'$ or in $v''$. Assume that
the $j$th factor of $v''$ is not equal $x_j$. By our assumption, $\rho'$
has fewer edges than $\rho$ (otherwise $v'=v''$ automatically). So there
is an edge $e$ such that $\rho'e$ is a subpath of some cyclic shift of
the path $\rho$. Since $\rho'$ was chosen maximal, the edge $e$
does not touch the $j$th factor. It also cannot touch a factor with a
number greater than $j$ because $j$ is maximal with this property. But it
also cannot touch a factor with a number less than $j$ because it belongs
to a subgraph of the form~(\ref{sbgrph}), and the $j$th factor of the
initial point of $e$ does not belong to $X$, a contradiction. So
$\rho'$ is a nontrivial cycle that belongs to a shift of the tree
${\cal T}_{x_j}$. However, this is impossible since a shift of a tree
is a tree itself. This contradiction shows that ${\cal T}$ has no
nontrivial cycles. Applying what we have said above, we conclude that
${\cal T}$ is a maximal subtree in ${\cal K}_w(\pp)$.
\vspace{0.5ex}

Now we need to calculate the fundamental group $G=\pi_1({\cal K}_w(\pp))$
by using the maximal subtree ${\cal T}$. All edges of the complex
${\cal K}_w(\pp)$ are regarded as elements of the group $G$, and the edges
from ${\cal T}$ equal the identity in $G$. Paths in this complex, regarded
as products of edges, are just elements of the group $G$. To understand
how the other relations in $G$ look like, we need to describe the $2$-cells
in ${\cal K}_w(\pp)$. First of all let us mention some important property.
Recall that $M$ is the free monoid over the alphabet of the presentation
$\pp$, and $M$ acts both from the left and from the right on the complex
${\cal K}(\pp)$. Let $s$, $t$, $u$, $v$ be elements in $M$, each decomposed
into the product of factors, and let $s$ equals $t$ modulo $\pp$, $usv$
equals $w$ modulo $\pp$. Let us take an arbitrary path $p$ in ${\cal K}(\pp)$
that connects vertices $s$ and $t$. It is clear that the paths
$u\ast p\ast v$ and $\pi(u)\ast p\ast\pi(v)$ belong to ${\cal K}_w(\pp)$. We
claim that the following equality holds
\be{rav}
u\ast p\ast v=\pi(u)\ast p\ast\pi(v)
\ee
in the group $G$. This is what we are going to prove. To prove that,
one can consider the contours of the corresponding $2$-cells
as words in the generators of the group $G$ which are equal to the identity
in $G$. However, we think that one can check equality~(\ref{rav}) easier,
using the Kilibarda Theorem. Namely, to prove the equality~(\ref{rav}), it
suffices to use the fact that $G$ is isomorphic to the diagram group over
$\pp$ with base $\pi(usv)$. So let us find the diagrams over $\pp$ that
represent the elements in $G$ from both sides of equality~(\ref{rav}), and
then let us check that the diagrams are equal.

Let $x=x_1\ldots x_m$ be any word in $M$ decomposed into the product of
its factors. For any $i$ from $1$ to $m$ let $q_i$ be the geodesic path from
$\pi(x_i)$ to $x_i$ in the tree ${\cal T}_{\pi(x_i)}$. Then
$$
p_x=(1\ast q_1\ast\pi(x_2\ldots x_m))(x_1\ast q_2\ast\pi(x_3\ldots x_m))
\ldots(x_1\ldots x_{m-1}\ast q_m\ast1)
$$
is a path from $\pi(x)$ to $x$. Let $\Delta_x$ be the diagram represented
by it. Let us consider such paths and diagrams for all $x\in\{s,t,u,v\}$.
Also let $q$ be the geodesic path from $\pi(usv)$ to $\pi(utv)$ in the tree
${\cal T}_Q$, and let $\Delta$, $\Psi$ be the diagrams represented by $p$,
$q$, respectively. To find the spherical diagram with base $\pi(usv)$
represented by the path $u\ast p\ast v$ (via the isomorphism of the diagram
group and the group $G$), one needs to concatenate three diagrams:
the $\Delta_1$ that corresponds to the path in the tree ${\cal T}$ from
$\pi(usv)$ to $usv$, the diagram $\Delta_2=\ve(u)+\Delta+\ve(v)$ (that
corresponds to the path $u\ast p\ast v$), and the diagram $\Delta_3$ that
corresponds to the path in the tree ${\cal T}$ from $utv$ to $\pi(usv)$. So
consider the path $(1\ast p_u\ast\pi(sv))(u\ast p_s\pi(v))(us\ast p_v\ast1)$.
It follows from the description of subgraphs~(\ref{sbgrph}) that this path
is contained in ${\cal T}$. It corresponds to the diagram
$$
\Delta_1=
(\Delta_u+\ve(\pi(sv)))(\ve(u)+\Delta_s+\ve(\pi(v)))(\ve(us)+\Delta_v).
$$
Further, the path
$(ut\ast p_v\iv\ast1)(u\ast p_t\iv\ast\pi(v))(1\ast p_u\iv\ast\pi(tv))$
is contained in ${\cal T}$. Multiplying it by the path $q\iv$ on the right,
we obtain the path in ${\cal T}$ from $utv$ to $\pi(usv)$. This path
is represented by the diagram
$$
\Delta_3=(\ve(ut)+\Delta_v\iv)(\ve(u)+\Delta_t\iv+\ve(\pi(v)))
(\Delta_u\iv+\ve(\pi(tv)))\Psi\iv.
$$
Let us now multiply the diagrams $\Delta_1$, $\Delta_2$ and $\Delta_3$. It is
easy to see that the subdiagram $\Delta_u$ cancels with $\Delta_u\iv$ in this
product, and $\Delta_v$ cancels with $\Delta_v\iv$ (see the picture below).

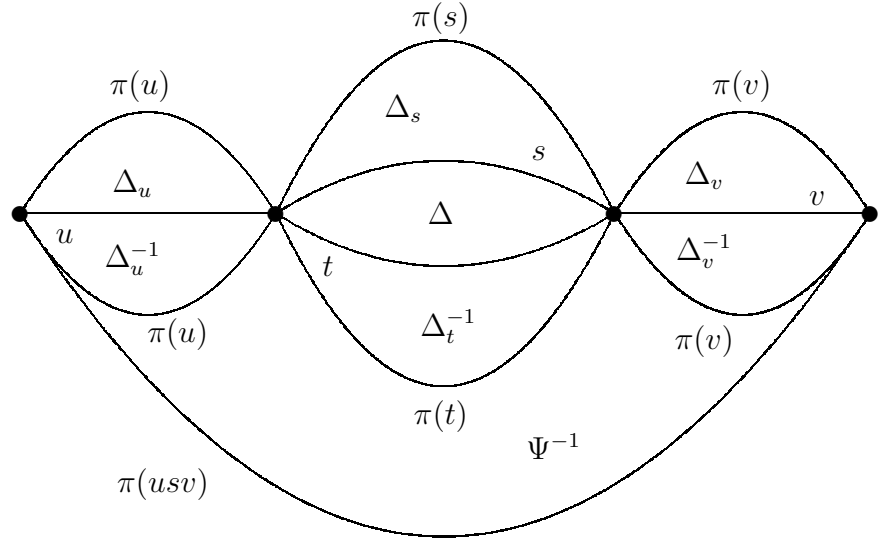
\begin{center} 
\unitlength=1mm
\special{em:linewidth 0.4pt}
\linethickness{0.4pt}
\begin{picture}(115.00,70.00)
\put(1.00,44.00){\circle*{2.00}}
\put(35.00,44.00){\circle*{2.00}}
\put(80.00,44.00){\circle*{2.00}}
\put(114.00,44.00){\circle*{2.00}}
\put(1.00,44.00){\line(1,0){34.00}}
\put(80.00,44.00){\line(1,0){34.00}}
\bezier{256}(1.00,44.00)(18.00,71.00)(35.00,44.00)
\bezier{256}(1.00,44.00)(18.00,17.00)(35.00,44.00)
\bezier{256}(80.00,44.00)(97.00,71.00)(114.00,44.00)
\bezier{256}(80.00,44.00)(97.00,17.00)(114.00,44.00)
\bezier{212}(35.00,44.00)(57.00,58.00)(80.00,44.00)
\bezier{212}(35.00,44.00)(57.00,30.00)(80.00,44.00)
\bezier{408}(35.00,44.00)(57.00,90.00)(80.00,44.00)
\bezier{408}(35.00,44.00)(57.00,-2.00)(80.00,44.00)
\bezier{824}(1.00,44.00)(57.00,-42.00)(114.00,44.00)
\put(57.00,70.00){\makebox(0,0)[cc]{$\pi(s)$}}
\put(17.00,61.00){\makebox(0,0)[cc]{$\pi(u)$}}
\put(97.00,61.00){\makebox(0,0)[cc]{$\pi(v)$}}
\put(22.00,28.00){\makebox(0,0)[cc]{$\pi(u)$}}
\put(92.00,27.00){\makebox(0,0)[cc]{$\pi(v)$}}
\put(7.00,41.00){\makebox(0,0)[cc]{$u$}}
\put(107.00,46.00){\makebox(0,0)[cc]{$v$}}
\put(57.00,17.00){\makebox(0,0)[cc]{$\pi(t)$}}
\put(70.00,52.00){\makebox(0,0)[cc]{$s$}}
\put(42.00,37.00){\makebox(0,0)[cc]{$t$}}
\put(52.00,58.00){\makebox(0,0)[cc]{$\Delta_s$}}
\put(57.00,44.00){\makebox(0,0)[cc]{$\Delta$}}
\put(58.00,29.00){\makebox(0,0)[cc]{$\Delta_t^{-1}$}}
\put(16.00,48.00){\makebox(0,0)[cc]{$\Delta_u$}}
\put(16.00,38.00){\makebox(0,0)[cc]{$\Delta_u^{-1}$}}
\put(92.00,49.00){\makebox(0,0)[cc]{$\Delta_v$}}
\put(92.00,39.00){\makebox(0,0)[cc]{$\Delta_v^{-1}$}}
\put(72.00,13.00){\makebox(0,0)[cc]{$\Psi^{-1}$}}
\put(20.00,8.00){\makebox(0,0)[cc]{$\pi(usv)$}}
\end{picture}
\end{center}

After cancelling $\Delta_u$ and $\Delta_u\iv$, $\Delta_v$ and
$\Delta_v\iv$, we obtain a diagram that is a product
\be{dddd}
(\ve(\pi(u))+\Delta_s+\ve(\pi(v)))(\ve(\pi(u))+\Delta+\ve(\pi(v)))
(\ve(\pi(u))+\Delta_t\iv+\ve(\pi(v)))\Psi\iv.
\ee
Repeating the arguments of the above paragraph for the path
$\pi(u)\ast p\ast\pi(v)$, it is easy to see that this path is represented
by the diagram~(\ref{dddd}) in the diagram group over $\pp$ with the base
$\pi(usv)$. This proves the equality~(\ref{rav}).

Let us consider an arbitrary edge $(u,s\to t,v)$ of the complex
${\cal K}_w(\pp)$. Let $(s=t)\in\sss$. The words $u$, $v$ can be decomposed
into products of factors, and the equality
$(u,s\to t,v)=(\pi(u),s\to t,\pi(v))$ holds in $G$. The right-hand side of
this equality can be regarded as an element of the group
$\pi_1({\cal K}_w(\qq))$, where ${\cal T}_Q$ is the maximal subtree in
${\cal K}_w(\qq)$. Now let $(s=t)\notin\sss$. In this case there exists
an element $x\in X$ and words $u_1$, $v_1$, $u_2$, $v_2$ such that
$u=u_1u_2$, $v=v_2v_1$, where $u_1$ ($v_1$) is the maximal prefix (suffix)
of the word $u$ (resp. $v$) that can be decomposed into a product of factors.
Here $u_2sv_2$ equals $x$ modulo $\pp$. Then by~(\ref{rav}) we have the
equality
$e=(u_1u_2,s\to t,v_2v_1)=\pi(u_1)\ast(u_2,s\to t,v_2)\ast\pi(v_1)$.
For any $x\in X$, let us consider the fundamental group
$\pi_1({\cal K}(\hat\pp_x),x)\cong G_x$ that can be calculated using the
maximal subtree ${\cal T}_x$ in the connected component of the Squier complex
over $\pp_x$ that contains $x$. The edges of this component will be thus
the elements of a group isomorphic to $G_x$ so the edge $e$ will belong to
an isomorphic copy of this group that is generated by edges obtained as
a result of shifts. Namely, let $U,V\in X^*$, $x\in X$. We consider the
group denoted by $U\ast G_x\ast V$. It is generated by edges of the form
$U\ast f\ast V$, where $f$ runs over edges that generate the group
$\pi_1({\cal K}(\hat\pp_x),x)\cong G_x$. In this sense, the edge $e$
belongs to the group $\pi(u_1)\ast G_x\ast\pi(v_1)$. The argument of this
paragraph can be summarized as follows: the group $G$ is generated by the
subgroup $\pi_1({\cal K}(\qq),w)\cong{\cal D}(\qq,w)$ and groups of the
form $u\ast G_x\ast v$, where $u,v\in X^*$, $x\in X$, and $uxv$ equals $w$
modulo $\pp$.

Now it remains to analyze all $2$-cells of ${\cal K}_w$ and to find out
what will be the relations between the generators of $G$ described above.
According to the description of $2$-cells in a Squier complex
given in Section~\ref{Prelim}, let a $2$-cell be given by a $5$-tuple
$(u,\ell_1\to r_1,z,\ell_2\to r_2,v)$, where $(\ell_1,r_1)$, $(\ell_2,r_2)$
belong to $\rr\cup\sss\cup{\cal W}$. Note that the word $u\ell_1z\ell_2v$
equals $w$ modulo $\pp$. Let us consider several cases depending on
the defining relations involved. The relation between edges that is obtained
from the given $2$-cell, has the form
\be{defpath}
(u,\ell_1\to r_1,z\ell_2v)(ur_1z,\ell_2\to r_2,v)=
(u\ell_1z,\ell_2\to r_2,v)(u,\ell_1\to r_1,zr_2v).
\ee

a) Let $(\ell_1,r_1)$, $(\ell_2,r_2)$ both belong to $\sss$. Then each of
the words $u$, $v$, $z$ can be decomposed into the product of factors.
Using equality~(\ref{rav}), one can replace in~(\ref{defpath}) the words
$u$, $v$, $z$ by their projections (we use the fact that each of the words
$\ell_j$, $r_j$ ($j=1,2$) coincides with its projection). Thus one can assume
that the words $u$, $v$, $z$ in~(\ref{defpath}) belong to $X^*$.
Then~(\ref{defpath}) is a defining relation of the group
$\pi_1({\cal K}(\qq),w)\cong{\cal D}(\qq,w)$ (calculated by using the
maximal subtree ${\cal T}_Q$).

b) Suppose that none of the relations $(\ell_1,r_1)$, $(\ell_2,r_2)$ belongs
to $\sss$. Suppose also that these relations are applied to different factors
of the word $u\ell_1z\ell_2v$. This means that there exist letters $x,y\in X$
and decompositions of the form $u=u_1u_2$, $z=z'z_0z''$, $v=v_2v_1$, where
$u_1$ is the maximal prefix of $u$ that is a product of factors, $v_1$ is the
maximal suffix of $v$ that is a product of factors, and $z_0$ is a maximal
subword of the word $z$ that is a product of factors. (It is not hard to
see that this word can be found uniquely.) Here $x$ equals $u_2\ell_1z'$ and
$y$ equals $z''\ell_2v_2$ (equalities are considered modulo $\hat\pp_x$ and
$\hat\pp_y$, respectively). Let us substitute the decompositions of the words
$u$, $v$, $z$ in the equality~(\ref{defpath}) using the fact that
$\pi(u_2\ell_1z')=\pi(u_2r_1z')=x$, $\pi(z''\ell_2v_2)=\pi(z''r_2v_2)=y$.
We obtain that the elements
$\pi(u_1)\ast(u_2,\ell_1\to r_1,z')\ast\pi(z_0)y\pi(v_1)$ ¨
$\pi(u_1)x\pi(z_0)\ast(z'',\ell_2\to r_2,v_2)\ast\pi(v_1)$ commute.
The first of them belongs to the group
$\pi(u_1)\ast G_x\ast\pi(z_0)y\pi(v_1)$, and the second one belongs to
the group $\pi(u_1)x\pi(z_0)\ast G_y\ast\pi(v_1)$. Let $U$, $V$, $Z$ be
arbitrary words over $X$ and let $x,y\in X$ be arbitrary letters such that
the word $UxZyV$ equals $w$ modulo $\pp$. We can conclude that any element
in $U\ast G_x\ast ZyV$ commutes with any element in $UxZ\ast G_y\ast V$
since for any edges $e$ and $f$ that belong to the generating sets of the
groups $G_x$ and $G_y$ respectively, one can find a suitable $2$-cell of
the form described above in such a way that the defining relations obtained
from it will be the relation of commutativity of $U\ast e\ast ZyV$ and
$UxZ\ast f\ast V$. Thus we get relations of the form
$[U\ast G_x\ast ZyV,UxZ\ast G_y\ast V]=1$, where $UxZyV$ equals $w$ modulo
$\pp$, $U,V,Z\in X^*$, $x,y\in X$.

c) Again, let none of the relations $(\ell_1,r_1)$, $(\ell_2,r_2)$ belong to
$\sss$ but assume now that the relations $(\ell_1=r_1)$ and $(\ell_2=r_2)$
are applied to the same factor of the word $u\ell_1z\ell_2v$.
This means that there exists a letter $x\in X$ and decompositions $u=u_1u_2$,
$v=v_2v_1$, where $u_1$ is the maximal prefix of the word $u$ that is a
product of factors, $v_1$ is the maximal suffix of the word $v$ that is a
product of factors. Now $x$ equals $u_2\ell_1z\ell_2v_2$ modulo $\hat\pp_x$.
Consider a $2$-cell of the Squier complex over $\hat\pp_x$ that corresponds
to the $5$-tuple $(u_2,\ell_1\to r_1,z,\ell_2\to r_2,v_2)$. All cells of
this form lead to the defining relations of a group isomorphic to $G_x$.
Acting on this $2$-cell by the element $\pi(u_1)$ on the left and by the
element $\pi(v_1)$ on the right, we get a $2$-cell of the complex
${\cal K}_w(\pp)$. The relation written on its contour is equivalent
to~(\ref{defpath}) if one takes the equality~(\ref{rav}) into account. Thus
we get the defining relations of all groups of the form $U\ast G_x\ast V$,
where $UxV$ equals $w$ modulo $\pp$, $U,V\in X^*$, $x\in X$.

d) Suppose that one of the relations $(\ell_1,r_1)$, $(\ell_2,r_2)$ belongs
to $\sss$ and the other one does not. First of all, let $(\ell_1,r_1)\in\sss$.
Then we have decompositions of the form $z=z_0z''$, $v=v_2v_1$, where $z_0$,
$v_1$ are products of factors that are chosen to be minimal with respect
to this property, as above. Now $z''\ell_2v_2$ equals $x$ modulo $\hat\pp_x$
for some letter $x\in X$. Let $f=(z'',\ell_2\to r_2,v_2)$. Substituting
the decompositions of words $v$, $z$ in~(\ref{defpath}), taking into
account that $\pi(z''\ell_2v_2)=\pi(z''r_2v_2)=x$ and applying~(\ref{rav}),
we obtain the following equality:
$$
(\pi(u),\ell_1\to r_1,\pi(z_0)x\pi(v_1))\cdot
(\pi(u)r_1\pi(z_0)\ast f\ast\pi(v_1))=
$$
$$
(\pi(u)\ell_1\pi(z_0)\ast f\ast\pi(v_1))\cdot
(\pi(u),\ell_1\to r_1,\pi(z_0)x\pi(v_1))\phantom{=}.
$$
Thus for each $x\in X$ and for any words $U$, $Z$, $V$ over $X$ such that
$(\ell_1,r_1)\in\sss$ and $U\ell_1ZxV$ equals $w$ modulo $\pp$, we obtain
the relations
\be{Soot1e}
Ur_1Z\ast f\ast V=(U\ell_1Z\ast f\ast V)^e,
\ee
where $e=(U,\ell_1\to r_1,ZxV)$, and $f$ runs over the generating set of
the group $\pi_1(\hat\pp_x,x)\cong G_x$.

Analogously, if $(\ell_2,r_2)\in\sss$, then we get the relations
\be{Soot2e}
U\ast f\ast Zr_2V=(U\ast f\ast Z\ell_2V)^e,
\ee
where $U$, $Z$, $V$ are words over $X$, $x\in X$, $e=(UxZ,\ell_2\to r_2,V)$,
and $f$ runs over the generating set of the group
$\pi_1(\hat\pp_x,x)\cong G_x$.
\vspace{1ex}

To compute the diagram product, let us define a structure of a graph of
groups on the $1$-skeleton of the Squier complex ${\cal K}(\qq)$. Our diagram
product is the fundamental group of the corresponding $2$-complex of groups.
We will apply the above described procedure of computing a fundamental
group of a $2$-complex of groups and we will then compare it with the
presentation of $G$. Let $U$, $V$ be words over $X$, $x\in X$. By
$H(U\cdot x\cdot V)$ we denote the group $U\ast G_x\ast V$ isomorphic to
$G_x$. Then, for any word $u=u_1\ldots u_m$, where $u_i\in X$ ($1\le i\le m$),
we denote by $H_u$ the free product of the groups of the form
$$
H(u_1\ldots u_{i-1}\cdot u_i\cdot u_{i+1}\ldots u_m)
$$
over all $i$ from $1$ to $m$. These groups are assigned to vertices of
${\cal K}(\pp)$. Now let us take an edge $e=(u,s\to t,v)$, where $u,v\in X^*$,
$(s=t)\in\sss$. The group $H_e=H_u\times H_v$ is assigned to it. The
maps $\iota_e$, $\tau_e$ naturally embed $H_e=H_u\times H_v$ into the
groups $H_{usv}\cong H_u\times H_s\times H_v$ and
$H_{utv}\cong H_u\times H_t\times H_v$, respectively (here $H_u$ maps onto
$H_u$,  nd $H_v$ maps onto $H_v$). It is clear that, instead of presenting
edge groups of the form $G_e$, one can present an isomorphism of some
subgroup of $H_{\iota(e)}$ to some subgroup of $H_{\tau(e)}$, for each egde
$e$. In our case this isomorphism is very simple: it maps the subgroup
$H_u\times\{\,1\,\}\times H_v$ ¢ $H_{usv}$ onto
$H_u\times\{\,1\,\}\times H_v$ ¢ $H_{utv}$. In this case we will speak
about the isomorphism induced by an edge $e$.

The groups of the form $H(U\cdot x\cdot V)$ will be presented as groups
generated by edges of the form $U\ast f\ast V$, where $f$ runs over the set
of edges of the connected component of the Squier complex
${\cal K}(\hat\pp_x)$ that contain vertex $x$. Here edges satisfy
the relations $U\ast f\ast V=1$ whenever $f$ belongs to the tree ${\cal T}_x$,
and also relations $U\ast r\ast V=1$, where $r$ is the defining path of
a $2$-cell of this complex. These relations of the group $G$ were
obtained in subsection c). For the direct products of the groups of the
form $H(U\cdot x\cdot V)$, we introduce relations of commutativity: each
element in the group $H(U\cdot x\cdot ZyV)$ commutes with each element in
the group $H(UxZ\cdot y\cdot V)$. For the group $G$, such relations were
obtained in subsection b).

Let us have an edge $e=(u,s\to t,v)$ that belongs to the complex ${\cal K}_w$.
Let $x,y\in X$, $u=u_1xu_2$, $v=v_1yv_2$. The isomorphism induced by the edge
$e$, takes $H(usv_1\cdot y\cdot v_2)$ to $H(utv_1\cdot y\cdot v_2)$, and
it takes $H(u_1\cdot x\cdot u_2yv)$ to $H(u_1\cdot x\cdot u_2yv)$. So,
according to ~(\ref{descr}, the conjugation by the edge $e$ leads to
the following relations
\be{Soot1E}
(usv_1\ast f\ast v_2)^e=(utv_1\ast f\ast v_2),
\ee
\be{Soot2E}
(u_1\ast f\ast u_2sv)^e=(u_1\ast f\ast u_2tv),
\ee
where $f$ runs over the set of edges that generate the corresponding group
in each of the cases. These relations coincide with relations~(\ref{Soot1e})
and (\ref{Soot2e}) of the group $G$ from subsection d). Finally, we
represent the group $\pi_1({\cal K}(\qq),w)$ as a group generated by edges
of ${\cal K}_w(\qq)$, claiming that the edges in ${\cal T}_Q$ are equal to
the identity and adding relations that correspond to the defining paths
of $2$-cells of this complex. Such relations of the group $G$ are described
in subsection a). Thus the quotient group of the free product of the group
$\pi_1({\cal K}(\qq),w)$ and groups of the form $H_u$ for all vertices $u$
of ${\cal K}_w(\qq)$, by the normal closure of relations~(\ref{Soot1E}) and
(\ref{Soot2E}), is given by the same generators and defining relations as
$G$. This means that the diagram product ${\cal D}(G_X;\sss,w)$ of the family
$G_X=\{\,G_x\ (x\in X)\,\}$ of groups over the presentation $\qq=\pres{X}\sss$
with base $w$ is isomorphic to the group $G=\pi_1({\cal K}(\pp,w)$, that
is, to the diagram group $\dg$.

The Theorem is proved.
\vspace{2ex}

Now let us consider a few applications of Theorem~\ref{DiagProd}. The first
three of them deal with already known constructions. We give them to
demonstrate that all group-theoretical constructions for diagram groups we
dealt with earlier (see~\cite[Section 8]{GuSa97}) are examples of diagram
products. Then we show that the class of diagram groups is closed under some
new operations: countable direct powers, restricted wreath products with the
group $\bf Z$, and also under some new special construction that will be used
in Section~\ref{SubConj}.

\begin{ex}
\label{ex-dir}
{\rm Let $X=\{\,x_1,\ldots,x_n\,\}$ be a finite alphabet. Consider the
presentation $\qq=\pres{X}{\emptyset}$ with empty set of defining relations
and let $w=x_1\ldots x_n$. To each letter $x_i$, we assign an arbitrary
group $G_i$ ($1\le i\le n$). It is obvious that the connected component
of the Squier complex of $\qq$, which contains $w$, consists of exactly
one vertex $w$. In the corresponding graph of groups, we have the group
$G_w=G_1\times\cdots\times G_n$. Obviously, it is the fundamental group
of the $2$-complex of groups from the definition of a diagram product.
Thus the diagram product ${\cal D}(G_X;\sss,w)$ of the family
$G_X=\{\,G_i\ (1\le i\le n)\,\}$ of groups over the presentation $\qq$ with
base $w$ is the direct product $G_1\times\cdots\times G_n$.
}
\end{ex}

\begin{ex}
\label{ex-free}
{\rm Let $I$ be a nonempty set and let $G_i$ ($i\in I$) be an arbitrary
family of groups. Let us consider an alphabet
$X=\{\,x\,\}\cup\{\,x_i\ (i\in I)\,\}$ and let $\qq=\pres{X}{\sss}$,
where $\sss$ consists of relations of the form $x=x_i$ for all $i\in I$.
Let $G_X$ be a family of groups that assigns the trivial group to the
letter $x$ and the group $G_i$ to the letter $x_i$ ($i\in I$). The
connected component of the Squier complex ${\cal K}(\qq)$ containing $x$
is a tree in which the vertex $x$ is connected by edges with all vertices
labelled by $x_i$ ($i\in I$). Let us consider the structure of a graph of groups
on the $1$-skeleton of the connected component of this Squier complex according to the
definition~\ref{diagprod}. It is easy to see that all edge groups are
trivial. From this, using description~(\ref{fgSq}), it is easy to see that
the fundamental group of the resulting $2$-complex of groups is the free
product of groups $G_i$ ($i\in I$). So the diagram product
${\cal D}(G_X;\sss,x)$ of the family $G_X=\{\,G_i\ (i\in I)\,\}$ of groups
over the presentation $\qq$ with base $x$ is the free product $\ast\,G_i$,
$i\in I$.
}
\end{ex}

\begin{ex}
\label{ex-bullet}
{\rm Let $G$, $H$ be any groups. Let us consider the presentation
$\qq=\pres{X}{\sss}$, where $X=\{\,x,y,z\,\}$, $\sss=\{\,x=xy,z=yz\,\}$.
Let $G_x=G$, $G_y=1$, $G_z=H$ and consider the diagram product
${\cal D}(G_X;\sss,xz)$ of the family $G_X=\{\,G_x,G_y,G_z\,\}$ of
groups over the presentation $\qq$ with base $xz$. We obtain that it is
isomorphic to the group $G\bullet H$, where $\bullet$ is the operation
defined in~\cite{GuSa97}. One can check this directly by comparing
the Theorem~\ref{DiagProd} and the definition of the operation $\bullet$
in~\cite{GuSa97}. Indeed, the group $G\bullet H$ can be described in the
following way. Consider countable number of copies $G_i$ of the group $G$
and countable number of copies $H_i$ of the group $H$ ($i\in\zz$).
An infinite cyclic group $\la z\ra$ acts on the group
\be{bul_}
(\ast\,G_i)\times(\ast\,H_i)
\ee
(free products are taken over all $i\in\zz$) permuting the factors: it takes
$G_i$ to $G_{i+1}$ and $H_i$ to $H_{i+1}$ for all integers $i$. The group
$G\bullet H$ is the semidirect product of the group~(\ref{bul_}) and the
group $\la z\ra$.
}
\end{ex}

\begin{ex}
\label{ex-DirPow}
{\rm
Let $G$ be an arbitrary group. Let us consider the presentation
$\qq=\pres{X}{\sss}$, where $X=\{\,x,y\,\}$, $\sss=\{\,x=xy\,\}$. Let
$G_x=1$, $G_y=G$. Consider the diagram product ${\cal D}(G_X;\sss,x)$ of the
family $G_X=\{\,G_x,G_y\,\}$ of groups over the presentation $\qq$ with base
$x$. Let us show that it is isomorphic to the countable direct power of the
group $G$.

The connected component ${\cal K}_x$ of the Squier complex over $\qq$
containing $x$ has the following form:

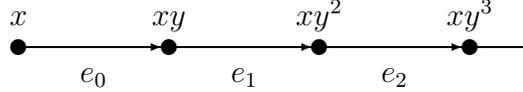
\begin{center} 
\unitlength=1mm
\special{em:linewidth 0.4pt}
\linethickness{0.4pt}
\begin{picture}(70.00,18.00)
\put(2.00,9.00){\circle*{2.00}}
\put(22.00,9.00){\circle*{2.00}}
\put(42.00,9.00){\circle*{2.00}}
\put(62.00,9.00){\circle*{2.00}}
\put(2.00,9.00){\vector(1,0){19.00}}
\put(22.00,9.00){\vector(1,0){19.00}}
\put(42.00,9.00){\vector(1,0){19.00}}
\put(62.00,9.00){\line(1,0){8.00}}
\put(2.00,12.00){\makebox(0,0)[cb]{$x$}}
\put(22.00,12.00){\makebox(0,0)[cb]{$xy$}}
\put(42.00,12.00){\makebox(0,0)[cb]{$xy^2$}}
\put(62.00,12.00){\makebox(0,0)[cb]{$xy^3$}}
\put(12.00,5.00){\makebox(0,0)[cc]{$e_0$}}
\put(32.00,5.00){\makebox(0,0)[cc]{$e_1$}}
\put(52.00,5.00){\makebox(0,0)[cc]{$e_2$}}
\end{picture}
\end{center}

Here vertices are all words of the form $xy^i$ ($i\ge0$), positive edges
have the form $e_i=(1,x\to xy,y^i)$ ($i\ge0$), and the maximal subtree
${\cal T}$ includes all these edges. This complex has no $2$-cells.
Thus it is obvious that its fundamental group is trivial. Using our
convention that (given a maximal subtree) all edges are regarded as elements
of the fundamental groups, we have equalities $e_i=1$ for all $i\ge0$.

Let us consider the structure of the graph of groups on the $1$-skeleton of
the complex ${\cal K}_x$ according to the definition of a diagram product.
We will obtain that the group $H_v=G_{n1}\times\cdots\times G_{nn}$, where
$G_{ni}$ ($1\le i\le n$) are groups isomorphic to $G$, corresponds to
the vertex $v=xy^n$ ($n\ge0$). Let us consider a positive edge
$e=e_n=(1,x\to xy,y^n)$ ($n\ge0$). The group $G_e$ is isomorphic to $G^n$,
the $n$th direct power of $G$. The embedding $\iota_e$ maps $G^n$
isomorphically onto $G_{n1}\times\cdots\times G_{nn}$, and the mapping
$\tau_e$ maps $G^n$ isomorphically onto the last $n$ factors of the
direct product $G_{n+1,1}\times G_{n+1,2}\cdots\times G_{n+1,n+1}$.
The relations from the description~(\ref{descr}), together with the equality
$e=1$, allow to identify corresponding elements of $G_{n1}$ and $G_{n+1,2}$,
\dots, $G_{nn}$ and $G_{n+1,n+1}$. Thus we can introduce the following
notation: $G_0=G_{11}=G_{22}=\cdots\,$, $G_1=G_{21}=G_{32}=\cdots\,$, \dots,
$G_n=G_{n+1,1}=G_{n+2,2}=\cdots\,$, \dots\,. Each of these groups is
isomorphic to $G$. They generate a countable direct power of $G$. So our
diagram product is the countable direct power of $G$.
}
\end{ex}

Theorem~\ref{DiagProd} implies the following result.

\begin{thm}
\label{dir-pow}
The class of diagram groups is closed under countable direct powers.
\end{thm}

\begin{ex}
\label{ex-wr}
{\rm
Let $G$ be arbitrary group. Consider the presentation $\qq=\pres{X}{\sss}$,
where $X=\{\,x,y,z\,\}$, $\sss=\{\,x=xy,z=yz\,\}$. Let $G_x=1$, $G_y=G$,
$G_z=1$ and consider the diagram product ${\cal D}(G_X;\sss,xz)$ of the
family $G_X=\{\,G_x,G_y,G_z\,\}$ of groups over the presentation $\qq$ with
the base $xz$. Let us show that it is isomorphic to the (restricted) wreath
product $G\wr\zz$.

The connected component ${\cal K}_{xz}$ of the Squier complex over $\qq$
that contains $xz$ has the following form:

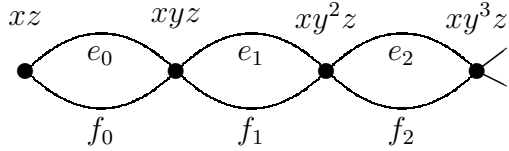
\begin{center} 
\unitlength=1mm
\special{em:linewidth 0.4pt}
\linethickness{0.4pt}
\begin{picture}(67.00,20.00)
\put(3.00,10.00){\circle*{2.00}}
\put(23.00,10.00){\circle*{2.00}}
\put(43.00,10.00){\circle*{2.00}}
\put(63.00,10.00){\circle*{2.00}}
\bezier{112}(3.00,10.00)(13.00,20.00)(23.00,10.00)
\bezier{112}(23.00,10.00)(33.00,20.00)(43.00,10.00)
\bezier{112}(43.00,10.00)(53.00,20.00)(63.00,10.00)
\bezier{112}(3.00,10.00)(13.00,0.00)(23.00,10.00)
\bezier{112}(23.00,10.00)(33.00,0.00)(43.00,10.00)
\bezier{112}(43.00,10.00)(53.00,0.00)(63.00,10.00)
\put(63.00,10.00){\line(4,3){4.00}}
\put(63.00,10.00){\line(2,-1){4.00}}
\put(3.00,17.00){\makebox(0,0)[cc]{$xz$}}
\put(23.00,17.00){\makebox(0,0)[cc]{$xyz$}}
\put(43.00,17.00){\makebox(0,0)[cc]{$xy^2z$}}
\put(63.00,17.00){\makebox(0,0)[cc]{$xy^3z$}}
\put(13.00,12.00){\makebox(0,0)[cc]{$e_0$}}
\put(33.00,12.00){\makebox(0,0)[cc]{$e_1$}}
\put(53.00,12.00){\makebox(0,0)[cc]{$e_2$}}
\put(13.00,2.00){\makebox(0,0)[cc]{$f_0$}}
\put(33.00,2.00){\makebox(0,0)[cc]{$f_1$}}
\put(53.00,2.00){\makebox(0,0)[cc]{$f_2$}}
\end{picture}
\end{center}

Here vertices are all words of the form $xy^iz$ ($i\ge0$), positive edges
have the form $e_i=(1,x\to xy,y^iz)$, $f_i=(xy^i,z\to yz,1)$ ($i\ge0$), and
the maximal subtree ${\cal T}$ consists of the edges $e_i$, $i\ge 0$. All
$2$-cells can be described as follows. Let $i\ge0$. Consider the vertex
$xy^iz$. The edges $e_i$, $f_i$ going out of it correspond to independent
transformations of words. So the given pair of edges defines two homotopic
paths $e_if_{i+1}$ and $f_ie_{i+1}$ that define a $2$-cell.

According to the convention that edges are regarded as elements of the
fundamental group $\pi_1({\cal K},xz)$, we have $e_i=1$ ($i\ge0$) in the
group. The equalities $e_if_{i+1}=f_ie_{i+1}$ that hold in this group imply
$f_i=f_{i+1}$ for all $i\ge0$. Let $f=f_0=f_1=f_2=\cdots\,$.

According to the definition of a diagram product, let us consider the
structure of the graph of groups on the $1$-skeleton of ${\cal K}_{xz}$.
The group $H_v=G_{n1}\times\cdots\times G_{nn}$ is assigned to
the vertex $v=xy^nz$ ($n\ge0$), where $G_{ni}$ ($1\le i\le n$) is a group
isomorphic to $G$. Let us consider relations~(\ref{descr}) that
correspond to positive edges. Let $e=e_n=(1,x\to xy,y^n)$ ($n\ge0$). As
in the previous example, using the equality $e=1$, we identify
corresponding elements of the groups $G_{ni}$ and $G_{n+1,i+1}$
($1\le i\le n$) and introduce the notation $G_0=G_{11}=G_{22}=\cdots\,$,
$G_1=G_{21}=G_{32}=\cdots\,$, \dots, $G_n=G_{n+1,1}=G_{n+2,2}=\cdots\,$,
\dots\,. As above, these groups generate a countable direct power of the
group $G$. Now let $e=f_n=(xy^n,z\to yz,1)$ ($n\ge0$). Consider relations
of the form~(\ref{descr}) that correspond to these edges. The group
$G_e$ is still the $n$th power of $G$. The embedding $\iota_e$ maps
$G^n$ onto $G_{n1}\times\cdots\times G_{nn}$ isomorphically, and the
embedding $\tau_e$ maps $G^n$ isomorphically onto the first $n$ factors
of the direct product $G_{n+1,1}\times G_{n+1,2}\cdots\times G_{n+1,n+1}$.
So the relations that correspond to the edge $f$ show that conjugation
by $f$ takes $G_{ni}$ to $G_{n+1,i}$ ($1\le i\le n$). Using our notation,
we obtain that the conjugation by $f$ takes the group $G_k$ to the group
$G_{k+1}$ for all $k\ge0$. Thus the diagram product we are considering
is generated by groups $G_0$, $G_1$, \dots\ and the element $f$. From this,
one can deduce that we have the restricted wreath product $G\wr\zz$.
}
\end{ex}

Applying Theorem~\ref{DiagProd}, we get one more result.

\begin{thm}
\label{wr}
The class of diagram groups is closed under restricted wreath products with
the infinite cyclic group $\zz$, that is, if $G$ is a diagram group, then
$G\wr\zz$ is also a diagram group.
\end{thm}

Note that if we take R.\,Thompson's group $F$ represented by diagrams over
$\pres{u}{uu=u}$ with base $u$ and consider the presentation
$\pres{x,u,z}{xu=x,uz=z,uu=u}$, then the diagram group over it with base
$xz$ will be isomorphic not to $F\wr\zz$ but to $F$. This can be checked
directly. To get the group $F\wr\zz$, one needs to represent the group $F$
by diagrams according to the statement of Theorem~\ref{DiagProd}. Namely,
one has to take the diagram group with base $y$ over the presentation
$\pres{u,a,y}{y=aua,uu=u}$. Then the diagram group with base $xz$ over
$\pp=\pres{x,y,z,a,u}{xy=x,yz=z,y=aua,uu=u}$ will be isomorphic to $F\wr\zz$.
The reader can easily list the presentations that lead to diagram groups of
the form $(\cdots((\zz\wr\zz)\wr\zz)\cdots)\wr\zz$.
\vspace{0.5ex}

Let us make one more remark. In~\cite{GuSa97} we constructed an example of
a diagram group that was finitely generated but not finitely presented
(Theorem~10.5). We took the group $\zz\bullet\zz$ for this purpose. It has
a presentation with three generators
$$
\zz\bullet\zz=\pres{a,b,t}{[a^{t^n},b]=1\ (n\ge0)}.
$$
Now we can also take the group $\zz\wr\zz$ as an example of a finitely
generated but not finitely presented group (the fact that $\zz\wr\zz$ has no
finite presentations can be easily proved using either HNN-extensions or
representation of groups by transformations). We have the following
presentation with two generators for this group:
$$
\zz\wr\zz=\pres{a,b}{[a^{b^n},a]=1\ (n\ge1)}.
$$

In the next example we deal with a more complicated construction.
At first sight one can think it is quite artificial. However, we will
efficiently use it later, in Section~\ref{SubConj}. Let us consider the
following group-theoretical construction. Take two groups, $G$ and
$H$. We assign to them a new group denoted by ${\cal O}(G,H)$. Let us
consider a countable family of copies $G_i$ of the group $G$,
and a coutable family of copies $H_i$ of $H$ ($i\in\zz$).
For any $i\in\zz$, let $g_i$ ($h_i$) be the element that corresponds to
$g\in G$ ($h\in H$). By $G^\infty$ ($H^\infty$) we denote a coutable direct
power of the groups $G_i$ ($H_i$) taken over all $i\in\zz$. Let
\be{OG}
{\cal O}(G,H)=\left.G^\infty\ast H^\infty\ast\la c\ra\right/{\cal N},
\ee
where ${\cal N}$ is the normal closure of the set of relations of the two
forms:
\be{soot1}
g_i^t=g_{i+1},\ \ h_i^t=h_{i+1}
\mbox{\ \ \ for all }i\in\zz,\ g\in G,\ h\in H;
\ee
\be{soot2}
[g_i,h_j]=1\mbox{\ \ for all }i,j\in\zz,\ g\in G,\ h\in H
\mbox{  such that  }i\le j.
\ee

\begin{ex}
\label{kh}
{\rm
Let $G$, $H$ be arbitrary groups. Let us consider the presentation
$\qq=\pres{X}{\sss}$, where $X=\{\,x,y,\by,z,p,q,r\,\}$,
$\sss=\{\,x=xyp,z=r\by z,pyq=q\by r\,\}$. Let $G_y=G$, $G_{\by}=H$,
$G_x=G_z=G_p=G_q=G_r=1$. Consider the diagram product ${\cal D}(G_X;\sss,w)$
of the family $G_X=\{\,G_x,G_y,G_{\by},G_z,G_p,G_q,G_r\,\}$ of groups over
the presentation $\qq$ with base $w=xyq\by z$. Let us show that it is
isomorphic to ${\cal O}(G,H)$.
\vspace{0.5ex}

The connected component ${\cal K}_w$ of the Squier complex over $\qq$
that contains $w$, has the following form:

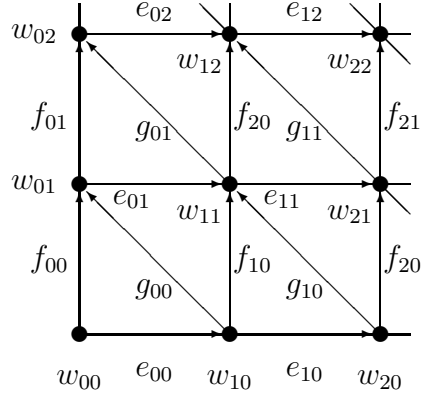
\begin{center} 
\unitlength=1mm
\special{em:linewidth 0.4pt}
\linethickness{0.4pt}
\begin{picture}(52.00,53.00)
\put(8.00,9.00){\circle*{2.00}}
\put(28.00,9.00){\circle*{2.00}}
\put(48.00,9.00){\circle*{2.00}}
\put(8.00,29.00){\circle*{2.00}}
\put(28.00,29.00){\circle*{2.00}}
\put(48.00,29.00){\circle*{2.00}}
\put(8.00,49.00){\circle*{2.00}}
\put(28.00,49.00){\circle*{2.00}}
\put(48.00,49.00){\circle*{2.00}}
\put(8.00,9.00){\vector(1,0){19.00}}
\put(28.00,9.00){\vector(1,0){19.00}}
\put(8.00,9.00){\vector(0,1){19.00}}
\put(8.00,29.00){\vector(0,1){19.00}}
\put(28.00,9.00){\vector(0,1){19.00}}
\put(28.00,29.00){\vector(0,1){19.00}}
\put(48.00,9.00){\vector(0,1){19.00}}
\put(48.00,29.00){\vector(0,1){19.00}}
\put(8.00,29.00){\vector(1,0){19.00}}
\put(28.00,29.00){\vector(1,0){19.00}}
\put(8.00,49.00){\vector(1,0){19.00}}
\put(28.00,49.00){\vector(1,0){19.00}}
\put(28.00,9.00){\vector(-1,1){19.00}}
\put(48.00,9.00){\vector(-1,1){19.00}}
\put(28.00,29.00){\vector(-1,1){19.00}}
\put(48.00,29.00){\vector(-1,1){19.00}}
\put(48.00,9.00){\line(1,0){4.00}}
\put(48.00,29.00){\line(1,0){4.00}}
\put(48.00,49.00){\line(1,0){4.00}}
\put(8.00,49.00){\line(0,1){4.00}}
\put(28.00,49.00){\line(0,1){4.00}}
\put(48.00,49.00){\line(0,1){4.00}}
\put(28.00,49.00){\line(-1,1){4.00}}
\put(48.00,49.00){\line(-1,1){4.00}}
\put(52.00,45.00){\line(-1,1){4.00}}
\put(52.00,25.00){\line(-1,1){4.00}}
\put(8.00,3.00){\makebox(0,0)[cc]{$w_{00}$}}
\put(28.00,3.00){\makebox(0,0)[cc]{$w_{10}$}}
\put(48.00,3.00){\makebox(0,0)[cc]{$w_{20}$}}
\put(2.00,29.00){\makebox(0,0)[cc]{$w_{01}$}}
\put(2.00,49.00){\makebox(0,0)[cc]{$w_{02}$}}
\put(24.00,25.00){\makebox(0,0)[cc]{$w_{11}$}}
\put(44.00,25.00){\makebox(0,0)[cc]{$w_{21}$}}
\put(24.00,45.00){\makebox(0,0)[cc]{$w_{12}$}}
\put(44.00,45.00){\makebox(0,0)[cc]{$w_{22}$}}
\put(18.00,4.00){\makebox(0,0)[cc]{$e_{00}$}}
\put(38.00,4.00){\makebox(0,0)[cc]{$e_{10}$}}
\put(18.00,15.00){\makebox(0,0)[cc]{$g_{00}$}}
\put(18.00,36.00){\makebox(0,0)[cc]{$g_{01}$}}
\put(38.00,15.00){\makebox(0,0)[cc]{$g_{10}$}}
\put(38.00,36.00){\makebox(0,0)[cc]{$g_{11}$}}
\put(4.00,19.00){\makebox(0,0)[cc]{$f_{00}$}}
\put(4.00,38.00){\makebox(0,0)[cc]{$f_{01}$}}
\put(51.00,19.00){\makebox(0,0)[cc]{$f_{20}$}}
\put(51.00,38.00){\makebox(0,0)[cc]{$f_{21}$}}
\put(18.00,52.00){\makebox(0,0)[cc]{$e_{02}$}}
\put(38.00,52.00){\makebox(0,0)[cc]{$e_{12}$}}
\put(31.00,19.00){\makebox(0,0)[cc]{$f_{10}$}}
\put(31.00,38.00){\makebox(0,0)[cc]{$f_{20}$}}
\put(15.00,27.00){\makebox(0,0)[cc]{$e_{01}$}}
\put(35.00,27.00){\makebox(0,0)[cc]{$e_{11}$}}
\end{picture}
\end{center}

Here the vertices are words $w_{ij}=x(yp)^iyq\by(r\by)^jz$ ($i,j\ge0$).
The positive edges have the form $e_{ij}=(1,x\to xyp,(yp)^iyq\by(r\by)^jz)$,
$f_{ij}=(x(yp)^iyq\by(r\by)^j,z\to r\by z,1)$ and
$g_{ij}=(x(yp)^iy,pyq\to q\by r,\by(r\by)^jz)$ ($i,j\ge0$). We choose the
maximal subtree ${\cal T}$ formed by the edges $e_{ij}$ for all
$i,j\ge0$ and also by the edges $f_{0j}$ for $j\ge0$. Thus our convention
that the choice of ${\cal T}$ makes the edges to be elements of the
fundamental group $\pi_1({\cal K},w)$, leads to the equalities $e_{ij}=1$
($i,j\ge0$), $f_{0j}=1$ ($j\ge0$) in the fundamental group of the complex.

Let us describe all $2$-cells of the complex ${\cal K}_w$. Recall that
there are defining relations of three types in $\qq$: $x=xyp$, $pyq=q\by r$,
$z=r\by z$. If we have two independent applications of elementary
transformations to words of the form $w_{ij}$, then it is easy to see that
they belong to different types because each of the letters $x$, $q$, $z$
occurs into the word $w_{ij}$ only once. Therefore, we have exactly three
situations.
\vspace{1ex}

1) The relations applied in the independent transformations are $x=xyp$ and
$z=r\by z$ (see the picture below).

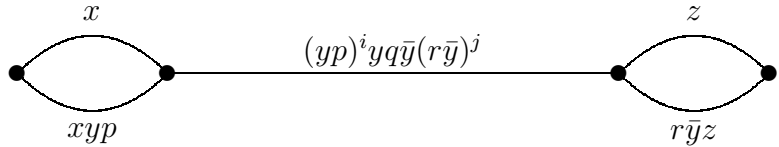
\begin{center} 
\unitlength=1mm
\special{em:linewidth 0.4pt}
\linethickness{0.4pt}
\begin{picture}(102.00,19.00)
\put(1.00,9.00){\circle*{2.00}}
\put(21.00,9.00){\circle*{2.00}}
\bezier{112}(1.00,9.00)(11.00,19.00)(21.00,9.00)
\bezier{112}(1.00,9.00)(11.00,-1.00)(21.00,9.00)
\put(81.00,9.00){\circle*{2.00}}
\put(101.00,9.00){\circle*{2.00}}
\bezier{112}(81.00,9.00)(91.00,19.00)(101.00,9.00)
\bezier{112}(81.00,9.00)(91.00,-1.00)(101.00,9.00)
\put(21.00,9.00){\line(1,0){60.00}}
\put(11.00,17.00){\makebox(0,0)[cc]{$x$}}
\put(91.00,17.00){\makebox(0,0)[cc]{$z$}}
\put(11.00,1.00){\makebox(0,0)[cc]{$xyp$}}
\put(91.00,1.00){\makebox(0,0)[cc]{$r\by z$}}
\put(51.00,12.00){\makebox(0,0)[cc]{$(yp)^iyq\by(r\by)^j$}}
\end{picture}
\end{center}

This diagram corresponds to the two paths in the Squier complex:
$e_{ij}f_{i+1,j}$ and $f_{ij}e_{i,j+1}$. This leads to relations
$e_{ij}f_{i+1,j}=f_{ij}e_{i,j+1}$. Simplifying, we have $f_{i+1,j}=f_{ij}$
for all $i,j\ge0$. It is obvious that $f_{ij}$ does not depend on $i$.
So $f_{0j}=1$ gives $f_{ij}=1$ for all $i,j\ge0$.
\vspace{1ex}

2) The relations are $x=xyp$ and $pyq=q\by r$ (see the picture below).

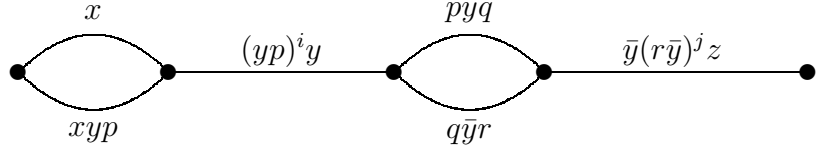
\begin{center} 
\unitlength=1mm
\special{em:linewidth 0.4pt}
\linethickness{0.4pt}
\begin{picture}(107.00,19.00)
\put(1.00,9.00){\circle*{2.00}}
\put(21.00,9.00){\circle*{2.00}}
\bezier{112}(1.00,9.00)(11.00,19.00)(21.00,9.00)
\bezier{112}(1.00,9.00)(11.00,-1.00)(21.00,9.00)
\put(51.00,9.00){\circle*{2.00}}
\put(71.00,9.00){\circle*{2.00}}
\bezier{112}(51.00,9.00)(61.00,19.00)(71.00,9.00)
\bezier{112}(51.00,9.00)(61.00,-1.00)(71.00,9.00)
\put(106.00,9.00){\circle*{2.00}}
\put(21.00,9.00){\line(1,0){30.00}}
\put(71.00,9.00){\line(1,0){35.00}}
\put(11.00,17.00){\makebox(0,0)[cc]{$x$}}
\put(11.00,1.00){\makebox(0,0)[cc]{$xyp$}}
\put(61.00,17.00){\makebox(0,0)[cc]{$pyq$}}
\put(61.00,1.00){\makebox(0,0)[cc]{$q\by r$}}
\put(36.00,12.00){\makebox(0,0)[cc]{$(yp)^iy$}}
\put(88.00,12.00){\makebox(0,0)[cc]{$\by(r\by)^jz$}}
\end{picture}
\end{center}

In this case, we have the equality $e_{i+1,j}g_{i+1,j}=g_{i,j}e_{i,j+1}$,
that is, $g_{i+1,j}=g_{ij}$ for all $i,j\ge0$. This means that $g_{ij}$
does not depend on $i$.
\vspace{1ex}

3) The relations are $pyq=q\by r$ and $z=r\by z$ (see the picture below).

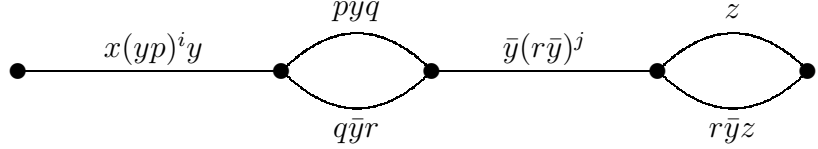
\begin{center} 
\unitlength=1.00mm
\special{em:linewidth 0.4pt}
\linethickness{0.4pt}
\begin{picture}(107.00,19.00)
\put(1.00,9.00){\circle*{2.00}}
\put(36.00,9.00){\circle*{2.00}}
\put(56.00,9.00){\circle*{2.00}}
\bezier{112}(36.00,9.00)(46.00,19.00)(56.00,9.00)
\bezier{112}(36.00,9.00)(46.00,-1.00)(56.00,9.00)
\put(86.00,9.00){\circle*{2.00}}
\put(106.00,9.00){\circle*{2.00}}
\bezier{112}(86.00,9.00)(96.00,19.00)(106.00,9.00)
\bezier{112}(86.00,9.00)(96.00,-1.00)(106.00,9.00)
\put(1.00,9.00){\line(1,0){35.00}}
\put(56.00,9.00){\line(1,0){30.00}}
\put(96.00,17.00){\makebox(0,0)[cc]{$z$}}
\put(96.00,1.00){\makebox(0,0)[cc]{$r\by z$}}
\put(46.00,17.00){\makebox(0,0)[cc]{$pyq$}}
\put(46.00,1.00){\makebox(0,0)[cc]{$q\by r$}}
\put(19.00,12.00){\makebox(0,0)[cc]{$x(yp)^iy$}}
\put(71.00,12.00){\makebox(0,0)[cc]{$\by(r\by)^j$}}
\end{picture}
\end{center}

Here we have the equality $g_{i,j}f_{i,j+1}=f_{i+1,j}g_{i,j+1}$. Hence
$g_{ij}=g_{i,j+1}$ for all $i,j\ge0$. Therefore, $g_{ij}$ depends on
neither $i$ nor $j$. For convenience, let $c=g_{ij}$ for all $i,j\ge0$.
\vspace{1ex}

Let us take an arbitrary vertex $v=w_{ij}=x(yp)^iyq\by(r\by)^jz$ for
some $i,j\ge0$. In the graph of groups that corresponds to the diagram
product, the product of $(i+1)$th power of $G$ and the $(j+1)$th power of
$H$ will correspond to the vertex $G_{w_{ij}}$ (the number of factors is just
the number of occurrences of $y$ and $\by$ in $v$, respectively). Thus
we can present the group $G_{w_{ij}}$ in the form
$$
K_{ij}=L_{iji}\times\cdots L_{ij0}\times H_{ij0}\times\cdots H_{ijj},
$$
where the factors of the form $L_{ijk}$ are isomorphic to $G$, and the
factors of the form $H_{ijk}$ are isomorphic to $H$. Relations~(\ref{descr})
that correspond to a positive edge $e=(u,s\to t,v)$ will be studied with
respect to the type of the involved defining relation $(s=t)\in\sss$ (there
are three types of them).
\vspace{1ex}

1) $s=x$, $t=xyp$. We have $e=(1,x\to xyp,v)$, where $v=(yp)^iyq\by(r\by)^jz$
for some $i,j\ge0$. The group $G_e$ is isomorphic to $G^{i+1}\times H^{j+1}$.
It maps isomorphically onto $K_{ij}$ under $\iota_e$. Note that $K_{i+1,j}$
is the direct product of $L_{i+1,j,i+1}$ and the isomorphic image of
$K_{ij}$ under $\tau_e$. Using the fact that $e=e_{ij}=1$ we see that
relation~(\ref{descr}) identifies some subgroups. Let us write down these
identifications as equalities. By these equalities we mean that the
corresponding elements of equal groups are identified.
We have: $L_{ijk}=L_{i+1,j,k}$ for $0\le k\le i$ and
$H_{ijk}=H_{i+1,j,k}$ for $0\le k\le j$.

2) $s=z$, $t=r\by z$. Now $e=(u,z\to r\by z,1)$, where
$u=x(yp)^iyq\by(r\by)^j$ for some $i,j\ge0$. Arguing analogously to the
previous case and taking into account that $e=f_{ij}=1$, we conclude that
relation~(\ref{descr}) leads to the following identifications of subgroups:
$L_{ijk}=L_{i,j+1,k}$ for $0\le k\le i$ and $H_{ijk}=H_{i,j+1,k}$ for
$0\le k\le j$.
\vspace{0.5ex}

Summarizing what we got in the first two cases of relations, we see that
groups $L_{ijk}$, $H_{ijk}$ depend of $k$ only. In other words, one can
introduce groups $L_k$, $H_k$ ($k\ge0$) in such a way that the equalities
$L_{ijk}=L_k$ for all $i\ge k$, $j\ge0$, and $H_{ijk}=H_k$ for all $i\ge0$,
$j\ge k$ hold in our diagram product.
\vspace{0.5ex}

3) $s=pyq$, $t=q\by r$. In this case $e=(u,pyq\to q\by r,v)$, where
$u=x(yp)^iy$, $v=\by(r\by)^j$ for some $i,j\ge0$. In the fundamental
group $\pi_1({\cal K},w)$, the equality $e=g_{ij}=c$ holds, as it was
shown above, the group $G_e$ is isomorphic to $G^{i+1}\times H^{j+1}$.
We have
$$
G_{\iota(e)}=
L_{i+1}\times\cdots\times L_0\times H_0\times\cdots\times H_j,
$$
$\iota_e$ maps $G^{i+1}$ onto the direct product
$L_{i+1}\times\cdots\times L_1$ and it maps $H^{j+1}$ onto the direct
product $H_0\times\cdots\times H_j$. Analogously,
$$
G_{\tau(e)}=
L_i\times\cdots\times L_0\times H_0\times\cdots\times H_{j+1},
$$
$\tau_e$ maps $G^{i+1}$ onto the direct product $L_i\times\cdots\times L_0$
and it maps $H^{j+1}$ onto the direct product $H_1\times\cdots\times H_{j+1}$.
Therefore, conjugating by the element $c$ takes $L_{i+1}$, \dots, $L_1$ to
$L_i$, \dots, $L_0$ respectively. The subgroups $H_0$, \dots, $H_j$ are
taken to $H_1$, \dots, $H_{j+1}$ respectively, under this conjugation.
Briefly, we can write $L_{i+1}^c=L_i$, $H_j^c=H_{j+1}$ for any $i,j\ge0$.

Thus the equalities $L_i=L_0^{c^{-i}}$, $H_j=H_0^{c^j}$ hold for any
nonnegative integers $i$, $j$. Let us extend these equalities to the case
of negative $i$, $j$ regarding these equalities as definitions. Note that
elements from different subgroups of the form $L_i$ ($i\ge0$) commute. So
the analogous fact is true for all integers $i$. The same fact is true
for subgroups $H_j$ for all $j\in\zz$. Let $G_i=L_{-i}$ ($i\in\zz$) by
definition. Obviously, $G_i^c=G_{i+1}$, $H_i^c=H_{i+1}$ for all $i\in\zz$.
This means that relations~(\ref{soot1}) hold. We also have conditions that
any element in $L_0$, $L_1$, \dots commutes with any element in $H_0$, $H_1$,
\dots\,. In particular, $[L_0,H_{j-i}]=1$ for any $j\ge i$. Taking into
account that $G_0=L_0$ and conjugating by the element $c^i$, we obtain
$[G_i,H_j]=1$ for $i\le j$, that is, relations~(\ref{soot2}) hold. It is
easy to see that these relations are in fact equivalent to the condition
that $[L_i,H_j]=1$ for any $i,j\ge0$. Indeed, the inequality $-i\le j$
and relations~(\ref{soot1}) allow us to conclude that $[G_{-i},H_j]=1$,
where $G_{-i}$ is $L_i$.

Thus we see that the diagram product we have calculated is in fact the
group given by relations~(\ref{soot1}) and~(\ref{soot2}), that is, it is
isomorphic to ${\cal O}(G,H)$.
}
\end{ex}

Using Theorem~\ref{DiagProd}, we have the following result.

\begin{thm}
\label{skew}
If $G$, $H$ are diagram groups, then ${\cal O}(G,H)$ is also a diagram
group.
\end{thm}

The previous example shows in details how, given two diagram groups $G$
and $H$, one can construct a presentation and a base, for which
${\cal O}(G,H)$ will be a diagram group.

\section{Nilpotent and Abelian Subgroups of Diagram Groups}
\label{NilAb}

We know from the previous section that soluble subgroups of any degree can
be subgroups of diagram groups. Contrary to that, we shall prove in this
Section that any nilpotent subgroup of a diagram group is abelian. We will
also establish the fact that all abelian subgroups of diagram groups are
free abelian. This will generalize the result that any abelian diagram group
is free abelian. Finally, we shall describe finite sets of pairwise
commuting diagrams, generalizing a description of pairs of commuting diagrams
from~\cite{GuSa97}.
\vspace{1ex}

We will use some concepts from combinatorics on diagrams
from~\cite[Section 15]{GuSa97}. For reader's convenience, let us recall
some definitions.

A spherical diagram is called {\em absolutely reduced} if any positive
integer power of it is reduced (does not contain dipoles). A spherical
diagram is called {\em normal} if it cannot be decomposed into a sum of
two non-spherical diagrams. We proved~\cite[Theorem  15.14]{GuSa97} that
for any spherical diagram $\Delta$ there exists an absolutely reduced
normal spherical diagram $\hat\Delta$ (that may have different base, in
general) and some (not necessarily spherical) diagram $\Psi$ such that
$\Delta=\Psi\iv\hat\Delta\Psi$.

\begin{thm}
\label{cen-der}
Let $H$ be an arbitrary subgroup of a diagram group $\dg$. Then the
centre of $H$ and the commutator subgroup of $H$ intersect trivially
that is, $Z(H)\cap H'=1$.
\end{thm}

\prf Let $G=\dg$ be a diagram group and let $H$ be a subgroup of $G$.
Suppose that $Z(H)\cap H'\ne1$. Consider a nontrivial element
$g\in Z(H)\cap H'$ and let $\Delta$ be a diagram representing it.
Applying~\cite[Lemma 15.10c]{GuSa97}, we find an absolutely reduced diagram
$\Delta_0$ that is conjugated to $\Delta$. Let $\Delta_0=\Psi\iv\Delta\Psi$,
where $\Psi$ is a $(\bs,\bs_0)$-diagram. Conjugation by $\Psi$ is an
isomorphism that takes the group $G$ to the group $G_0={\cal D}(\pp,\bs_0)$.
Under this isomorphism, the subgroup $H$ is taken to a subgroup $H_0$, where $g_0\in Z(H_0)\cap H_0'$,
and the element $g_0$ is represented by an absolutely reduced
$(w_0,w_0)$-diagram $\Delta_0$.

Let us decompose the diagram $\Delta_0$ into a sum of components:
$\Delta_0=A_1+\cdots+A_m$, where $A_i$ is a spherical $(w_i,w_i)$-diagram
($1\le i\le m$). As in~\cite[Theorem 15.35]{GuSa97}, we conclude that
the centralizer of $g_0$ is the direct sum of centralizers of the elements
represented by diagrams $A_1$, \dots, $A_m$. More precisely, if $\Gamma$ is a
spherical $(w_i,w_i)$-diagram that commutes with $\Delta_0$ in the group
$G_0$, then $\Gamma=B_1+\cdots+B_m$, where $B_i$ is a $(w_i,w_i)$-diagram
that commutes with $A_i$. By the assumption, any diagram representing an
element in $H_0$, commutes with $\Delta_0$ since $g_0$ belongs to the
centre of $H_0$.

Since $\Delta_0$ represents a nontrivial element, there exists an integer
$i$ between $1$ and $m$ such that the diagram $A_i$ is nontrivial and so it is
a simple absolutely reduced diagram. Its centralizer is cyclic
(see the proof of Theorem 15.35 in~\cite{GuSa97}). Let us now take two
diagrams $\Gamma$, $\Xi$ that represent elements in $H_0$. By the arguments
of the above paragraph, there are decompositions of the form
$\Gamma=B_1+\cdots+B_m$, $\Xi=C_1+\cdots+C_m$, where $B_i$, $C_i$ are
$(w_i,w_i)$-diagrams that commute with $A_i$. Cyclicity of the centralizer
of $A_i$ implies that $B_i$ and $C_i$ commute, that is, $[B_i,C_i]=\eps(w_i)$.
Therefore $[\Gamma,\Xi]=[B_1,C_1]+\cdots+[B_m,C_m]=
\Delta'+\eps(w_i)+\Delta''$, where $\Delta'$, $\Delta''$ are spherical
diagrams with bases $w_1\ldots w_{i-1}$, $w_{i+1}\ldots w_m$. It is clear
that the product of diagrams of the form $\Delta'+\eps(w_i)+\Delta''$
is again a diagram of the same form. Hence any element of the
commutator subgroup of the group $H_0$ has the form $\Delta'+\eps(w_i)+\Delta''$.
This contradicts the condition $\Delta_0=A_1+\cdots+A_m$, where
$A_i\ne\eps(w_i)$.

The Theorem is proved.

\begin{cy}
\label{nilp}
Any nilpotent subgroup of a diagram group is abelian.
\end{cy}

\prf Let $\dg$ be a diagram group and let $K$ be its nilpotent subgroup.
If $K$ is not abelian then $K$ has a (non-abelian) nilpotent subgroup $H$
of degree $2$. This means that the commutator subgroup of $H$ is contained
in its centre, that is, $H'\subseteq Z(H)$. Theorem~\ref{cen-der} claims
that the centre of $H$ and its commutator subgroup have trivial intersection
so $H'=1$, that is, $H$ is abelian,  a contradiction.
\vspace{2ex}

Let us now describe all abelian subgroups of diagram groups. It turns out
that all of them are free abelian.
Note that if there
was an abelian subgroup in a diagram group that was not free abelian, then
we would immediately disprove the Subgroup Conjecture, because it is known
from~\cite{GuSa97} that the quotient of any diagram group by its commutator
subgroup is free abelian.

\begin{thm}
\label{FrAb}
Any abelian subgroup of a diagram group is free abelian.
\end{thm}

\prf Let $\pp=\rav$ be a semigroup presentation, let $G=\dg$ be a
diagram group and let $H\le G$ be an abelian subgroup of $G$. If $H=1$
then we have nothing to prove. Let $H\ne1$. Consider an element $h\in H$,
$h\ne1$. Using certain conjugation and replacing the subgroup by an
isomorphic one, we can assume by ~\cite[Lemma 15.14]{GuSa97}
that $h$ is represented by an absolutely reduced normal diagram
$\Delta$ that is decomposed into the sum of components:
$\Delta=\Delta_1+\cdots+\Delta_m$, where $\Delta_i$ is a spherical diagram
with base $u_i$ ($1\le i\le m$). Since $H$ is abelian, it is contained
in the centralizer of $h$. Thus any element $g$ in $H$ decomposes into
a sum of $(u_i,u_i)$-diagrams. Denote by $\psi_i(g)$ the $i$th summand of
this decomposition. It is easy to see that $\psi_i$ is a homomorphism
of the group $H$ into the diagram group over $\pp$ with base $u_i$.

Let $1\le k\le m$ be a number such that $\Delta_k$ is nontrivial. The
centralizer of $\Delta_k$ is cyclic~\cite[Theorem 15.35]{GuSa97}.
Consider the homomorphism $\psi_k$. Firstly, $\psi_k(h)=\Delta_k\ne\ve(u_i)$.
Secondly, the image of $\psi_k$ is contained in the centralizer of $\Delta_k$,
that is, in a cyclic group. We thus proved that for any $h\in H$, $h\ne1$,
there exists a homomorphism $\psi\colon H\to\zz$ such that $\psi(h)\ne1$.
This means that $H$ is residually cyclic, that is, it embeds into a
Cartesian power of the infinite cyclic group. An easy argument in spirit
of linear algebra (using the Choice Axiom) shows that a Cartesian power
$\zz$ is a free abelian group. Therefore, $H$ is also free abelian.
\vspace{2ex}

In conclusion of this Section we give a simple but useful generalization
of Theorem 15.34 from~\cite{GuSa97}. We will need it later in
Section~\ref{Disto}.

\begin{thm}
\label{prws-cmmt}
Let $\pp=\rav$ be a semigroup presentation and let $G=\dg$ for some
$\bs\in\Sigma^+$. Suppose that $A_1$, \dots, $A_m$ are spherical diagrams
with base $\bs$ that pairwise commute in $G$. Then there exist a word
$v=v_1\ldots v_n$, spherical $(v_j,v_j)$-diagrams $\Delta_j$
$(1\le j\le n)$, integers $d_{ij}$ $(1\le i\le m$, $1\le j\le n)$ and
some $(\bs,v)$-diagram $\Gamma$ such that
$$
\Gamma\iv A_i\Gamma=\Delta_1^{d_{i1}}+\cdots+\Delta_n^{d_{in}}
$$
for all $1\le i\le m$. We can additionally assume that each of the diagrams
$\Delta_1$, \dots, $\Delta_n$ is either trivial or simple absolutely
reduced.
\end{thm}

\prf First of all let us show that the additional assumption about $\Delta_i$
can be proved provided the main statement of the theorem is proved.
If we have already found the
decompositions of diagrams from the main statement of the Theorem, then
by Lemma 15.14 from~\cite{GuSa97} we can find diagrams $\Psi_j$
such that diagarams $\Psi_j\Delta_j\Psi_j\iv$ are normal and absolutely
reduced. Each of these diagarams decomposes into a sum of components
that are either trivial or simple absolutely reduced. After that we apply
additional conjugation by $\Psi_1+\cdots+\Psi_n$ replacing each of the
$n$ summands by a sum of components.

The rest will be proved by induction on $m$. If $m=1$ or $m=2$ then it is
proved in~\cite{GuSa97} (Lemma 15.10c and Theorem 15.34). So we assume that
$m>2$, and the statement is true for all values less than $m$.
Consider the diagram $A_m$. Applying Lemma 15.10c, we find a
$(\bs,u)$-diagram $\Psi$ such that the diagram $A_m'=\Psi\iv A_m\Psi$ will
be absolutely reduced and normal. It can be decomposed into the sum of
components: $A_m'=B_1+\cdots+B_k$, where $B_j$ is a $(u_j,u_j)$-diagram
for some $u_j$ ($1\le j\le k$) and $u=u_1\ldots u_k$. Each of the diagrams
$B_j$ ($1\le j\le k$) is either trivial or simple and all of them are
absolutely reduced and normal. Let $A_i'=\Psi\iv A_i\Psi$ for all
$1\le i\le m-1$. It is clear that all diagrams $A_1'$, \dots, $A_m'$
pairwise commute. Since $A_1'$, \dots, $A_{m-1}'$ also commute with $A_m'$,
each of them can be decomposed into a sium of $(u_j,u_j)$-diagrams
($1\le j\le n)$ by~\cite[Theorem 15.35]{GuSa97}. We have
$A_i'=C_{i1}+\cdots+C_{in}$, where $1\le i\le m-1$, each of the diagrams
$C_{ij}$ is a $(u_j,u_j)$-diagram and $C_{ij}$ commutes with $B_j$
($1\le i\le m-1$, $1\le j\le n$). Suppose that $r$ ($1\le r\le k$) is a
number such that the component $B_r$ is nontrivial. Since $B_r$ is a simple
absolutely reduced diagram, it has a cyclic centralizer (see the proof
of Theorem 15.35). In particular, there exists a spherical diagram $\Delta_r$
with base $u_r$ such that any of the diagrams $C_{1r}$, \dots, $C_{m-1,r}$,
$B_r$ is a power of $\Delta_r$. Let $C_{ir}=\Delta_r^{d_{ir}}$
($1\le i\le m-1$), $B_r=\Delta_r^{d_{mr}}$. We do that with each $r$
($1\le r\le k$), for which the component $B_r$ is nontrivial.

Now let $1\le r\le k$ be such that the component $B_r$ is trivial,
that is, $B_r=\ve(u_r)$. Applying the inductive assumption to diagrams
$C_{1r}$, \dots, $C_{m-1,r}$, we find a word
$v_r=v_{1r}\ldots v_{n_r,r}$, spherical $(v_j,v_j)$-diagrams $\Delta_{jr}$
$(1\le j\le n_r)$, integers $d_{ijr}$ $(1\le i\le m-1$, $1\le j\le n_r)$
and some $(u_r,v_r)$-diagram $\Gamma_r$ such that
\begin{equation}
\label{decomp}
\Gamma_r\iv C_{ir}\Gamma_r=
\Delta_{1r}^{d_{i1r}}+\cdots+\Delta_{n_r,r}^{d_{i,n_r,r}}
\end{equation}
for all $1\le i\le m-1$. Now $\Gamma_r\iv B_r\Gamma_r=\ve(v_r)=
\ve(v_{1r})+\cdots+\ve(v_{n_r,r})$, and one can put $d_{m1r}=\cdots=
d_{m,n_r,r}=0$. Then equality~(\ref{decomp}) will be true also for
$i=m$, if we put $C_{mr}=B_r$.

For numbers $r$ such that $B_r$ is nontrivial we put $v_r=u_r$, and take the
trivial diagram for $\Gamma_r$. In this case, we also need to put $n_r=1$,
$\Delta_{1r}=\Delta_r$, $d_{i1r}=d_{ir}$. Then the equalities~(\ref{decomp})
are true for all $1\le r\le k$, $1\le i\le m$. Putting
$\Gamma=\Psi(\Gamma_1+\cdots+\Gamma_k)$, we see that
$$
\Gamma\iv A_i\Gamma=
\sum\limits_{r=1}^k\sum\limits_{j=1}^{n_r}\Delta_{jr}^{d_{ijr}},
$$
for $1\le i\le m$, that is, we have got the required decomposition of
diagrams into a sum.

The proof is complete.
\vspace{2ex}

In~\cite{GuSa97} we established that the conjugacy problem is decidable for
any diagram group $\dg$, where $\pp=\rav$ is a semigroup presentation
with decidable word problem. In particular, this implies the
decidability of the conjugacy problem in R.\,Thompson's group $F$. We can
pose a more general problem --- a uniform conjugacy problem for sequences.

\begin{prob}
\label{MCP}
Let $\pp=\rav$ be a semigroup presentation with decidable word problem,
$\bs\in\Sigma^+$, $G=\dg$. Does there exist an algorithm that decides,
given two sequences of elements $x_1$, \dots, $x_n$ and $y_1$, \dots, $y_n$
of the group $G$ $($elements are represented by diagrams$)$, whether there
is an element $z\in G$ such that $x_i^z=y_i$ for all $i$ from $1$ to $n$?
In particular, is this problem decidable for R.\,Thompson's group $F$?
\end{prob}

Note that there is some analogy between diagram groups and matrix groups
(we remarked about this in~\cite{GuSa97}). The corresponding question
for matrix groups was solved positively in~\cite{Sark} and independently
in~\cite{GrSeg}.

\section{Soluble Subgroups in Diagram Groups}
\label{Sol}

In this section we shed some light on the structure of soluble subgroups in
diagram groups. First of all let us consider an example that demonstrates
that there exist soluble subgroups of any degree in R.\,Thompson's group $F$.
Let $\pp=\{\,x\mid xx=x\}$. All groups ${\cal D}(\pp,x^k)$, where
$k=1,2,\ldots$\,, are isomorphic to $F$. Consider any nontrivial
$(x,x)$-diagram $\Delta$. Then the diagrams
$\Delta_1=\eps(x^2)+\Delta+\eps(x)$ and $\Delta_2=\eps(x)+\Delta+\eps(x^2)$
are conjugate because they are sums of components that conjugate
respectively. Secondly, the diagrams commute which can be seen directly and
are conjugated by the diagram $\Gamma=(x^2\to x)+\eps(x)+(x\to x^2)$.
Denoting by $a$, $b$ the elements of ${\cal D}(\pp,x^4)$ that represent
diagrams $\Delta_1$ and $\Gamma$, respectively, we can see that $a$ and $b$
generate the group $\zz\wr\zz$, a (restricted) wreath product of two infinite
cyclic groups. In the following theorem we present a more general form of the
above construction.

\begin{thm}
\label{xyz-wr}
Let $\pp=\rav$ be a semigroup presentation. Suppose that there exist
nonempty words $x$, $y$, $z$ over $\Sigma$ such that $xy=x$, $yz=z$ modulo
$\pp$, and suppose that the diagram group ${\cal D}(\pp,y)$ is nontrivial.
Then the group $G={\cal D}(\pp,xyz)$ contains a subgroup isomorphic to
$\zz\wr\zz$. Namely, let $\Delta$ be any nontrivial $(y,y)$-diagram, let
$\Gamma_1$ be arbitrary $(xy,x)$-diagram, and let $\Gamma_2$ be arbitrary
$(z,yz)$-diagram. Then elements $a$ and $b$, represented by diagrams
$\eps(x)+\Delta+\eps(z)$ and $\Gamma_1+\Gamma_2$, respectively, generate in
$G$ a subgroup isomorphic to $\zz\wr\zz$.
\end{thm}

\prf First of all, let us mention that elements $xz$, $xyz$, $xy^2z$, \dots\
are equal modulo $\pp$ so the diagram groups with these bases over $\pp$
will be isomorphic to each other. To get the diagrams from the above
example, one needs to put $x=y=z$ and then go from the group with base
$xyz=x^3$ to the group with base $xy^2z=x^4$ using conjugation by the
element $(x\to x^2)+\eps(x^2)$.

To show that elements $a$, $b$ of some group $G$ generate $\zz\wr\zz$, it
suffices to show that the elements $a_i=a^{b^i}$ ($i\in\zz$) form a free
basis of a free abelian group. To check this, it suffices to show that
for any positive integer $n$, the elements $a_0$, $a_1$, \dots, $a_n$
form a basis of the free abelian group they generate. We will
explicitly find the elements $a_i$ ($0\le i\le n$) of the corresponding
diagram group ${\cal D}(\pp,xyz)$. For convenience, we will go to the
diagram group over $\pp$ with base $xy^{n+1}z$ using conjugation by
the diagram
$$
\Psi=
(\Gamma_1\iv+\ve(yz))(\Gamma_1\iv+\ve(y^2z))\cdots(\Gamma_1\iv+\ve(y^nz)).
$$
One can check directly that $c_i=\Psi\iv a_i\Psi=
\ve(xy^{n-i})+\Delta+\ve(y^iz)$ for all $0\le i\le n$. It is clear that
the elements $c_0$, $c_1$, \dots, $c_n$ pairwise commute. The obvious formula
$$
c_0^{d_0}c_1^{d_1}\ldots c_n^{d_n}=
\ve(x)+\Delta^{d_0}+\Delta^{d_1}+\cdots+\Delta^{d_n}+\ve(z)
$$
shows that the elements $c_0$, $c_1$, \dots, $c_n$ form a basis of the free
abelian subgroup in ${\cal D}(\pp,xy^{n+1}z)$. So the elements $a_0$,
$a_1$, \dots, $a_n$ also form a basis of a free abelian subgroup of
${\cal D}(\pp,xyz)$ as desired.
\vspace{2ex}

We will return to the group $\zz\wr\zz$ later. Now we shall prove a
simple fact that will imply that $F$ contains soluble subgroups of any
degree.

\begin{lm}
\label{FwrZ}
$($Restricted$)$ wreath product $F\wr\zz$ is a subgroup of $F$.
\end{lm}

\prf We will use some known properties of R.\,Thompson's group $F$
mentioned in Section~\ref{Prelim}. As we mentioned above, there are several
representations of $F$ by piecewise linear functions. Let us consider the
representation by functions on $[0,\infty)$. For any positive integer $k$ we
consider the functions from $F$ that are identical outside $[k,k+1]$. By
$\Phi_k$ we denote the set of all these functions. It is obvious that they
form a group isomorphic to the group of all piecewise linear functions on
$[0,1]$ (with the properties mentioned in Section~\ref{Prelim}), that is, it
is isomorphic to $F$. It is also easy to see that elements in different
subgroups $\Phi_k$ commute with each other. Therefore, the groups $\Phi_k$
($k\ge1$) generate a direct power of the group $F$. Conjugation by the
element $x_0$, represented by the function given by $tx_0=2t$ ($t\in[0,1]$),
$tx_0=t+1$ ($t\ge1$), takes $\Phi_k$ to $\Phi_{k+1}$. It is now clear
that $t$ and $\Phi_k$ ($k\ge1$) generate the restricted wreath product
$F\wr\zz$ in $F$.
\vspace{0.5ex}

It is not hard to find generators of the subgroup $F\wr\zz$ of $F$ in
a diagram form and also in a normal form. In particular, the subgroup in $F$
generated by elements $x_0$, $x_1x_2x_1^{-2}$, $x_1^2x_2x_1^{-3}$ will be
isomorphic to $F\wr\zz$. The reader can easily draw the diagrams representing
these elements.
\vspace{1ex}

Let us define a sequence of groups by induction: $H_1=\zz$,
$H_{n+1}=H_n\wr\zz$. Thus groups $H_n=(\cdots(\zz\wr\zz)\wr\cdots\,)\wr\zz$,
where $\zz$ occurs $n$ times, are diagram groups by Theorem~\ref{wr}.
The group $H_n$ is soluble of degree $n$. Using Lemma~\ref{FwrZ} and
elementary properties of wreath products, we have the following result that
can be proved by induction on $n$. This statement was obtained by M.\,Brin
(private communication), see also~\cite{Brin}.

\begin{cy}
\label{ZwrZetc}
For any $n$, the group $H_n=(\cdots(\zz\wr\zz)\wr\cdots\,)\wr\zz$ is a
soluble subgroup of degree $n$ in R.\,Thompson's group $F$.
\end{cy}

Talking about wreath products, we would like to mention a fact about
subgroups of R.\,Thompson's group $F$. It was shown in~\cite{BrSq85} that
any subgroup of $F$ is either metabelian or contains an infinite direct
power of the group $\zz$. (In fact, one can replace the word ``metabelian"
by ``abelian", see~\cite{CFP}.) The proof given in~\cite{BrSq85,CFP} uses
representations of $F$ by piecewise linear functions. Actually, the result
is obtained for subgroups of some group which is bigger than $F$. It turns
 out that one can extract a stronger fact from this proof. Consider all
piecewise linear continuous transformations of the unit interval $I=[0,1]$
onto itself. We consider only mappings that preserve orientation and have
finitely many breaks of the derivative. All these functions form a group with
respect to composition. Let us denote this group by $PL_0(I)$. It contains
$F$ as a subgroup. We have the following alternative for subgroups of the
group $PL_0(I)$.

\begin{thm}
\label{PLIZwrZ}
Any subgroup of $PL_0(I)$ is either abelian, or contains an isomorphic
copy of $\zz\wr\zz$.
\end{thm}

\prf Our proof basicaly follows the proof or a weaker alternative
from~\cite{BrSq85,CFP}. For $f\in PL_0(I)$ by $\supp f$ we denote the set of
all $t\in I$, for which $tf\ne t$. Let $G$ be a non-abelian subgroup of
$PL_0(I)$. Consider functions $f,g\in G$ such that $fg\ne gf$. Let
$J=\supp f\,\cup\,\supp g$. It is obvious that $J$ is a union of finitely
many disjoint intervals $J_k=(a_k,b_k)$, $1\le k\le m$. By definition,
$[f,g]\ne1$ in $PL_0(I)$. Then on some of intervals $J_1$, \dots, $J_m$ our
function $[f,g]$ is not the identity. Denote by $\nu(f,g)$ the number of such
intervals. Without loss of generality, we can assume that the elements
$f,g\in G$ which do not commute are chosen in such a way that the value
$\nu(f,g)$ is the smallest possible. Let $H$ be a subgroup of $PL_0(I)$
generated by $f$ and $g$. By definition, the endpoints of $J_1$, \dots, $J_m$
are stable under $f$ and $g$ so each of these intervals is $H$-invariant.

An easy argument shows that for any $x,y\in J_k$ ($1\le k\le m$),
where $x<y$, there exists a function $w\in H$ such that $xw>y$. Let us take
the greatest upper bound $z$ of the set $\{\,xh\mid h\in H\,\}$. It is clear
that $a_i<z\le b_i$. If $z\ne b_i$ then either $zf\ne z$ or $zg\ne z$ by
definition of the set $J$. Without loss of generality let $zf\ne z$.
This inequality also holds in a small neighbourhood of the point $z$.
Therefore one of the numbers $zf$ or $zf\iv$ is greater than $z$, a
contradiction. Thus $z=b_i$. This implies that acting by some element of $H$
one can make the image of $x$ as close to $b_i$ as one wishes which is what
we had to prove.

Let us take an interval $(a_i,b_i)$ ($1\le i\le m$) such that $[f,g]$ is
not identical on it. It is easy to ee that the function $[f,g]$ is identical
in some neighbourhood of each of the points $a_i$, $b_i$. Thus $\supp[f,g]$
is nonempty and it is contained in $[c_0,d_0]$, where $a_i<c_0<d_0<b$.
According to the above, there exists a function $w\in H$ such that
$d_0<c_0w<b$. Let us denote $[f,g]$ by $h_0$. For any $n\ge1$, let
$c_n=c_0w^n$, $d_n=d_0w^n$, $h_n=h_0^{w^n}$. It is obvious that
$c_0<d_0<c_1<d_1<\cdots\,$, and $\supp h_n\cap J_i\subseteq[c_n,d_n]$.
Therefore, for any $i,j\ge0$, the commutator $[h_i,h_j]$ is identical on
$J_i$. Suppose that $[h_i,h_j]\ne1$ for some $i$, $j$. Since all the
intervals $J_1$,\dots, $J_m$ are $H$-invariant, it is clear that
all the functions $h_1$, $h_2$, \dots\ are identical on all the intervals
$J_k$ ($1\le k\le m$) where the function $h_0$ is identical. Therefore,
we can replace $f$, $g$ by $h_i$, $h_j$ obtaining $\nu(h_i,h_j)<\nu(f,g)$.
This contradicts the choice of $f$, $g$. This proves that $[h_i,h_j]=1$ for
any $i,j\ge0$.

So far we very closely followed the proof given in~\cite{CFP}. The conclusion
in~\cite{CFP}, is that elements $h_0$, $h_1$, $h_2$, \dots\ form a basis of
a free abelian group. To prove a stronger statement of our theorem, it
remains now to add that $h_n^w=h_{n+1}$ for all $n\ge0$, so the elements
$h_0$, $w$ generate $\zz\wr\zz\subseteq G$.

The Theorem is proved.
\vspace{1ex}

Using Theorem~\ref{FrAb}, we have the following alternative for subgroups
of R.\,Thompson's group $F$.

\begin{cy}
\label{Alter}
Any subgroup of R.\,Thompson's group $F$ is either free abelian or
contains the restricted wreath product $\zz\wr\zz$.
\end{cy}

Note that in the group of all piecewise linear functions, not every abelian
subgroup is free abelian.
\vspace{0.5ex}

We can extract one more corollary from Theorem~\ref{PLIZwrZ}.

\begin{cy}
\label{OneRelPL}
A non-abelian group with one defining relation cannot be a subgroup of
the group $PL_0(I)$ $($in particular, it cannot be a subgroup of
R.\,Thompson's group $F)$.
\end{cy}

\prf In~\cite{Cheb}, A.\,A.\,Chebotar described all subgroups
of one-relator groups that do not contain free subgroups of rank $2$.
They are: a) abelian subgroups, b) free product $\zz_2\ast\zz_2$ of the
cyclic group of order $2$ by itself, and c) Baumslag -- Solitar groups
$B_{1k}=\pres{a,b}{b\iv ab=a^k}$. It is clear that the group $\zz\wr\zz$
does not occur in this list. Thus a non-abelian one-relator group cannot
be a subgroup of $PL_0(I)$ by Theorem~\ref{PLIZwrZ}.

The Corollary is proved.
\vspace{2ex}

Now let us consider the following interesting question: under what condition
on a semigroup presentation $\pp$ a diagram group over this presentation
contains $zz\wr\zz$ as a subgroup? The answer is given in the following
theorem.

\begin{thm}
\label{3equiv}
Let $\pp=\rav$ be a semigroup presentation and let $G=\dg$ be a diagram
group. Then the following three conditions are equivalent.
\begin{enumerate}
\item
The group $G=\dg$ contains $\zz\wr\zz$ as a subgroup.
\item
The group $G$ contains elements $a$, $b$ such that $[a,b]\ne1$, $[a,a^b]=1$
$($in other words, there are two distinct elements in $G$ that are conjugate
and commute\/$)$.
\item
There are words $x$, $y$, $z$ such that the equalities $xy=x$, $yz=z$,
$xz=\bs$ hold modulo $\pp$, and ${\cal D}(\pp,y)\ne1$.
\end{enumerate}
\end{thm}

\prf The proof uses the following scheme:
$(1)\Rightarrow(2)\Rightarrow(3)\Rightarrow(1)$. The implication
$(1)\Rightarrow(2)$ is obvious and holds for any group $G$. The implication
$(3)\Rightarrow(1)$ was proved in Theorem~\ref{xyz-wr}. It remains to show
that $(2)\Rightarrow(3)$.

Suppose that $G=\dg$ has elements $a$, $b$ such that $[a,b]\ne1$,
$[a,a^b]=1$. By Theorem~\ref{prws-cmmt}, one can pass from the base $\bs$ to
some base $v$ that equals $\bs$ modulo $\pp$ in such a way that the
diagarams representing the two given commuting elements will be absolutely
reduced and normal. Without loss of generality we can assume that $a$
is represented by a diagram $A=A_1+\cdots+A_m$ decomposed into the sum of
components, and the element $a^b$ is represented by a diagarm $C$ that
has a decomposition into the sum of the same number of components:
$C=C_1+\cdots+C_m$ by~\cite[Lemma 15.15]{GuSa97}. By the same lemma,
the element $b$ which conjugates $A$ and $C$ is represented by a diagram $B$
of the form $B_1+\cdots+B_m$, where $C_i=B_i\iv A_iB_i$. Let
$v=v_1\ldots v_m=v_1'\ldots v_m'$, where $A_i$, $C_i$ are spherical
diagrams with bases $v_i$, $v_i'$, respectively, and let $B_i$ be a
$(v_i,v_i')$-diagram ($1\le i\le m$). It is obvious that $v_i=v_i'$ modulo
$\pp$ for all $i$.

We would like to prove that there exists an $i$ from $1$ to $m$ such that
the diagram $A_i$ (and $C_i$ as well) is nontrivial and the occurrences
of the words $v_i$, $v_i'$ in the word $v$ do not overlap and do not
contain each other. This would imply condition $(3)$. Indeed, without loss
of generality, let $v=pv_iqv_i'r$. Then the equalities
$p=v_1\ldots v_{i-1}=v_1'\ldots v_{i-1}'=pv_iq$,
$r=v_{i+1}'\ldots v_m'=v_{i+1}\ldots v_m=qv_i'r=qv_ir$ hold modulo $\pp$.
One can put $x=p$, $y=v_iq$, $z=v_ir$, and then the equalities $x=p=pv_iq=xy$,
$z=v_ir=v_iqv_ir=yz$, $\bs=v=pv_iqv_i'r=pv_iqv_ir=pv_ir=xz$ will hold modulo
$\pp$ (that is, in the semigroup $S$). The diagram group over $\pp$
with base $y=v_iq$ will be definitely nontrivial because there exists a
nontrivial spherical diagram $A_i+\ve(q)$ with this base.

Let us prove the existence of an $i$ such that $A_i$ is nontrivial, and
the occurrences of the word $v_i$, $v_i'$ in the word $v$ have no common
letters. Let us consider an arbitrary $1\le i\le m$ such that $A_i$ is
nontrivial. Since $C$ commutes with $A$, the diagram $C$ can be presented as
a sum $C'+D+C''$, where $C'$, $D$, $C''$ are spherical diagrams with bases
$v_1\ldots v_{i-1}$, $v_i$, $v_{i+1}\ldots v_m$, respectively, and $D$
commutes with $A_i$. Therefore, the diagram $D$ consists of one component.
The same is true for the diagram $C_i$. So, if the occurrences of $v_i$,
$v_i'$ have common letters, then diagrams $D$ and $C_i$ must coincide.
This implies that $A_i$, $C_i$ are powers of the same element $\Delta_i$,
and are conjugate. It follows from results of~\cite[Section 15]{GuSa97}
that $A_i=C_i$, and the occurences of $v_i$ and $v_i'$ coincide. It is
now obvious that $A=C$. But this contradicts the assumption that $a\ne a^b$.

The proof is complete.
\vspace{1ex}

Let us make two remarks about the theorem we have just proved. Firstly,
in the third condition we cannot avoid the condition that the
diagram over $\pp$ with base $y$ is nontrivial. Without this condition, the
diagram group may coincide with $\zz$. Note that the algorithm to verify
whether the diagram group over a given finite presentation with given base
is nontrivial, is unknown. Secondly, we have to note that if elements $a$,
$b$ of a diagram group are such that $[a,b]\ne1$, $[a,a^b]=1$, then the
subgroup isomorphic to $\zz\wr\zz$, is not necessarily contained in the
subgroup generated by $a$, $b$. An example illustrating that is given below
in Section~\ref{SubConj}.
\vspace{2ex}

Finishing this Section, let us give a sufficient condition for a
diagram group to contain R.\,Thompson's group $F$ as a subgroup.

\begin{thm}
\label{FSub}
Let $\pp=\rav$ be a semigroup presentation, and let the semigroup $S$
presented by $\pp$ contain an idempotent. Then there is a word $\bs$ such
that the diagram group $G=\dg$ contains R.\,Thompson's group $F$ as a
subgroup. Moreover, one can take any word $\bs$ that represents an idempotent
in $S$ for such a base.
\end{thm}

\prf The proof is quite easy. It is based on the fact that all proper
homomorphic images of the group $F$ are abelian (see~\cite{CFP}).
Let us take a word $\bs$ such that $\bs\bs=\bs$ modulo $\pp$. Let us
consider a reduced $(\bs^2,\bs)$-diagram $\Delta$ over $\pp$. Now we
construct a homomorphism from $F$ to $\dg$ in the following way.
We use the fact that $F$ is a diagram group over $\qq=\pres{x}{x^2=x}$.
It is convenient to take the element $x^5$ as a base. To any diagram over
$\qq$, we assign a diagram over $\pp$, replacing the label $x$ by $\bs$
and filling in the cells of the form $x^2=x$ by diagrams $\Delta$. This rule
defines a homomorphism from $F$ to $\dg$. Taking into account what we have
said above, it is enough to check that the image of this homomorphism
is not abelian. To show this, it suffices to compute the image of the
commutator $[x_0,x_1]$. In the group ${\cal D}(\qq,x^5)$, this commutator
is represented by the diagram $\ve(x)+(x^2=x)+(x=x^2)+\ve(x)$. It is obvious
that after we replace all $(x^2,x)$-cells in this diagram by copies of
$\Delta$, we get a diagram without dipoles. Thus the image of this commutator
is not equal to the identity, so the image of $F$ under the homomorphism is
not abelian.

The Theorem is proved.
\vspace{2ex}

We do not know whether the condition on $S$ to have an idempotent is also
sufficient.

\begin{prob}
Let $\pp=\rav$ be a semigroup presentation and let $G=\dg$ be a diagram
group. Suppose that $G$ contains R.\,Thompson's group $F$ as a subgroup.
Is it true that the semigroup $S$ presented by $\pp$ contains an
idempotent?
\end{prob}

\section{The Subgroup Conjecture}
\label{SubConj}

In this Section, we construct a counterexample to the Subgroup
Conjecture, that is, we will construct a subgroup in a diagram
group that is not a diagram group itself. Note that the first
candidate to disprove the Subgroup Conjecture was the group $F'$
--- the commutator subgroup of R.\,Thompson's group $F$. However,
it turned out that $F'$ is a diagram group. This answered
some open questions about diagram groups. Before proving the
corresponding Theorem, let us make the following two remarks.
\vspace{1ex}

The first remark is about properties of semigroup diagrams over certain
special presentations. Let we have a semigroup presentation of the form
$\pres{X}{\rr}$, where all relations in $\rr$ have the form $u=V$, where
$u\in X$, $V\in X^+$. We assume that all left-hand sides of the relations
are distinct and the right-hand sides contain more than
one letter. It is known (see, for instance~\cite{Gu97}) that any
reduced diagram $\Delta$ over such a presentation can be
uniquely decomposed into a concatenation:
$\Delta=\Delta_1\circ\Delta_2$, where $\Delta_1$ corresponds to a derivation
where only applications of relations of the form $u=V$ from $\rr$ are used,
and $\Delta_2$ corresponds to a derivation where only applications of
relations of the form $V=u$ are used, $(u=V)\in\rr$. This fact can be easily
proved by choosing the longest positive path from $\iota(\Delta)$ to
$\tau(\Delta)$. It is easy to see that all cells ``above" this path will
correspond to relations $u=V$, and all cells ``below" the path will
correspond to $V=u$. We will call $\Delta_1$ (resp. $\Delta_2$) the
{\em positive\/} (resp. {\em negative\/}) part of diagram $\Delta$.

Note that the presentation $\pres{x}{x=x^2}$ satisfies the above conditions.
The same holds for the presentation below from the statement of
Theorem~\ref{CommF}.
\vspace{1ex}

The second remark is about the structure of a commutator subgroup of
a diagram group. It is described in~\cite[Theorem 11.3]{GuSa97}. Let
us recall this description. Let $\Delta$ be a $(w,w)$-diagram over
$\pp=\rav$. By $M$, we denote the monoid presented by $\pp$. We consider
the free abelian group ${\cal A}$ with $M\times\rr\times M$ as a free
basis. For each vertex $\mu$ of diagram $\Delta$ we take any positive path
from $\iota(\Delta)$ to $\mu$. Its label defines an element in
$\ell(\mu)\in M$. It is easy to show that this element does not depend
on the choice of a path. Analogously, we define the element $r(\mu)\in M$
as the value of the label of any positive path from $\mu$ to $\tau(\Delta)$.
Now to each cell $\pi$ of the diagram $\Delta$ we assign an element
$\delta\cdot(\ell(\iota(\pi)),u=v,r(\tau(\pi)))$, where $\delta=1$, if
$u=\varphi(\topp(\pi))$, $v=\varphi(\bott(\pi))$, $(u,v)\in\rr$ and
$\delta=-1$, if $v=\varphi(\topp(\pi))$, $u=\varphi(\bott(\pi))$,
$(u,v)\in\rr$. By $\rho(\Delta)$ we denote the sum of elements assigned
to all the cells of diagram $\Delta$. Thus $\rho$ defines a homomorphism
from the group $G=\dg$ into ${\cal A}$. As it is shown in~\cite{GuSa97},
the kernel of $\rho$ is exactly $G'$ --- the commutator subgroup of $G$.

\begin{thm}
\label{CommF}
The commutator subgroup of R.\,Thomp\-son's group $F$ is a diagram
group. Namely, $F'\cong{\cal D}(\qq,a_0b_0)$, where
$$
\qq=\pres{x,a_i,b_i\ (i\ge0)}{x=xx,a_i=a_{i+1}x,b_i=xb_{i+1}\ (i\ge0)}.
$$
\end{thm}

\prf Here is the direct proof of this proposition. Let us construct a mapping
from $H={\cal D}(\qq,a_0b_0)$ to ${\cal D}(\pp,x^2)$, where
$\pp=\pres{x}{x=x^2}$. To each spherical diagram over $\qq$ with base
$a_0b_0$ we assign a diagram that is obtained from the previous one
replacing all its labels by $x$. It is clear that we get a spherical
diagram over $\pp$ with base $x^2$. Obviously, this induces a
homomorphism $\psi\colon H\to F$ because ${\cal D}(\pp,x^2)\cong F$. Our
aim is to prove that the homomorphism $\psi$ is injective and its image
is $F'$.

Let $\Delta$ be a nontrivial reduced $(a_0b_0,a_0b_0)$-diagram over
$\qq$. Its image under $\psi$ cannot contain dipoles. Otherwise, the
preimages of the cells that form a dipole in $\Delta$, would form a
dipole themselves. This implies that $\psi$ is injective.

Let us check that $\psi(\Delta)\in F'$ for any reduced diagram $\Delta$ in
$H$. The monoid $M$ presented
by $\pres{x}{x=xx}$ consists of two elements $1$ and $x$. A cell
$\pi$ of diagram $\Delta'=\psi(\Delta)$ satisfies $\ell(\iota(\pi))=1$
if and only if $\iota(\pi)=\iota(\Delta')$. Analogously, $r(\tau(\pi))=1$
if and only if $\tau(\pi)=\tau(\Delta')$. As we know,
the diagram $\Delta'$ is reduced. It can be decomposed into a product
$\Delta_1'\circ\Delta_2'$ of its positive and negative part according to
the first remark before the statement of the Theorem. It is easy to see
that there are no cells $\pi$ of the diagram $\Delta'$ can satisfy both
of the conditions $\iota(\Delta)=\iota(\pi)$, $\tau(\Delta)=\tau(\pi)$
simultaneously (recall that the base of $\Delta'$ is $x^2$).

Consider the diagram $\Delta$ and decompose it into a concatenation of
positive and negative part: $\Delta=\Delta_1\circ\Delta_2$ (this is possible
because $\qq$ satisfies conditions of the first remark before the statement
of the theorem). Let $a_n$ ($n\ge0$) be the label of the first edge of
the path that cuts $\Delta$ into a positive part and a negative part.
Then it is easy to extract from the form of the defining relations that all
labels of edges which start at $\iota(\Delta)$, if one reads them from the
top to the bottom of the diagram, are $a_0$, $a_1$, \dots, $a_n$, \dots,
$a_1$, $a_0$. From this, it follows that the number of cells $\pi$ such that
$\iota(\pi)=\iota(\Delta)$, is the same for $\Delta_1$ and $\Delta_2$.
>From this fact, we immediately conclude that it is the same for $\Delta'$,
if we compare the number of these cells in $\Delta_1'$ and $\Delta_2'$.
All these cells from $\Delta_1'$ map to $(1,x=x^2,x)$ under $\rho$
and all cells from $\Delta_2'$ map to $-(1,x=x^2,x)$ (we emphasize the fact
that we are talking about the cells whose initial vertices coincide with
the initial vertex of the diagram).

An analogous argument can be applied to the cells whose terminal
vertices coincide with the terminal vertex the diagram. Here the list of
labels of the edges that come into $\tau(\Delta)$, if one reads them from top
to bottom, is $b_0$, $b_1$, \dots, $b_m$, \dots, $b_1$, $b_0$ for some
$m\ge0$. Now one can use the fact that the cells of $\Delta_1'$
which we deal with in this paragraph map to $(x,x=x^2,1)$ and the cells of
$\Delta_2'$ map to $-(x,x=x^2,1)$. It remains to note that $\Delta_1'$ has
the same number of cells as $\Delta_2'$ since $\Delta'$ is spherical and all
relations have the same form. Therefore, the number of cells $\pi$ of
$\Delta_1'$ that satisfy $\iota(\pi)\ne\iota(\Delta')$ and
$\tau(\pi)\ne\tau(\Delta')$  is the same as the number of corresponding
cells in $\Delta_2'$. However, the first ones map to $(x,x=x^2,x)$ and
the second ones map to $-(x,x=x^2,x)$. Hence $\rho(\Delta')=0$. This
proves that $\psi(\Delta)\in F'$.

It remains to show that every element in $F'$ belongs to the image of $\psi$.
In order to do that, let us take a reduced spherical diagram $\Delta'$ with
base $x^2$ over $\pres{x}{x=x^2}$. It follows from $\Delta'\in F'$ that
$\rho(\Delta')=0$. We consider separately the cells of three types:
those that map into\ \ a) $\pm(1,x=x^2,x)$,\ \ b) $\pm(x,x=x^2,1)$,\ \
c) $\pm(x,x=x^2,x)$ under $\rho$, respectively.

Each cell belongs to exactly one of the three types. So it is clear
that the sum over all cells of each of the types equals zero. This
means that in the decomposition $\Delta'=\Delta_1'\circ\Delta_2'$ into
positive and negative part, the number of cells of each of the types
in $\Delta_1'$ will be the same as the number of cells of the same type in
$\Delta_2'$. Let us have $n\ge0$ cells of the first type in each of the
parts. Let us rename labels of the edges that go out of $\iota(\Delta')$,
replacing them by $a_0$, $a_1$, \dots, $a_n$, \dots, $a_1$, $a_0$,
respectively, from top to bottom. Analogously, let us have $m\ge0$ cells
of the second type in each of the parts. We rename labels of the edges
that come into $\tau(\Delta')$, replacing them in the same way by
$b_0$, $b_1$, \dots, $b_m$, \dots, $b_1$, $b_0$, respectively. The diagram
we get as a result will be denoted by $\Delta$.

It remains to note that $\Delta$ will be a spherical diagram over $\qq$
with base $a_0b_0$. Indeed, any cell that has the same initial vertex as
the one of $\Delta$, corresponds to a relation of the form $a_i=a_{i+1}x$
($i\ge0$) or its inverse. If a cell $\pi$ of the positive part is taken,
then $\topp(\pi)=e$, $\bott(\pi)=e_1e_2$, where $e$, $e_1$, $e_2$ are
edges of the diagram. By our construction, the label of $e$ equals $a_i$
for some $i\ge0$. It follows from the way we renamed the labels that
$e_1$ has label $a_{i+1}$. Note that the initial vertex of $e_2$ is not
$\iota(\Delta)$ because $e_1$ cannot be a loop. Also the terminal vertex
of $e_2$ is not $\tau(\Delta)$. Otherwise the edge $e$ connects the initial
and the terminal vertex of $\Delta$ but this is impossible. Therefore,
the label of $e_2$ is $x$. The arguments for the negative part of the
diagram are analogous. Of course, any cell that has the same terminal vertex
as $\Delta$, corresponds to a relation of the form $b_i=xb_{i+1}$
($i\ge0$) or its inverse. It is clear that $\psi$ takes $\Delta$ into
$\Delta'$. This completes the proof.

\begin{cy}
\label{simp}
A diagram group can be simple. In particular, there exist nontrivial
diagram groups that coincide with their commutator subgroups and so
they do not admit an LOG-presentation.
\end{cy}

The group $F'$ is simple (see~\cite{CFP}). We proved in~\ref{CommF}
that $F'$ is a diagram group. In~\cite[Section 17]{GuSa97} we asked if
a nontrivial diagram group may coincide with its commutator subgroup.
We have given a positive answer. This is interesting if to compare
this result with~\cite[Theorem 12.1]{GuSa97}. It was proved there that
if all diagram groups over a semigroup presentation coincide with their
commutator subgroups, then all of them are trivial. As we see, certain
diagram groups may coincide with their commutator subgroups. As a
by-product, we gave an answer to Problem~17.1 of the same paper: is it
true that any diagram group admits an LOG-presentation? Recall that
an LOG-presentation is a group presentation such that all defining relations
have form $a=b^c$, where $a$, $b$ and $c$ are generators. The groups
that admit such a presentation are called LOG-groups (this concept was
introduced in~\cite{Bogl87}, where these groups were characterized in
terms of labelled oriented graphs). In Russian papers, one can often
meet an equivalent terminology ``C-group". Some interesting characterization
of these groups was recently obtained by Yu.\,V.\,Kuzmin~\cite{Kuz95,Kuz96}.
We have already shown that any diagram group over a complete presentation
(see~\cite{GuSa97}) admits an LOG-presentation (cf R.\,Thompson's group $F$).
We also proved that any diagram group is a retract of an LOG-group.
Since any LOG-group has $\zz$ as its homomorphic image, it cannot coincide
with its commutator subgroup. Thus we proved that a diagram group may
not have an LOG-presentation.
\vspace{2ex}

To construct a counterexample to the Subgroup Conjecture, we strongly use
results of Section~\ref{Constr}. In particular, we need Theorem~\ref{skew}
and Example~\ref{kh}.

\begin{thm}
\label{dispr}
There exist subgroups of digaram groups that are not diagram groups
themselves. For instance, the following one-relator group
$$
\pres{x,y}{xy^2x=yx^2y}
$$
can be isomorphically represented by diagrams over a semigroup
presentation but it is not a diagram group itself.
\end{thm}

\prf Let $L=\pres{x,y}{xy^2x=yx^2y}$. We shall prove that $L$ is not a
diagram group. Consider the group $L_\infty$ given by
$$
L_\infty=\pres{z_i\ (i\in\zz)}{[z_i,z_{i+1}]=1\ (i\in\zz)}.
$$
The mapping $\psi$ that takes $z_i$ to $z_{i+1}$ for all $i\in\zz$,
obviously induces an automorphism of the group $L_\infty$. Consider
HNN-extension of the group $L_\infty$ with stable letter $t$ via
automorphism $\psi$ (this will be also a semidirect product of $L_\infty$
and $\zz$). We obtain the group $\pres{L_\infty,t}{t\iv z_it=z_{i+1}
\ (i\in\zz)}$ that can be simplified to one-relator group
$\pres{t,z_0}{[z_0,z_0^t]=1}$. Using Tietze transformations, one can
transform it into $L$\ \,($x=z_0t$, $y=t\iv$). So the group $L_\infty$
is a subgroup of $L$.

Note that we could consider the group
$$
L_0=\pres{z_i\ (i=0,1,2,\ldots)}{[z_i,z_{i+1}]=1\ (i=0,1,2,\ldots)}
$$
instead of $L_\infty$. The mapping $\psi$, where $\psi(z_i)=z_{i+1}$
($i=0,1,2,\ldots$), induces a monomorphism of the group $L_0$ into
itself. Indeed, the mapping $\theta$ such that $\theta(z_0)=1$,
$\theta(z_i)=z_{i-1}$ ($i=1,2,\ldots$) induces an endomorphism
of the group $L_0$ and $\theta(\psi(z))=z$ for any $z\in L_0$. If we
take an HNN-extension of the group $L_0$ with stable letter $t$ via
monomorphism $\psi$, then we get the group
$\pres{L_0,t}{t\iv z_it=z_{i+1}\ (i=0,1,2,\ldots)}$, that can be
transformed into $L$ after simplifications.
Remark that $L_0$ is obviously non-abelian (it maps onto a free group of rank $2$ by the homomorphism which maps $z_1$ to $1$, and $z_i$ to $1$ for
$i=3,4,\ldots\,$). Hence $L$ is also non-abelian, that is, $xy\ne yx$.

For $a=yx$, $b=x$ we have equality $[a,a^b]=1$ in the group $L$ and
$[a,b]\ne1$. If $L$ is a diagram group then it satisfies Condition~$2$
of Theorem~\ref{3equiv}. Thus it also satisfies Condition~$1$, that is it
contains $\zz\wr\zz$ as a subgroup. As we have mentioned above,
in the proof of Corollary~\ref{OneRelPL}, the group $\zz\wr\zz$ cannot be a
subgroup of a one-relator group because of the result of~\cite{Cheb}. The
contradiction we have obtained shows that $L$ is not a diagram group.

It remains to show that $L$ can be isomorphically embedded into a
diagram group. We apply Theorem~\ref{skew}. It follows from it that
the group $K={\cal O}(\zz,\zz)$ is a diagram group. It follows from the
description given in Section~\ref{Constr} that $K$ has a presentation
in terms of generators $g_i$, $h_i$ ($i\in\zz$), $t$ and defining
relations $[g_i,g_j]=[h_i,h_j]=1$, $g_i^t=g_{i+1}$, $h_j^t=h_{j+1}$,
where $i,j\in\zz$ and $[g_i,h_j]=1$ for $i\le j$, $i,j\in\zz$.
So it suffices to prove that the group $L=\pres{x,y}{[xy,yx]=1}$ is a
subgroup in the diagram group $K={\cal O}(\zz,\zz)$. Consider the group
$$
K_0=\pres{g_i,h_i\ (i\ge0)}
{[g_i,g_j]=[h_i,h_j]=1\ (i,j\ge0),\ [g_i,h_j]=1\ (j\ge i\ge0)}.
$$
The map $g_i\mapsto g_{i+1}$, $h_j\mapsto h_{j+1}$ ($i,j\ge0)$
can be extended to an endomorphism $\psi$ of the group $K_0$. It is a
monomorphism because the map $g_0,h_0\mapsto1$, $g_i\mapsto g_{i-1}$,
$h_j\mapsto h_{j-1}$ ($i,j\ge1)$ can be also extended to an endomorphism
$\theta$ and $\theta(\psi(z))=z$ for any $z\in K_0$. Therefore, one can
consider an HNN-extension of the group $K_0$ with a stable letter $t$ via
monomorphism $\psi$. We obtain the group
$\pres{K_0,t}{\psi(z)=z^t\ (z\in K_0)}$ that has almost the same presentation as
$K$ with the only difference that the subscripts of the presentation of $K$
run over all $\zz$. Adding new generators $g_i=g_0^{t^i}$, $h_j=h_0^{t^j}$
for negative $i$, $j$, we easily transform the presentation obtained above to
the presentation of $K$. So it suffices to prove the following Lemma.

\begin{lm}
\label{LinK}
The subgroup in $K$ generated by elements $g_ih_{i+1}$ $(i\ge0)$ and $t$
is isomorphic to $L$.
\end{lm}

\prf The group $K$ is an HNN-extension of the group $K_0$, that is,
we add the letter $t$ and relations $g_i^t=g_{i+1}$, $h_j^t=h_{j+1}$
($i,j\ge0$) to its presentation. Define the sets $\rr_k$ ($k\ge0$) of
defining relations over the alphabet $\{\,z_0,h_0,z_1,h_1,\ldots\,\}$.
For $\rr_0$ we take the set of relations of the group $L_0$, that is,
$\rr_0=\{\,[z_i,z_{i+1}]=1\ (i\ge0)\}$. Further, for $k\ge1$ we put
$$
\rr_k=\{\,z_i^{h_k}=z_i\ (0\le i<k)\,\}
\cup\{\,z_i^{h_k}=z_i^{z_{k-1}}\ (i\ge k)\,\}
\cup\{\,h_j^{h_k}=h_j\ (1\le j<k)\,\}.
$$
It is clear that $L_0=\pres{z_i\ (i\ge0)}{\rr_0}$. Let
$$
L_k=\pres{z_i\ (i\ge0),\,h_j\ (1\le j\le k)}{\rr_0\cup\rr_1\ldots\cup\rr_k}
$$
for $k\ge1$. We shall prove that for any $k\ge0$, the group $L_{k+1}$ can
be obtained from $L_k$ by a suitable HNN-extension.
\vspace{1ex}

Consider the map $\psi_k$ given by the following rules: $\psi_k(z_i)=z_i$
for $0\le i\le k$, $\psi_k(z_i)=z_i^{z_k}$ for $i>k$, $\psi_k(h_j)=h_j$
for $1\le j\le k$. Let us extend it to a homomorphism of the corresponding
free group into the group $L_k$. Let us check that all relations of the
group $L_k$ will be equalities in $L_k$ under $\psi_k$.

First of all we shall check that $\psi_k([z_i,z_{i+1}])=1$ for all $i\ge0$.
If $0\le i<k$, then $\psi_k([z_i,z_{i+1}])=[z_i,z_{i+1}]=1$ in $L_k$.
The equality $\psi_k(z_i)=z_i^{z_k}$ holds also for $i=k$. Hence
for all $i\ge k$ we also have
$\psi_k([z_i,z_{i+1}])=[z_i^{z_k},z_{i+1}^{z_k}]=[z_i,z_{i+1}]^{z_k}=1$.
Now consider the other relations of $L_k$. They have one of the following
three forms: $z_i^{h_j}=z_i$ for $0\le i<j\le k$; $z_i^{h_j}=z_i^{z_{j-1}}$
for $i\ge j$, $1\le j\le k$; $[h_i,h_j]=1$ for $1\le i<j\le k$.
Considering three cases, we map each of the relations by $\psi_k$.

If $0\le i<j\le k$ then we have
$\psi_k(z_i^{h_j})=\psi_k(z_i)^{\psi_k(h_j)}=z_i^{h_j}=z_i=\psi_k(z_i)$.

In the second case we will consider two subcases: $i\le k$ and $i>k$.
In the first subcase, that is, $1\le j\le i\le k$, we get
$\psi_k(z_i^{h_j})=\psi_k(z_i)^{\psi_k(h_j)}=z_i^{h_j}=z_i^{z_{j-1}}=
\psi_k(z_i)^{\psi_k(z_{j-1})}=\psi_k(z_i^{z_{j-1}}$. In the second subcase,
that is, $1\le j\le k<i$, we have: $\psi_k(z_i^{h_j})=
\psi_k(z_i)^{\psi_k(h_j)}=(z_i^{z_k})^{h_j}=(z_i^{h_j})^{z_k^{h_j}}=
(z_i^{z_{j-1}})^{z_k^{z_{j-1}}}=(z_i^{z_k})^{z_{j-1}}=
\psi_k(z_i)^{\psi_k(z_{j-1})}=\psi_k(z_i^{z_{j-1}})$ (equality
$z_k^{h_j}=z_k^{z_{j-1}}$ in the group $L_k$ we used in these calculations,
is a partial case of the relation of the second form for $i=k$).

In the third case everything is easy:
$\psi_k([h_i,h_j])=[\psi_k(h_i),\psi_k(h_j)]= [h_i,h_j]=1$ for
$1\le i<j\le k$.

So $\psi_k$ induces an endomorphism of the group $L_k$. Let us also
introduce the map $\theta_k$ given by the rules $\theta_k(z_i)=z_i$
for $0\le i\le k$, $\theta_k(z_i)=z_i^{z_k\iv}$ for $i>k$, $\psi_k(h_j)=h_j$
for $1\le j\le k$. One can analogously check that $\theta_k$ induces an
endomorphism of the group $L_k$. It is obvious that $\theta_k(\psi_k(z))=
\psi_k(\theta_k(z))$ for any $z\in L_k$. This means that $\psi_k$ and
$\theta_k$ are mutually inverse automorphisms of the group $L_k$.

Consider an HNN-extension of the group $L_k$ with stable letter $h_{k+1}$
via automorphism $\psi_k$ of the group $L_k$. Its presentation is
obtained from the one of $L_k$ by adding $h_{k+1}$ to the set of generators
and adding relations of the form $\psi_k(z)=z^{h_{k+1}}$ to the set of
defining relations, where $z$ runs over all generators of $L_k$. These
new relations form exactly the set $\rr_{k+1}$. Therefore, this
HNN-extension is the group $L_{k+1}$. In addition, we also have a
natural embedding of $L_k$ into $L_{k+1}$ for $k\ge0$.
\vspace{1ex}

We have a sequence of embedded subgroups
$$
L_0<L_1<\cdots<L_k<L_{k+1}<\cdots\,,
$$
that give the group
$$
\hat L=\pres{z_i\ (i\ge0),\,h_j\ (j\ge1)}
{\rr_0\cup\rr_1\cup\cdots\rr_k\cup\cdots\,}
$$
as a union of them. Let $H$ be a subgroup generated by $z_0$ and all $h_j$
($j\ge1$). Adding $h_0$ as a stable letter, we construct an HNN-extension
of the group $\hat L$ via identical endomorphism of $H$ onto itself. That is,
we add a new generator $h_0$ and relations $[z_0,h_0]=1$, $[h_j,h_0]=1$ for
all $j\ge1$. Let us describe explicitly the group $\bar L$ that we get as a
result. It has generators $z_i$, $h_i$ ($i\ge0$) subject to the following
defining relations:
$$
[z_i,z_{i+1}]=1\ \ (i\ge0),
$$
$$
[h_i,h_j]=1\ \ (i,j\ge0),
$$
$$
[z_i,h_j]=1\ \ (0\le i<j),
$$
$$
[z_i,z_{j-1}h_j\iv]=1\ \ (1\le j\le i),
$$
$$
[z_0,h_0]=1.
$$
(we took the relations $\rr_k$ for all $k\ge0$ together with the relations
added at the last step).

Let us introduce new generators $g_i=z_ih_{i+1}\iv$ ($i\ge0$). The
elements $g_i$, $h_i$ generate $\bar L$ so our aim is to describe
relations of the group $\bar L$ in terms of these generators. Replacing
elements $z_i$ by $g_ih_{i+1}$ in the defining relations of the group
$\bar L$, we get:
\be{R1}
[g_ih_{i+1},g_{i+1}h_{i+2}]=1\ \ (i\ge0),
\ee
\be{R2}
[h_i,h_j]=1\ \ (i,j\ge0),
\ee
\be{R3}
[g_ih_{i+1},h_j]=1\ \ (0\le i<j),
\ee
\be{R4}
[g_ih_{i+1},g_{j-1}]=1\ \ (1\le j\le i),
\ee
\be{R5}
[g_0h_1,h_0]=1.
\ee

Since elements of the form $h_i$ ($i\ge0$) pairwise commute, we can
simplify~(\ref{R3}), getting $[g_i,h_j]=1$ for $0\le i<j$.
Then in~(\ref{R4}) the elements $g_{j-1}$ ¨ $h_{i+1}$ commute for
$1\le j\le i$ and so~(\ref{R4}) reduces to $[g_i,g_{j-1}]=1$ for all
$1\le j\le i$. This means that all elements $g_i$ ($i\ge0$) pairwise
commute. Let us simplify~(\ref{R1}). Note that $h_{i+2}$ commutes
with the other three elements so it can be excluded. The equality
$[g_ih_{i+1},g_{i+1}]=1$ we get is equivalent to $[h_{i+1},g_{i+1}]=1$
since $g_i$ commutes with the other elements. Thus, simplifying~(\ref{R5}),
we obtain equalities $[g_i,h_i]=1$ for all $i\ge0$.

Let us summarize the above. The group $\bar L$ has generators $g_i$, $h_i$
($i\ge0$), where elements $g_i$ pairwise commute. Elements of the form $h_i$
also pairwise commute and $g_i$ commutes with $h_j$ whenever $i\le j$.
This means that the group $\bar L$ coincides with $K_0$. The group $L_0$,
naturally embedded into $\bar L$, is generated by elements $z_i$ ($i\ge0$).
So the subgroup of $\bar L$ generated by $g_ih_{i+1}=z_i$ ($i\ge0$) is
isomorphic to $L_0$.

The group $K$ is an HNN-extension of the group $K_0$ via the monomorphism
$\psi\colon K_0\to K_0$ with $t$ as a stable letter. Let us have a
subgroup $M_0$ of $K_0$ such that $\psi(M_0)\subseteq M_0$. In this case,
it is easy to show that the subgroup generated by $M_0$ and $t$ will be
the HNN-extension of $M$ via the restriction of $\psi$ on $M_0$.

Indeed, let us take such an HNN-extension. It has  the form
$M=\pres{M_0,t}{z^t=\psi(z)\ (z\in M_0)}$. The map $t\mapsto t$, $z\to z$
for $z\in M_0$ induces a homomorphism $\phi$ from $M$ to $K$. Since
$zt=t\psi(z)$ we can represented any element of the group $M$ in the form
$t^\alpha zt^{-\beta}$, where $z\in M_0$, $\alpha,\beta\in\zz$, $\beta\ge0$.
Therefore, any element $m$ in $M$ is conjugated to an element of the form
$t^\gamma z$ for some $\gamma\in\zz$, $z\in M_0$. If $m\ne1$, then either
$\gamma\ne0$ or $z\ne1$. The element $t^\gamma z$ maps to an element in
$K$ of the same form under $\psi$. It follows from the elementary
properties of HNN-extensions that it is not equal to $1$ in $K$. Therefore,
$\phi$ is an embedding of $M$ into $K$ and its image is exactly the
subgroup of $K$ generated by $M_0$ and $t$.

Note that the subgroup of $K_0$ generated by elements $z_i=g_ih_{i+1}$
($i\ge0$) is invariant under $\psi$ because $\psi(z_i)=z_{i+1}$
for all $i\ge0$. So one can regard this subgroup (isomorphic to $L_0$) as
$M_0$ and apply the arguments from the above paragraph. The corresponding
HNN-extension of it is the subgroup of $K$ generated by $t$ and $g_ih_{i+1}$
($i\ge0$). On the other hand, as we have mentioned in the beginning, this
HNN-extension is exactly $L$.
\vspace{1ex}

The Lemma and Theorem~\ref{dispr} are proved.
\vspace{2ex}

There are not many known counterexamples to the Subgroup Conjecture.
So it is natural to try to prove this conjecture under some
restrictions on the subgroup. With respect to Theorem~\ref{CommF}, we
would like to ask a few questions.

\begin{prob}
Is it true that any subgroup of R.\,Thompson's group $F$ is a diagram
group?
\end{prob}

\begin{prob}
Is it true that the commutator subgroup of any diagram group is
a diagram group?
\end{prob}

It is easy to see that the commutator subgroup of the group $F$
satisfies the following condition.
Let $\Delta$ be a diagram representing an element in $F'$, and suppose that
a conjugate diagram $\Psi\iv\Delta\Psi$ is a sum $\Gamma_1+\Gamma_2$
of two nontrivial spherical diagrams with bases $v_1$, $v_2$, respectively.
Then the diagrams $\Delta_1=\Psi(\Gamma_1+\ve(v_2))\Psi\iv$ and
$\Delta_2=\Psi(\ve(v_1)+\Gamma_2)\Psi\iv$ also belong to $F'$.
Consider any
subgroup $H$ of a diagram group $G$ that satisfies this condition.
We shall say that $H$ is {\em closed\/} in $G$.
It is easy to see that a subgroup $H$ is closed if
whenever $\Delta\in H$ and $\Delta=\Delta_1\Delta_2$, where $\Delta_1$ and $\Delta_2$ commute and do not belong to the same cyclic subgroup, we have $\Delta_1,\Delta_2\in H$.

\begin{prob}
Let $H$ be a closed subgroup in a diagram group $G$. Is it true that
$H$ is a diagram group?
\end{prob}

If the answer to the next problem is positive, then this would imply
that all word hyperbolic diagram groups are free.

\begin{prob}
Let $H$ be a subgroup in a diagram group $G$. Suppose that for any $h\in H$,
$h\ne1$, the centralizer $C_G(h)$ of $h$ in $G$ is cyclic. Does it imply
that $H$ is free $($at least for the particular case $H=G)$?
\end{prob}

At the end of this Section let us consider an interesting family
of groups. Let
$$
G_n=\pres{x_1,\ldots,x_n}
{[x_1,x_2]=[x_2,x_3]=\cdots=[x_{n-1},x_n]=[x_n,x_1]=1}.
$$
It is easy to see that $G_1=\zz$, $G_2=\zz\times\zz$,
$G_3=\zz\times\zz\times\zz$, $G_4={\cal F}_2\times{\cal F}_2$, where
${\cal F}_2$ is the free group of rank $2$. All these groups can be obtained
from $\zz$ using finite direct and free products. So these are diagram groups.
However, the group $G_5$ is not a diagram group.

\begin{thm}
\label{PartCom}
The groups $G_n$ are not diagram groups for odd $n\ge5$.
\end{thm}

\prf Let $\pp=\rav$ be a semigroup presentation and let $G=\dg$ be a
diagram group. Consider any element $g\in G$ presented by a diagram
$\Delta$. Let us decompose $\Delta$ into the sum of spherical components.
It was proved in~\cite{GuSa97} that the number of these components is an
invariant of a diagram with respect to conjugation: see the remark after
the proof of Lemma 15.15. One can see from the same Lemma that the
number of nontrivial components is also an invariant. Thus one can
introduce a function $\comp$, denoted by $\comp(g)$, the number of
nontrivial components of a diagram that represents an element $g\in G$.

Let us introduce a partial binary relation $\prec$ on $G$. Let $g_1,g_2\in G$
be such that $\comp(g_1)=\comp(g_2)=1$ (in particular, $g_1$, $g_2$ are
nontrivial). We put $g_1\prec g_2$ whenever diagrams $\Delta_1$, $\Delta_2$
that represent elements $g_1$, $g_2$ respectively, satisfy the following
condition: there are words $x,y,z\in\Sigma^*$, $u,v\in\Sigma^+$, some
$(\bs,xuyvz)$-diagaram $\Gamma$, simple absolutely reduced spherical
diagrams $\Psi_1$ and $\Psi_2$ with bases $u$, $v$ respectively such that
$\Gamma\iv\Delta_1\Gamma=\ve(x)+\Psi_1+\ve(yvz)$,
$\Gamma\iv\Delta_2\Gamma=\ve(xuy)+\Psi_2+\ve(z)$. It follows from this
definition that if $g_1\prec g_2$, then elements $g_1$, $g_2$ commute and
generate a subgroup isomorphic to $\zz\times\zz$ in $G$. In particular,
they do not belong to the same cyclic subgroup. Thus the relation $\prec$
is antireflexive, that is, condition $g\prec g$ never holds for $g\in G$.
Let us establish a few properties of $\prec$.

\begin{lm}
\label{trans}
Let $\pp=\rav$ be a semigroup presentation and let $G=\dg$ be a diagram
group. The relation $\prec$ is transitive, that is, for any $g_1,g_2,g_3\in G$
such that $\comp(g_1)=\comp(g_2)=\comp(g_3)=1$, conditions $g_1\prec g_2$
and $g_2\prec g_3$ imply $g_1\prec g_3$.
\end{lm}

\prf Let element $g_i$ in $G$ be represented by a diagram $\Delta_i$ ($i=1,2,3$).
Since $g_1\prec g_2$, there are words $x,y,z\in\Sigma^*$, $u,v\in\Sigma^+$,
$(\bs,xuyvz)$-diagram $\Gamma$, simple absolutely reduced spherical
diagrams $\Psi_1$ and $\Psi_2$ with bases $u$, $v$ respectively, such that
$\Gamma\iv\Delta_1\Gamma=\ve(x)+\Psi_1+\ve(yvz)$,
$\Gamma\iv\Delta_2\Gamma=\ve(xuy)+\Psi_2+\ve(z)$. Since $g_2\prec g_3$,
there are words $x',y',z'\in\Sigma^*$, $u',v'\in\Sigma^+$,
$(\bs,x'u'y'v'z')$-diagram $\Gamma'$, simple absolutely reduced spherical
diagrams $\Psi_2'$ and $\Psi_3'$ with bases $u'$, $v'$ respectively,
such that $(\Gamma')\iv\Delta_2\Gamma'=\ve(x')+\Psi_2'+\ve(y'v'z')$,
$(\Gamma')\iv\Delta_3\Gamma'=\ve(x'u'y')+\Psi_3'+\ve(z')$. It follows
from these conditions that
$\ve(xuy)+\Psi_2+\ve(z)=\Theta\iv(\ve(x')+\Psi_2'+\ve(y'v'z'))\Theta$, where
$\Theta=(\Gamma')\iv\Gamma$. Diagrams $\ve(xuy)+\Psi_2+\ve(z)$ and
$\ve(x')+\Psi_2'+\ve(y'v'z')$ are conjugate by an element $\Theta$.
It follows from~\cite[Lemma 15.15]{GuSa97} that the components of these
diagrams are conjugate respectively. In particular, words $x'$ and $z$ are
nonempty. Applying this Lemma, we conclude that
$\Theta=\Theta_1+\Theta_2+\Theta_3$, where $\Theta_1$, $\Theta_2$, $\Theta_3$
are $(x',xuy)$-, $(u',v)$- and $(y'v'z',z)$-diagrams, respectively.

Let us now take the diagram $\Xi=\Gamma(\ve(xuyv)+\Theta_3\iv)$. It is clear
that $\Xi=\Gamma'\Theta(\ve(xuyv)+\Theta_3\iv)=
\Gamma'(\Theta_1+\Theta_2+\Theta_3)(\ve(xuyv)+\Theta_3\iv)=
\Gamma'(\Theta_1+\Theta_2+\ve(y'v'z'))$. We have $\Xi\iv\Delta_1\Xi=
(\ve(xuyv)+\Theta_3\iv)\iv(\Gamma\iv\Delta_1\Gamma)(\ve(xuyv)+\Theta_3\iv)=
(\ve(xuyv)+\Theta_3\iv)\iv(\ve(x)+\Psi_1+\ve(yvz))(\ve(xuyv)+\Theta_3\iv)=
\ve(x)+\Psi_1+\ve(yvy'v'z')$, and $\Xi\iv\Delta_3\Xi=
(\Theta_1+\Theta_2+\ve(y'v'z'))\iv(\Gamma')\iv\Delta_3\Gamma'
(\Theta_1+\Theta_2+\ve(y'v'z'))=(\Theta_1\iv+\Theta_2\iv+\ve(y'v'z'))
(\ve(x'u'y')+\Psi_3'+\ve(z'))(\Theta_1+\Theta_2+\ve(y'v'z'))=
\ve(xuyvy')+\Psi_3'+\ve(z')$. Thus conjugating diagrams $\Delta_1$
and $\Delta_3$ by a $(\bs,xuyvy'v'z')$-diagram $\Xi$, we represent them in
the form enabling us to conclude that $g_1\prec g_3$.

The proof is complete.
\vspace{1ex}

Lemma~\ref{trans} implies that the relation $\prec$ is also antisymmetric,
that is $g_1\prec g_2$ excludes $g_2\prec g_1$. Let us establish one more
property of $\prec$.

\begin{lm}
\label{LinOrd}
Let $\pp=\rav$ be a semigroup presentation and let $G=\dg$ be a diagram
group. We claim that for any commuting elements $g_1,g_2\in G$ that do not
belong to the same cyclic subgroup and satisfy $\comp(g_1)=\comp(g_2)=1$,
exactly one of the following conditions holds: $g_1\prec g_2$ or
$g_2\prec g_1$.
\end{lm}

\prf Let $A_i$ be a diagram that represents an element $g_i\in G$
($i=1,2$). Since $[g_1,g_2]=1$, we can apply Theorem~\ref{prws-cmmt}
and find a word $v=v_1\ldots v_n$, spherical $(v_j,v_j)$-diagrams
$\Delta_j$ ($1\le j\le n$), integers $d_{ij}$ ($1\le i\le 2$, $1\le j\le n$)
and some $(\bs,v)$-diagram $\Gamma$ such that
$\Gamma\iv A_i\Gamma=\Delta_1^{d_{i1}}+\cdots+\Delta_n^{d_{in}}$,
where diagrams $\Delta_j$ ($1\le j\le n$) are either trivial or simple
absolutely reduced. The condition $\comp(g_1)=1$ means that there is
exactly one number $j$ from $1$ to $n$ such that diagram $\Delta_j^{d_{1j}}$
is not trivial. Analogously, condition $\comp(g_2)=1$ means that there
exists exactly one number $k$ from $1$ to $n$ such that diagram
$\Delta_k^{d_{1k}}$ is not trivial. If $j=k$ then diagrams $A_1$, $A_2$
belong to the same cyclic subgroup of $G$ but this is impossible.
If $j<k$, then $g_1\prec g_2$ by definition. If $k<j$ then $g_2<g_1$.

The proof is complete.
\vspace{1ex}

Let us continue the proof of Theorem~\ref{PartCom}. Let $n=2k+1$, $k\ge2$.
Suppose that $G_n=\dg$ is a diagram group over $\pp=\rav$ with base $\bs$.
First of all let us prove that $\comp(x_i)=1$ for all generators $x_i$
of $G_n$. Let us establish that the centralizer $C(x_i)$ of the element $x_i$
($1\le i\le n$) in $G_n$ is the subgroup generated by elements $x_{i-1}$,
$x_i$ and $x_{i+1}$, isomorphic to the direct product ${\cal F}_2\times\zz$
(subscripts are taken modulo $n$). By symmetry, it suffices to consider
the centralizer of $x_n$.

Suppose that $C(x_n)\ne\grp{x_1,x_{n-1},x_n}$. Consider a group word
$W$ of minimal length in $x_1$, \dots, $x_n$ such that $W\in C(x_n)$,
$W\notin\grp{x_1,x_{n-1},x_n}$. In particular, the word $W$ is nonempty
and it has neither nonempty initial nor nonempty terminal segment that
belongs to $\grp{x_1,x_{n-1},x_n}$.
The group $G_n$ is an HNN-extension with base
$$
\Gamma=\pres{x_1,\ldots,x_{n-1}}
{[x_1,x_2]=[x_2,x_3]=\cdots=[x_{n-2},x_{n-1}]=1},
$$
and stable letter $x_n$, with respect to the identical automorphism of
the subgroup $\grp{x_1,x_{n-1}}$ of $\Gamma$. Consider the element
$x_n\iv Wx_nW\iv$ of this HNN-extension. It equals $1$ in the group $G_n$
since $W\in C(x_n)$. By Britton's Lemma (see~\cite{LS80rus}), the word
$x_n\iv Wx_nW\iv$ has a subword of the form $U=x_n^{-\delta}Vx_n^\delta$,
where $\delta=\pm1$, $V\in\grp{x_1,x_{n-1}}$ is a word that does not
contain $x_n^{\pm1}$. Since the word $W$ is chosen to have minimal length,
$U$ is not contained in $W^{\pm1}$. Otherwise the occurrence of $U$ can
be replaced by an occurrence of the word $V$ that is equal to $U$ in $G_n$,
decreasing the length of $W$. It is clear that $V$ is nonempty since
$W$ cannot begin or end with $x_n^{\pm1}$. Thus $V$ is neither initial
nor terminal segment of $W^{\pm1}$. So it is clear that $U$ does not
occur in $x_n\iv Wx_nW\iv$. We got a contradiction.

Thus $C(x_n)=\grp{x_1,x_{n-1},x_n}$. Consider a mapping of the alphabet
$\{\,x_1,\ldots,x_n\,\}$ into $G_n$, sending each of the elements $x_1$,
$x_{n-1}$, $x_n$ to itself and senging all the other elements to $1$.
Extending this mapping to a homomorphism of the corresponding free group
into $G_n$, we see that all defining relations of the group $G_n$ are sent
to $1$. Thus we have an induced homomorphism $\phi\colon G_n\to G_n$. It is
obvious that it is a retraction, that is, $\phi^2=\phi$. On the one hand,
the subgroup $\phi(G_n)$ of $G_n$ equals $\grp{x_1,x_{n-1},x_n}$; on the
other hand, this group is presented by relations of the group $G_n$ with
additional conditions $x_2=\cdots=x_{n-2}=1$. Thus for any $n>3$ we have
$$
\grp{x_1,x_{n-1},x_n}=\phi(G_n)=
\pres{x_1,x_{n-1},x_n}{[x_{n-1},x_n]=[x_n,x_1]=1}\cong{\cal F}_2\times\zz,
$$
as desired.
\vspace{0.5ex}

So $C(x_i)\cong{\cal F}_2\times\zz$ for all $i$ from $1$ to $n$;
in particular, the centre of $C(x_i)$ is cyclic. If $x_i$ were represented
by a diagram with more than one nontrivial component, then its centralizer
would have at least two direct summands isomorphic to $\zz$
by~\cite[Theorem 15.35]{GuSa97}. So its centre would not be cyclic.
Taking into account that $x_i$ is nontrivial, we conclude that $\comp(x_i)=1$.
(Notice that we have not used yet that $n$ is odd.)

Apply Lemma~\ref{LinOrd} and suppose without loss of generality that
$x_1\prec x_2$. Suppose that $x_2\prec x_3$. Then Lemma~\ref{trans} would
imply that $x_1\prec x_3$, so elements $x_1$ and $x_3$ commute. It is clear
that these elements do not commute in $G_n$. The contradiction we obtain
allows to apply Lemma~\ref{LinOrd} again and to conclude that $x_3\prec x_2$.
We will obtain a contradiction again if we suppose that $x_4\prec x_3$.
So in fact $x_3\prec x_4$. Continuing in this way, we shall conclude that
$x_{2k+1}\prec x_{2k}$, $x_{2k+1}\prec x_1$, $x_2\prec x_1$. We have
a contradiction.
\vspace{1ex}

The Theorem is proved.
\vspace{2ex}

It is reasonable to pose a question with respect to Theorem~\ref{PartCom}.

\begin{prob}
\label{probGn}
For which $n$ the groups $G_n$ are diagram groups? For which $n$ they are
isomorphically representable by diagrams?
\end{prob}

If there is an odd $n\ge5$ such that the group $G_n$ is representable by
diagrams, then we have one more counterexample to the Subgroup Conjecture.
Otherwise we would have a generalization of Theorem~\ref{PartCom}.

\section{Distortion of Subgroups in Diagram Groups}
\label{Disto}

The problems that concern distortion in groups form a branch of geometric
group theory under development (see~\cite{Grom93,Olsh97,Olsh98}). Let us
recall some definitions.

Let $A$ be a group with finite set of generators $X$. In this case, for
any $g\in G$ there exists an $n\ge0$ and $x_1,\ldots,x_n\in X^{\pm1}$
such that $g=x_1\ldots x_n$. The least $n$ with this property is
called the {\em length\/} of the element $g$ with respect to the
generating set $X$ and it is denoted by $|g|_X$.

If there are two functions $\phi$, $\psi$ from $G$ to the set of all
nonnegative integers, then we shall write $\phi\preceq\psi$, whenever
there is a positive integer constant $C$ such that $\phi(g)\le C\psi(g)$
for all $g\in G$. If it holds $\phi\preceq\psi$ and $\psi\preceq\phi$
for the two functions simultaneously, then we call these functions
{\em equivalent\/} and denote this fact by $\phi\sim\psi$. Obviously,
$\sim$ is in fact the equivalence relation (one does not have to mix
it with another equivalence relation that is often used when the Dehn
functions are discussed). So, for the two functions one has $\phi\sim\psi$
if and only if there exisets a positive integer constant $C$ such that
$$
\frac{\phi(g)}C\le\psi(g)\le C\phi(g)\quad\mbox{ for all }g\in G.
$$

If $X$ and $Y$ are finite sets of generators of the same group $A$,
then elementary arguments show that functions $|\ |_X$ and $|\ |_Y$
are equivalent.

Let we have two finitely generated groups $A$ and $B$ such that $A$ is
a subgroup of $B$. Let us fix some finite system of generators $X$ for
the group $A$ and some finite system of generators $Y$ for the group $B$.
For any element $g\in A$ we define two numbers:  $|g|_X$ and $|g|_Y$.
Functions $|\ |_X$, $|\ |_Y$ can be regarded as functions on $A$.
Later we will compare functions on two groups one embedded into
another, with respect to $\preceq$, taking the corresponding restrictions
of these functions. It follows from elementary reasons that
$|\ |_Y\preceq|\ |_X$. If the converse is true, that is,
$|\ |_X\preceq|\ |_Y$ holds, the we say that a subgroup $A$ embeds into
$B$ {\em quasiisometrically\/} or {\em without distortion\/} (this
happens, if $|\ |_X\sim|\ |_Y$). Note that the equivalence of two
length functions $|\ |_X$ and $|\ |_Y$ does not depend on the choice
of finite systems of generators $X$ and $Y$.  If to consider functions up
to equivalence, then one can introduce length functions $\ell_A$ and $\ell_B$
in finitely generated groups $A$ and $B$, respectively, that will depend
of $A$ and $B$ only. The quasiisometricity of an embedding of $A$ into $B$
means that $\ell_A\sim\ell_B$.

Now consider a more general situation of an embedding of $A$ into $B$ for
two finitely generated groups, $A\le B$. Let we distinguish some
finite generating sets $X$ and $Y$ in groups $A$ and $B$
respectively.  One can consider the function
$$
\disto(n)=\max\limits_{|g|_Y\le n}|g|_X,
$$
that describes distortion of the subgroup $A$ embedded into $B$. It is
called the {\em distortion function\/} of the subgroup $A$ in $B$. It is
easy to find out that if we change the generating sets, then the
distortion function $\disto(n)$ is not essentially changed. The reader
can easily write down the corresponding inequalities. Thus we can talk
about linear, quadratic, polynomial, exponential etc distortion. The
question about distortion is aslo interesting with respect to the so
called membership problem. Let we have two finitely generated groups
$A$ and $B$, where $A\le B$. The {\em membership problem\/} of elements
of the group $B$ into the subgroup $A$ is the question on the existence
of an algorithm that decides, given a word on the generators of $B$,
whether the element of $B$ presented by this word belongs to $A$.
The membership problem of elements of $B$ into a subgroup $A$ is
decidable if and only if the distortion function $\disto(n)$ defined
above is recursive (equivalently, has a recursive upper bound).

It is interesting to find the conditions under which all finitely
generated subgroups of a given group will embed into it (or any its
finitely generated subgroup) without distortion. Free groups and
abelian groups have this property. Even for the case of nilpotent
groups the situation is quite different: in any nilpotent (non-abelian)
torsion-free group there are cyclic subgroups that have distortion in
them. (Note that diagram groups may have distorted subgroups in general:
due to the classical result of Mikhailova~\cite{Mikh58}, the group
${\cal F}_2\times{\cal F}_2$ has finitely generated subgroups with
undecidable membership problem, so they are distorted.)

Let us mention two recent results of Burillo~\cite{Bur}: he proved that
every cyclic subgroup of R.\,Thompson's group $F$ is embedded into it
without distortion. Also he gave examples of quasi-isometric embeddings
of groups $F\times\zz^n$ ($n\ge1$) and $F\times F$ into $F$.

It is thus natural to ask whether every finitely generated subgroup
embeds quasi-isometrically into $F$. We give a negative answer. Namely,
for any integer $d\ge2$ we construct a finitely generated subgroup of
$F$ with distortion at least $n^d$. The fact that any cyclic subgroup
of {\bf every} finitely generated group representable by diagrams
(including the case of $F$) embeds quasi-isometrically into it, follows
easily from~\cite[Lemma 15.29]{GuSa97}.  We shall prove a more general
result.

\begin{thm}
\label{Znqi}
Let $B$ be a finitely generated subgroup of a diagram group $G$ and let
$A$ be a finitely generated abelian subgroup of $B$. Then $A$
embeds into $B$ quasi-isometrically.
\end{thm}

\prf Note that diagram groups are torsion-free~\cite[Theorem 15.11]{GuSa97}
and so $A$ is isomorphic to $\zz^m$ for some integer $m$. (This also
follows from Theorem~\ref{nilp}.) Let $G=\dg$, where $\pp=\rav$ is a
semigroup presentation. By $A_1$, \dots, $A_m$ we denote spherical
diagrams with base $\bs$ presenting free generators of $A\cong\zz^m$.
Since these elements are pairwise commutative in $G$, we can apply
Theorem~\ref{prws-cmmt} to them. Thus we have a word $v=v_1\ldots v_n$,
sperical $(v_j,v_j)$-diagrams $\Delta_j$ $(1\le j\le n)$, integers
$d_{ij}$ $(1\le i\le m$, $1\le j\le n)$ and some $(\bs,v)$-diagram $\Gamma$
such that
$$
\Gamma\iv A_i\Gamma=\Delta_1^{d_{i1}}+\cdots+\Delta_n^{d_{in}}
$$
for all $1\le i\le m$. Each of the diagrams $\Delta_1$, \dots, $\Delta_n$
is either trivial or simple absolutely reduced. Conjugation by diagram
$\Gamma$ is an isomorphism of groups $G=\dg$ and ${\cal D}(\pp,v)$.
Since the property of a subgroup to be embeddable quasi-isometrically
is an invariant under isomorphism, we can assume without loss of generality
that $B$ is a subgroup of ${\cal D}(\pp,v)$.

It suffices to prove that $\ell_A\preceq\ell_B$. Consider the diagrams
$$
\Delta_j'=\ve(v_1\ldots v_{j-1})+\Delta_j+\ve(v_{j+1}\ldots v_n)
$$
for all $j$ from $1$ to $n$. By $\Psi_1$, \dots, $\Psi_r$ we denote
those of diagrams $\Delta_1'$, \dots, $\Delta_n'$ that are nontrivial.
The form a basis $X$ of a free abelian group $C$. Since $A$ embeds into
$C$ quasi-isometrically, we have $\ell_A\sim\ell_C$. For any element in
the group ${\cal D}(\pp,v)$ presented by a reduced diagram $\Delta$,
we denote by $\#(\Delta)$ the number of cells in $\Delta$. Thus $\#$ is
a function on the diagram group. Let $Y$ be a finite generating set of
the group $B$ and let $K$ be the greatest number of cells for diagrams
in $Y$. Then it is obvious that $\#(\Delta)\le K|\Delta|_Y$ for any
diagram $\Delta$ in $B$. So we have $\#\preceq\ell_B$. Let $s_1$, \dots,
$s_r$ be arbitrary integers. Consider the element
$\Delta=\Psi_1^{s_1}\ldots\Psi_r^{s_r}$ in $C$. All diagrams $\Delta_j$
($1\le j\le n$) are absolutely reduced. So it follows easily from the
definition of diagrams $\Psi_k$ ($1\le k\le r$) that
$\#(\Delta)=|s_1|\#(\Psi_1)+\cdots+|s_r|\#(\Psi_r)$. Since
$|\Delta|_X=|s_1|+\cdots+|s_r|$, we can deduce inequality
$|\Delta|_X\le\#(\Delta)\le K'|\Delta|_X$, where $K'$ is the greatest
number of cells for the diagrams in $X$. Therefore, $\ell_C\sim\#$.

Summarizing what we have said above, we conclude that
$\ell_A\sim\ell_C\sim\#\preceq\ell_B$. Now obvious inequality
$\ell_B\preceq\ell_A$ gives us the equivalence $\ell_A\sim\ell_B$.

The Theorem is proved.
\vspace{2ex}

Now consider R.\,Thompson's group $F$, take any its element $g\in F$
and its centralizer $C_F(g)$ in $F$. In~\cite[Corollary 15.36]{GuSa97}
we gave the description of centralizers in $F$: they are finite direct
products of groups that isomorphic to either $F$ or $\zz$. In particular,
all of them are finitely generated. Remark that if an element $g\in F$
is presented by a diagram $\Delta$, then to find its centralizer, one
needs to find an absolutely reduced diagram $\Delta'$ conjugated to $\Delta$
(this can be done effectively by~\cite[Lemma 15.14]{GuSa97}) and then
decompose $\Delta'$ into a sum of (spherical) components. To each trivail
component we assign the group $F$ and to each nontrivial one we assign
$\zz$. Then we take direct product of these groups. It is easy to see
that the groups we get in this way are exactly groups of the form
$F^m\times\zz^n$, where $0\le m\le n+1$.

The Theorem below generalizes Burillo's results from~\cite{Bur}, where
it is shown that $F$ has quasi-isometrically embedded subgroups
isomorphic to $F\times F$ (Proposition~9), and for every $n\ge1$ there
are quasi-isometrically embedded subgroups isomorphic to $F\times\zz^n$
(Corollary~6). (Although the group $F\times F$ cannot be a centralizer
of an element in $F$, it is embeddable without distortion into
$F^2\times\zz$, which is a centralizer of some element in $F$. This
implies the first of results quoted above.)

\begin{thm}
\label{q-i-C}
For any element $g$ in R.\,Thompson's group $F$, the centralizer
$C_F(g)$ of this element embeds into $F$ quasi-isometrically.
\end{thm}

First of all we need a lemma that can be deduced easily as a consequence
of~\cite[Proposition 2]{Bur}. But we give a direct proof of the fact
we need.

Let $\pp=\pres{x}{x^2=x}$. For any $k\ge1$, by $\#_k(g)$ we denote
the number of cells in in the (reduced) spherical diagram with base
$x^k$ that presents the element $g\in F\cong{\cal D}(\pp,x^k)$. The
number $|g|$ denotes the length of $g\in F$ with respect to the set
$\{x_0,x_1\}$ of generators.

\begin{lm}
\label{estim}
For any $g\in F$, the following inequalities hold:
$$
\frac{|g|}3\le\#_3(g)\le2|g|.
$$
For any $k$ the function $\#_k$ is equivalent to the length
function $|\ |$.
\end{lm}

\prf Diagrams that correspond to the paths $(x^2,x\to x^2,1)(1,x^2\to x,x^2)$
and $(x,x\to x^2,x)(1,x^2\to x,x^2)$ have two cells each. They present
the elements $x_0$, $x_1$ of R.\,Thompson's group $F\cong{\cal D}(\pp,x^3)$.
If $g$ is an element of length $n$, then it can be presented by a diagram
with base $x^3$ that has at most $2n$ cells. The second inequality is
thus proved.

Let us prove the first inequality. Let an element $g$ is presented by
a diagram $\Delta$ with base $x^3$. The longest positive path from
$\iota(\Delta)$ to $\tau(\Delta)$ cuts this diagram into two parts:
positive and negative one. The number of cells in each of the parts
is
the same, let it be equal to $m$. Then $\#_3(g)=2m$. It is easy to see
that the longest positive path has length $m+3$. Represent $g$ as a
normal form $g=g_1g_2^{-1}$, where each of the elements $g_1$, $g_2$
is positive, that is, it is a product of positive exponents of
generators. Since $|g|\le|g_1|+|g_2|$, it suffices to estimate the
length of $g_1$. (The lenght of $g_2$ can be estimated analogously.)
So let $g_1=x_0^jx_{i_1}\ldots x_{i_s}$, where $s\ge0$ and
$1\le i_1\le\ldots\le i_s$ is the normal form of $g_1$. For any $i\ge1$,
replace $x_i$ by $x_0^{1-i}x_1x_0^{i-1}$, which is equal to it in $F$.
Then we have that $g_1$ equals in $F$ to the word
$$
x_0^{j-i_1+1}x_1x_0^{i_1-i_2}x_1\ldots x_0^{i_{s-1}-i_s}x_1x_0^{i_s-1},
$$
that has length
$$
|j-i_1+1|+(i_2-i_1)+\cdots+(i_s-i_{s-1})+i_s-1+s=2i_s-i_1+|j-i_1+1|+s-1.
$$
(If $s=0$, then the length is just $j$.)

according to the procedure described in~\ref{Prelim} (Example~ref{TGNF}),
we have inequality $s+j\le m$. If $s\ge1$, then the element $x_{i_s}$
corresponds to the edge $(x^t,x\to x^2,x^{i_s})$ in the Squier complex,
where $t\ge1$, so $\Delta$ ¥áâì has a positive path labelled by
$x^{t+2+i_s}$. This implies $t+2+i_s\le m+3$ hence $i_s\le m$.

Let us consider two cases.

 ) $j\ge i_1-1$ or $s=0$.\ \ We have $|g_1|\le2i_s+j+s-2i_1\le3m-2$ for
$s\ge1$. If $s=0$, then $|g_1|=j\le m$.

¡) $s\ne0$, $j<i_1-1$.\ \ In this case $|g_1|\le2i_s+s-j-2\le3m-2$.

Summarizing, we conclude that $|g_1|\le3m$ for all cases. Also
$|g_2|\le3m$. Therefore, $|g|\le6m=3\#_3(g)$, what we had to prove.

Now it remains to note that $|\#_k(g)-\#_3(g)|\le2|k-3|$ since the
diagram that presents $g$ in ${\cal D}(\pp,x^k)$ can be obtained
from $\Delta$ conjugating it by a diagram of $|k-3|$ cells. From
what follows that functions $\#_k$ and $\#_3$ are equivalent (one needs
to use that $\#_k(g)=0$ if and only if $g=1$.)

The Lemma is proved.
\vspace{2ex}

{\bf Proof of Theorem~\ref{q-i-C}.}\ Let $g\in F$ be an arbitrary
element. Let us present it by an $(x,x)$-diagram and reduce this
diagram to absolutely reduced form by conjugation. We get some
diagram $\Delta$ with base $x^k$. By Lemma~\ref{estim}, the length
function in $F$ is equivalent to $\#_k$. Let $\Delta=\Delta_1+\cdots+\Delta_m$
is a decomposition of $\Delta$ into the sum of components, where
$\Delta_i$ is a spherical diagram with base $z_i$ ($1\le i\le m$).
Any element in the centralizer of $\Delta$ is equal to a sum of
$(z_i,z_i)$-diagrams, and the $i$th summand commutes with $\Delta_i$
($1\le i\le m$). For any $i$ from $1$ to $m$, let $G_i=F$, if $\Delta_i$
is trivial and $G_i=\zz$, if $\Delta_i$ is nontrivial. Thus
$C_F(g)\cong G_1\times\cdots\times G_m$. Choose a system of generators
in each of the groups $G_i$: if $G_i=F$, then the system consists of
$x_0$, $x_1$, and for $G_i=\zz$ the system consists of one element.
these systems of generators form a generating set $Y$ of the centralizer
of $\Delta$. It is clear that any element $h$ in the centralizer can be
uniquely presented in the form $h_1\ldots h_m$, where $h_i\in G_i$ for
$1\le i\le m$, and $|h|_Y=|h_1|_1+\cdots+|h_m|_m$ (by $|\ |_i$ we
denote the length in $G_i$ with respect to the generating set we
have chosen). The number of cells in $\Delta$ is equal to the sum
of numbers of cells in diagrams $\Delta_i$ ($1\le i\le m$). So it follows
from the equivalence of the length function and the number of cells
that the function $|\ |_Y$ is equivalent to the length function in
$F\cong{\cal D}(\pp,x^k)$. This means that the embedding of $C_F(g)$ into
$F$ is quasi-isometric. (ɇǬ $G_i=F$, then we apply Lemma~\ref{estim}.
In the case $G_i=\zz$ the equivalence of the length function
and the number of cells is obvious.)

The Theorem is proved.
\vspace{2ex}

Before going to the proof of the next result about distorted subgroups
in $F$, let us consider the following construction that has its
preimage in~\cite{Mikh58}. Let $H$ be a group generated by a finite set
$X$ and let $R$ be a finite subset in $H$. By $N$ we denote the normal
closure of the set $R$ in $H$. Consider the subgroup $K$ in $H\times H$
generated by all elements of the form $(x,x)$, where $x\in X$, and also
all elements of the form $(r,1)$, where $r\in R$. it is easy to see that
for any $g,h\in H$, the element $(g,h)$ is in $K$ if and only if the
cosets of $g$ and $h$ by the subgroup $N$ are equal. (This is proved
in the same way as in~\cite{Mikh58}; see also~\cite{LS80rus}.)

It is possible to consider an analog of the Dehn function in this
situation. For any element $g\in N$ by $k(g)$ we denote the least $k$
such that the element $g$ is equal in $H$ to a product of $k$ elements
conjugated in $H$ to elements in $R$ or their inverses.  Let
$$
\Phi(n)=\max\limits_{|g|\le n}k(g),
$$
where $|g|$ is the length of $g$ with respect to the set $X$ of
generators. This function can be call a (relative) {\em Dehn function\/}
of presentation $\pres{X}{R}$ {\em with respect\/} to $H$; it is clear
that if $H$ is free, then we have standard Dehn function.

Let $Y$ denote the above set of generators of $K$. Suppose that the element
$(g,1)\in K$ can be presented as a product of $m$ elements from $Y^{\pm1}$.
Then we have an equality
$$
(g,1)=(u_0,u_0)(r_1,1)^{\ve_1}(u_1,u_1)\ldots(r_m,1)^{\ve_m}(u_m,u_m),
$$
that holds in $K$, where $u_0,u_1,\ldots,u_m\in H$, $r_1,\ldots,r_m\in R$,
$\ve_1,\ldots,\ve_m=\pm1$. Therefore, equalities
$g=u_0r_1^{\ve_1}u_1\ldots r_m^{\ve_m}u_m$, $1=u_0u_1\ldots u_m$ hold
in $H$. Then
$$
g=u_0u_1\ldots u_mr_1^{\ve_1u_1\ldots u_m}\ldots r_m^{\ve_mu_m}=
r_1^{\ve_1u_1\ldots u_m}\ldots r_m^{\ve_mu_m}
$$.
So the inequality $k(g)\le m$ holds. In particular, representing $(g,1)$
as a product of the least number of generators in $Y^{\pm1}$, we get
the inequality $k(g)\le|(g,1)|_K$. For each positive integer $n$ we
have an element $g\in H$ such that $|g|_X\le n$ and $\Phi(n)=k(g)$.
The group $H\times H$ has the following natural set of generators:
$Z=(X\times\{1\})\cup(\{1\}\times X)$. It is clear that
$|(g,1)|_Z\le|g|_X\le n$, but we have $|(g,1)|_K\ge k(g)=\Phi(n)$.
It follows from the definition of the distortion function that
$\disto(n)\ge\Phi(n)$, where we embed $K$ into $H\times H$. We proved
the following lemma.

\begin{lm}
\label{Phi-Dist}
Let $H$ be a group generated by a finite set $X$, let $R$ be a finite
subset of $H$. Consider the subgroup $K$ of $H\times H$ generated by
the set $Y$ that consists of all elements $(x,x)$ $(x\in X)$ and all
elements of the form $(r,1)$ $(r\in R)$. Then inequality
$\disto(n)\ge\Phi(n)$ holds, where $\disto(n)$ is the distortion function
for the embedding of the group $K$ generated by $Y$ into the group
$H\times H$ generated by $Z=(X\times\{1\})\cup(\{1\}\times X)$. Here
$\Phi(n)$ is the relative Dehn function of presentation $\pres{X}{R}$
with respect to $H$.
\end{lm}

An important property of R.\,Thompson's group $F$ is that $F\times F$
can be embedded into $F$. So in order to obtain distorted subgroups in $F$
we need to take such a subgroup $H$ with finite generating set $X$ and
a finite subset $R$ of $H$ such that the Dehn function of $\pres{X}{R}$
with respect to $H$ will be overlinear. Then, in the above notation,
we shall get an embedding of $K$ into the group $H\times H$, which is in
turn embeddable into $F\times F$ (and so it embeds into $F$). The
embedding of $K$ into $F$ will not be quasi-isometric. Note that if to
take $H=F$, $R=\{\,[x_0,x_1]\,\}$, then the {\bf relative\/} Dehn
function will be linear though the standard Dehn function (with respect
to the free group on $\{x_0,x_1\}$) will be quadratic. Now we will give
an example how to construct a subgroup in $F$ with at least quadratic
distortion. The last Theorem in this Section will be a generalization
of this example.

\begin{ex}
\label{Disto-n2}
{\rm Let $H=\zz\wr\zz$ be a subgroup of $F$ constructed in Section~\ref{Sol}.
Denote its generators by $a$ and $b$. Let the elements $a_n=a^{b^n}$
($n\in\zz$) form a basis of the free abelian subgroup. For $R$ we take
the set of a single element $[a,b]=a_0\iv a_1$. Conjugating this
element by all elements in $H$, we shall get all elements of the form
$c_i=a_i\iv a_{i+1}$ ($i\in\zz$). It is obvious that all elements of the
form $c_i$ are also a basis of the free abelian group. Let
$g_n=[a^n,b^n]=a_0^{-n}a_n^n\in H$. The length of $g_n$ with respect
to $\{a,b\}$ does not exceed $4n$. At the same time, we have equality
$g_n=c_0^nc_1^n\ldots c_{n-1}^n$, which shows that $g_n$ can be
presented as a product of $n^2$ elements of the form $c_i$ ($i\in\zz$).
Since the elements $c_i$ form a basis of a free abelian subgroup, it
follows that $g_n$ cannot be presented as a product of less than $n^2$
elements of the form $c_i^{\pm1}$. Therefore, the Dehn function $\Phi(n)$
of presentation $\pres{a,b}{[a,b]}$ with respect to $H$ satisfies
inequality $\Phi(4n)\ge n^2$. Let $K$ be a subgroup of $H\times H$
generated by $(a,a)$, $(b,b)$, $([a,b],1)$. Lemma~\ref{Phi-Dist}
shows that the distortion function $\disto$ that characterizes the
embedding of $K$ into $H\times H$, is at least quadratic. In particular,
$K$ embeds into $H\times H$ with distortion (that is, the embedding is
not quasi-isometric). It remains to embed $H\times H$ into $F\times F$
and then into $F$. Taking into account that $\ell_F\preceq\ell_{H\times H}$,
we obtain that $K$ embeds into $F$ with distortion.

One can give explicit expressions (in terms of normal forms) of the
generators of $K$ as a subgroup in $F$. The elements $a=x_1x_2x_1^{-2}$
and $b=x_0$ generate in $F$ a subgroup isomorphic to $\zz\wr\zz$. The rules
$(x_0,1)\mapsto x_1x_2x_1^{-2}$, $(x_1,1)\mapsto x_1^2x_2x_1^{-3}$,
$(1,x_0)\mapsto x_2x_3x_2^{-2}$, $(1,x_1)\mapsto x_2^2x_3x_2^{-3}$ give
an embedding of $F\times F$ into $F$. Using that, it is easy to compute
the generators of $K$. The following elements of $F$ generate the
subgroup isomorphic to $K$:
$$
x_1^2x_2^2x_6^2x_7^2x_8\iv x_7\iv
x_6^{-2}x_3\iv x_2\iv x_1^{-2},
$$
$$
x_1x_2x_4x_5x_4^{-2}x_1^{-2},
$$
$$
x_1^3x_2^2x_5x_6x_5^{-2}x_3\iv x_2\iv x_1^{-3}.
$$
}
\end{ex}

Now let us prove the result about distorted subgroups of $F$ in its
general form.

\begin{thm}
\label{Disto-nd}
For any $d\ge2$, there exists a finitely generated subgroup $K_d$ of
R.\,Thomp\-son's group $F$ such that the corresponding distortion
function satisfies inequality $n^d\preceq\disto(n)$.
\end{thm}

\prf Define the groups $H_k$ ($k\ge0$) by induction on $k$ in the
following way. Let $H_0=1$, $H_{k+1}=H_k\wr\la a_{k+1}\ra$ for $k\ge0$,
where all groups $\la a_k\ra$ are infinite cyclic. According to
Corollary~\ref{ZwrZetc}, all of them are embeddable into $F$. Let us
fix an integer $d\ge2$ and consider the group $H_d\times H_d$, which
is also embeddable into $F$. For any integers $k$, $n$ define an
element $g_k(n)$ as a left-normalized commutator
$$
g_k(n)=[a_1^n,a_2^n,\ldots,a_k^n],
$$
defined by induction on $k$: $g_1(n)=a_1^n$, $g_{k+1}(n)=[g_k(n),a_{k+1}^n]$
for $k\ge1$.

The element $g_k(1)$ will be denoted by $g_k$. For $R_d$ we take the
set of a single element $g_d=[a_1,a_2,\ldots,a_d]$.
\vspace{0.5ex}

The elements in $H_d$ of the form
$$
a_i(t_1,\ldots,t_r)=a_i^{a_{i+1}^{t_1}\ldots\,a_{i+r}^{t_r}},
$$
where $1\le i\le d$, $0\le r\le d-i$, $t_1,\ldots,t_r\in\zz$, will be
called {\em basic}. If $r=0$ then we just have $a_i$. Obviously,
$a_i=a_i(0)=a_i(0,0)=\cdots\,$ and so on. In general, we can add zeroes
on the right to the sequence $t_1$, \dots, $t_r$ in such a way that
the total number of arguments in brackets after $a_i$ do not exceed
$d-i$.

Consider two elements $w_i=a_i(s_{i+1},\ldots,s_d)$ and
$w_j=a_j(t_{j+1},\ldots,t_d)$, where $1\le i\le j\le d$, and $s_k$
($i<k\le d$), $t_k$ ($j<k\le d$) are integers. (It is easy to see
that each pair of basic elements can be presented in this form.)
It is clear that if the sequence $t_{j+1}$, \dots, $t_d$ is not an end
of the sequence $s_{i+1}$, \dots, $s_d$, then the elements $w_i$ and $w_j$
commute. Indeed, in this case one can choose the biggest $k$ such that
$s_k\ne t_k$. Conjugation by the inverse element to
$a_{k+1}^{s_{k+1}}\ldots a_d^{s_d}=a_{k+1}^{t_{k+1}}\ldots a_d^{t_d}$
takes elements $w_i$, $w_j$ into the elements
$w_i'=a_i(s_{i+1}\ldots s_k)\in H_{k-1}^{a_k^{s_k}}$,
$w_j'=a_j(t_{j+1}\ldots t_k)\in H_{k-1}^{a_k^{t_k}}$, respectively.
But it is clear from the elementary properties of wreath products
that the subgroups $G^{z^s}$ and $G^{z^t}$ of $G\wr\la z\ra$, where
$z$ generates $\zz$, commute elementwise for any $s\ne t$. Now, if the
sequence $t_{j+1}$, \dots, $t_d$ is the end of $s_{i+1}$, \dots, $s_d$,
that is, $s_k=t_k$ for $j<k\le d$, then elements $w_i$ and $w_j$
coincide in the case $i=j$; in the case $i<j$ one can write them as
$w_i=a_i(s_{i+1},\ldots,s_j)^v$, $w_j=a_j^v$, where
$v=a_{j+1}^{s_{j+1}}\ldots a_d^{s_d}$. Then for each $\ell\in\zz$ one has
equalities
\begin{eqnarray}
w_i^{w_j^\ell}&=&\left(a_i(s_{i+1},\ldots,s_j)^v\right)^{a_j^{\ell v}}=
\left(a_i(s_{i+1},\ldots,s_j)^{a_j^\ell}\right)^v=
a_i^{a_{i+1}^{s_{i+1}}\ldots\,a_j^{s_j+\ell}v}\nonumber\\
&=&a_i(s_{i+1},\ldots,s_j+\ell,s_{j+1},\ldots,s_d).
\end{eqnarray}
So we have a rule how to conjugate one basic element by another
basic element.
\vspace{0.5ex}

Let $1\le k\le d$. Consider the normal closure $M_k$ of the element
$a_1$ in $H_k$. It follows from the above that $M_k$ is an abelian group
freely generated by the set of elements
$$
a_1(s_2,\ldots,s_k)=a_1^{a_2^{s_2}\ldots\,a_k^{s_k}},
$$
where $s_2,\ldots,s_k\in\zz$. It is possible to define a homomorphism
$\phi_k\colon M_k\to\zz$ from $M_k$ into the {\bf additive} group
$\zz$ as follows: $\phi_k(a_1(s_2,\ldots,s_k))=s_2\ldots s_k$. From this
definition we have that for any $k>1$, $h\in M_{k-1}$ and for any
$\ell\in\zz$ the following equality holds:
$$
\phi_k(h^{a_k^\ell})=\ell\phi_{k-1}(h).
$$
In particular, $\phi_k$ equals zero on $M_{k-1}$.
\vspace{0.5ex}

Our aim is to establish the two facts.
\vspace{0.5ex}

1) For any $h\in H_k$ the equality $\phi_k(g_k^h)=1$ holds.
\vspace{0.5ex}

Note that $g_k$ obviously belongs to $M_k$ so we can apply $\phi_k$ to
any element conjugated to $g_k$.
\vspace{0.5ex}

2) If $1\le k\le d$, then the element $g_k(n)$ belongs to the
normal closure of the element $g_k$ and $\phi_k(g_k(n))=n^k$.
\vspace{0.5ex}

First we shall deduce the conclusion of our Theorem from these facts.
The elements $a_1$, \dots, $a_d$ generate the subgroup $H_d$. The
length of $g_d(n)$ with respect to these generators does not exceed
$Dn$, where $D=3\cdot2^{d-1}-2$ is a constant that does not depend on
$n$. From the above two facts it is clear that the element $g_k(n)$,
being a product of conjugates to $g_d$, cannot be presented as a product
of less than $n^d$ factors that are conjugates to $g_d$ or their
inverses. In the notation of Lemma~\ref{Phi-Dist}, this gives inequality
$\Phi(Dn)\ge n^d$. Applying this Lemma, we get $n^d\preceq\disto(n)$.
\vspace{1ex}

So let us prove the first of the above facts. We proceed by induction
on $k$. If $k=1$, then $g_1=a_1\in M_1$ and $g_1^h=a_1$ for any
$h\in H_1=\la a_1\ra$. By definition, $\phi_1(a_1)=1$. Let $k>1$, $h\in H_k$.
Then $g_k=[g_{k-1},a_k]=g_{k-1}^{-1}g_{k-1}^{a_k}\in M_k$ since
$g_{k-1}\in M_{k-1}$ by the inductive assumption. We have equalities
$$
\phi_k(g_k^h)=\phi_k\left([g_{k-1},a_k]^h\right)=
\phi_k(g_{k-1}^{-h}g_{k-1}^{a_kh})=
\phi_k(g_{k-1}^{-h})+\phi_k(g_{k-1}^{a_kh}).
$$
Since $\phi_k=0$ on $M_{k-1}$, the first summand equals zero. Further, the
elements $g_{k-1}^{a_k}\in H_{k-1}^{a_k}$ and $h\in H_{k-1}$
commute, what follows from the definition of a wreath product. Therefore,
$g_{k-1}^{a_kh}=g_{k-1}^{a_k}$. It follows from the above properties of
$\phi_k$ that for any $g\in M_{k-1}$ we have $\phi_k(g^{a_k})=\phi_{k-1}(g)$.
So the second summand equals $\phi_k(g_{k-1}^{a_k})=\phi_{k-1}(g_{k-1})=1$
because $g_{k-1}\in M_{k-1}$. As a result, $\phi_k(g_k^h)=1$, what we
had to prove.
\vspace{1ex}

Let us prove the second fact. By $N_k$ we denote the normal closure of
$g_k$ in $H_k$. Let us prove by induction on $k$ that $g_k(n)\in N_k$.
This is obvious for $k=1$ since $g_1(n)=a_1^n=g_1^n$. Let $k>1$, and let
the fact is true for all values of the parameter less than $k$. Since
$g_k=[g_{k-1},a_k]$, we have equality $g_{k-1}^{a_k}=g_{k-1}$ modulo
$N_k$. In the group $H_k$, any element in $H_{k-1}^{a_k}$ commutes with any
element in $H_{k-1}$. Therefore, $g_{k-1}$ centralizes $H_{k-1}$ modulo
$N_k$. Then, modulo $N_k$, any element in the normal closure of $g_{k-1}$
is some power of $g_{k-1}$. In particular, this is true for the element
$g_{k-1}(n)$ by the inductive assumption. Since $g_{k-1}$ and $a_k$
commute modulo $N_k$, we deduce that $g_k(n)=[g_{k-1}(n),a_k^n]$ equals
$1$ in the quotient group $H_k/N_k$, that is,  $g_k(n)\in N_k$.
\vspace{0.5ex}

Now we prove that $\phi_k(g_k(n))=n^k$ for $1\le k\le d$ by induction
on $k$. For $k=1$ we get $\phi_1(g_1(n))=\phi_1(a_1^n)=n\phi_1(a_1)=n$.
Let $k>1$; suppose that $\phi_{k-1}(g_{k-1}(n))=n^{k-1}$. Then
$\phi_k(g_k(n))=\phi_k([g_{k-1}(n),a_k^n])=\phi_k(g_{k-1}(n)\iv)+
\phi_k(g_{k-1}(n)^{a_k^n})=n\phi_{k-1}(g_{k-1}(n))=n\cdot n^{k-1}=n^k$
(we have used the properties of $\phi_k$, the fact that
$g_{k-1}(n)\in M_{k-1}$ and the inductive assumption).
\vspace{0.5ex}

The Theorem is proved.
\vspace{2ex}

It is an interesting question what else functions may be distortion
functions of finitely generated subgroups of $F$. In particular, it is
very interesting if such a distortion function may not have a recursive
upper bound. Let us give an equivalent form of this problem.

\begin{prob}
\label{MembF}
Does R.\,Thompson's group $F$ have a finitely generated subgroup
with unsolvable membership problem?
\end{prob}


\begin{thebibliography}{99}

\newcommand{\bi}{\bibitem}
\newcommand{\nb}{\newblock}

\bi{Bogl87}
W.~A.~Bogley.
\nb Retractive maps and local collapsibility.
\nb PhD thesis, Univ. of Oregon, 1987.

\bi{Brin}
M.~Brin.
\nb The ubiquity of Thompson's group $F$ in groups of piecewise linear
ho\-meo\-mor\-phisms of the unit interval.
\nb J. London Math. Soc. (to appear).

\bi{BrSq85}
M.~G.~Brin and C.~C.~Squier.
\nb Groups of piecewise linear homeo\-mor\-phisms of the real line.
\nb Invent. Math., 79 (1985), 485--498.

\bi{Bro87}
K.~S.~Brown.
\nb Finiteness properties of groups.
\nb J. of Pure and Applied Algebra, 44:45--75, 1987.

\bi{Bur}
J.~Burillo.
\nb Quasi-isometrically embedded subgroups of Thompson's group $F$.
\nb J. Algebra 212 (1999) No.\,1, pp. 65--78.

\bi{CFP}
J.~W.~Cannon, W.~J.~Floyd and W.~R.~Parry.
\nb Introductorary notes on Richard Thompson's groups.
\nb L'Enseignement Math\'ematique (2) 42 (1996), 215--256.

\bi{Cheb}
A.~A.~Chebotar.
\nb Subgroups of one-relator groups that do not contain free subgroups
of rank $2$.
\nb Algebra i Logika 10: 5 (1971), 570--586 (Russian).

\bi{Grom93}
M.~Gromov.
\nb Asymptotic invariants of infinite groups.
\nb London Math. Soc. Lect. Notes Ser. 1993, V. 182.
\nb Geometric group theory, V.2, P. 1-125.

\bi{GrSeg}
F.Grunewald and D. Segal.
\nb Some general algorithms. 1. Arithmetic groups.
\nb Ann.~Math., 112 (1980), 531--583.

\bi{Gu97}
V.~S.~Guba.
\nb On the relationship between the problems of equality and divisibility of
words for semigroup with a single defining relation.
\nb Izvesiya RAN: Ser. Mat. 61: 6 (1997) 27--58 (Russian).
\nb English transl. in: Izvestiya Mathematics 61: 6 (1997) pp. 1137--1169.

\bi{Gu98}
V.~Guba.
\nb Polynomial upper bounds for the Dehn function of R.\,Thompson's
group $F$.
\nb Journal of Group Theory 1 (1998), 203--211.

\bi{GuSa97}
V.~S.~Guba, M.~V.~Sapir.
\nb Diagram groups.
\nb Memoirs of the Amer. Math. Soc. 130, N~620, 1997, 1--117.

\bi{GuSa97a}
V.~S.~Guba, M.~V.~Sapir.
\nb The Dehn function and a regular set of normal forms for
R.~Thompson's group $F$.
\nb J.~Austral.~Math.~Soc. (Ser. A) 62 (1997), 315--328.

\bi{Haef}
A.~Haefliger.
\nb Complexes of groups and orbihedra.
\nb In: {\em Group Theory From a Geometrical Viewpoint}, (E.\,Ghys,
A.\,Haefliger, A.\,Verjovsky ed.), ICTP, Trieste, Italy, World Scientific,
1991, 504--540

\bi{Higg92}
P.~M.~Higgins.
\nb Techniques of Semigroup Theory.
\nb Oxford University Press, New York, 1992.

\bi{Kash70}
E.~V.~Kashintsev.
\nb Graphs and the word problem for finitely presented semigroups.
\nb Uch. Zap. Tul. Ped. Inst. 2 (1970), 290--302 (in Russian).

\bi{KilDiss}
V.~Kilibarda.
\nb On the algebra of semigroup diagrams.
\nb PhD Thesis, Univ.~of~Nebraska--Lincoln, 1994.

\bi{Kil}
V.~Kilibarda.
\newblock On the algebra of semigroup diagrams.
\newblock Int. J. of Alg. and Comput. 7 (1997), 313--338.

\bi{Kuz95}
Yu.~V.~Kuzmin.
\nb On one way for constructing $C$-groups.
\nb Izvestiya RAN: Ser. Mat. 59:4 (1995), 105--124 (Russian).

\bi{Kuz96}
Yu.~V.~Kuzmin.
\nb Groups of knoted compact surfaces and central extensions.
\nb Matem. Sb. 187: 2 (1996), 81--102 (Russian).

\bi{Loth83}
M.~Lothaire.
\nb Combinatorics on Words.
\nb Volume~17 of {\it Encyclopedia of Mathematics and its Applications},
Addison-Wesley, 1983.

\bi{LS80rus}
R.~Lyndon, P.~Schupp.
\nb Combinatorial group theory.
\nb Springer--Verlag, 1977.

\bi{Mikh58}
K.~A.~Mikhailova.
\nb The membership problem for direct products of groups.
\nb DAN SSSR 119 (1958), 1103--1105 (Russian).

\bi{Olsh97}
A.~Yu.~Ol'shanskii.
\nb Distortion functions for subgroups. In: {\em Geometric Group Theory
Down Under}, Proc. of a special year in geometric group theory, Canberra,
Australia, 1996. Eds. J. Cossey et al, Walter de Gruyter (1999) pp. 281--291.

\bi{Olsh98}
A.~Yu.~Ol'shanskii.
\nb On the subgroup distortion in finitely presented groups.
\nb Matem. Sb. 188: 11 (1997), 51--98 (Russian).
\nb English transl. in: Sbornik: Mathematics 188: 11 (1997),
pp. 1617--1664.

\bi{Pr95a}
S.~J.~Pride.
\nb Geometric methods in combinatorial group theory.
\nb In: J.~Fountain ed., {\em Semigroups, Formal Languages and Groups}.
\nb Kluwer Acad. Publ., Dordrecht, 1995, pp. 215--232.

\bi{Pr95b}
S.~J.~Pride.
\nb Low-dimensional homotopy theory for monoids.
\nb Int. J. of Alg. and Comput., v.5, \No 6, 1995, pp. 631--649.

\bi{Rem}
J.~H.~Remmers.
\nb On the geometry of semigroup presentations.
\nb Advances in Math., 36(3):283--296, 1980.

\bi{Sark}
R.~A.~Sarkisyan.
\nb The conjugacy problem for sequences of integer matrices.
\nb Matem. Zametki 25: 6 (1979), 811--824 (Russian).

\bi{Se80}
J.-P.~Serre.
\nb Trees.
\nb Springer--Verlag, 1980.

\bi{Sq94}
C.~C.~Squier.
\nb A finiteness condition for rewriting systems, revision by F.~Otto
and Y.~Kobayashi.
\nb Theoret. Comput. Sci., 131:271--294, 1994.

\bi{Sta87}
J.~R.~Stallings.
\nb A graph-theoretic lemma and group embeddings.
\nb Ann.Math., 111 (1987), 145--155.

\bi{Still80}
J.~Stillwell.
\nb Classical topology and combinatorial group theory.
\nb Springer-Verlag, 1980.

\end{thebibliography}
\end{document}